\documentclass[12pt,oneside]{article}
\usepackage[cp1251]{inputenc}
\usepackage[english]{babel}
\usepackage{amsmath}
\usepackage{amsthm}
\usepackage{amssymb}
\usepackage{amsfonts}
\usepackage{latexsym}
\usepackage{calc}
\usepackage{indentfirst}
\usepackage{color}
\usepackage[breaklinks]{hyperref}
\newtheorem*{proposition*}{Proposition}
\newtheorem{theorem}{Theorem}
\newtheorem*{theorem*}{Theorem}
\newtheorem{remark}{Remark}

\newtheorem*{definition*}{Definition}
\newtheorem{lemma}{Lemma}
\newtheorem*{lemma*}{Lemma}
\newtheorem{corollary}{Corollary}
\newtheorem*{corollary*}{Corollary}

\newcommand{\sep}{/\kern-2pt/ }

\allowdisplaybreaks

\begin{document}
% \newtheorem{theorem}{Theorem}
% \newtheorem{lemma}[thm]{Lemma}
% \newdefinition{remark}{Remark}
% \newproof{proof}{Proof}
% \newproof{pot}{Proof of Theorem \ref{thm2}}
%  \newtheorem{corollary}{Corollary}
% \begin{frontmatter}
\begin{center}
\Large\bf ANALYTIC IN THE UNIT BALL FUNCTIONS OF BOUNDED $L$-INDEX IN
DIRECTION
\end{center}

\begin{center}
\large\bf \MakeUppercase{BANDURA ANDRIY$^1$ and SKASKIV OLEH$^2$}
\end{center}

$^1$ Department of Higher Mathematics, Ivano-Frankivsk National Technical University of Oil and Gas, 
15 Karpatska st., Ivano-Frankivsk, 76019, Ukraine\\
E-mail: andriykopanytsia@gmail.com

$^2$ Department of Theory of Functions and Probability Theory, Ivan Franko National University of Lviv, 
1 Universytetska st., Lviv, 79000, Ukraine\\
olskask@gmail.com

\begin{abstract}
We propose a generalisation of analytic in a domain function of bounded index, which
was introduced by J. G. Krishna and S. M. Shah \cite{krishna}. In fact, analytic in the unit ball
function of bounded index by Krishna and Shah is an entire function. Our
approach allows us to explore properties of analytic in the unit ball functions.

We proved the necessary and sufficient conditions of bounded $L$-index in direction for
analytic functions. As a result, they are applied to study partial differential equations and get sufficient
conditions of bounded $L$-index in direction for analytic solutions. Finally, we estimated
growth for these functions.
\end{abstract}

\textbf{Keywords:}  analytic function, unit ball, bounded $L$-index in direction, growth estimates, partial differential equation, several complex variables
 
 \textbf{MSC}[2010] 32A10, 32A17, 35B08

%  \end{frontmatter}
%  \linenumbers

% \begin{center}
\section{Introduction}
% \end{center}

B. Lepson \cite{lepson} %\cite{lepson0,lepson}
introduced a class of
entire functions of bounded index.
He raised the problem to characterise entire functions of bounded index. An entire function $f$ is said
to be of bounded index if there exists an integer $N > 0$ that
\begin{equation}\label{deflepson}
(\forall z\in\mathbb{C})(\forall n\in\{0,1,2,\ldots\})\colon\
\frac{|f^{(n)}(z)|}{n!}\leq\max\Big\{|f(z)|,\frac{|f^{(j)}(z)|}{j!}\colon
1\leq j\leq N \Big\}.
\end{equation}
The least such integer $N$ is called the index of $f.$

Afterwards, S. Shah \cite{shah} and W. Hayman \cite{Hayman}
independently   proved that every entire function of bounded index
is a function of exponential type. 
 Namely, its growth is at most the first order and normal type. Further, W. Hayman showed that an entire function is
of bounded value distribution if and only if its derivative is of bounded index.
 An entire function $f$ is
said to be of bounded value distribution if for every $r>0$ there
exists a fixed integer $p(r) > 0$ such that the equation $f(z) =
w$ has never more than $p(r)$ roots in any disc of radius $r$ and
for any $w\in\mathbb{C}$.  The functions of bounded index have been used in the
theory value distribution and differential equations (see
bibliography in \cite{shah}).

T. Lakshminarasimhan \cite{lakshmi1974} generalised a 
 bounded index. He introduced entire functions of $L$-bounded index, where $L(r)$ is a positive continuous slowly
increasing function. D. Somasundaram and R. Thamizharasi \cite{somathami1983}-\cite{somathami1988} continued his 
investigations of entire functions of $L$-bounded index. They studied growth properties and
characterisations of these functions.

B. C. Chakraborty, Rita Chanda and Tapas Kumar Samanta \cite{Chakraborty1995}-\cite{Chakraborty2001} introduced
bounded index and $L$-bounded index for entire functions in $\mathbb{C}^n.$ They found a necessary
and sufficient condition for an entire function to be of $L$-bounded index and proved some
interesting properties.

J. Gopala Krishna and S. M. Shah \cite{krishna} studied of the existence and analytic continuation of the local solutions of partial differential equations. They introduced an analytic in a domain (a nonempty connected open set)
 $\Omega\subset \mathbb{C}^n$ $(n\in\mathbb{N})$ function of bounded index for $\alpha=(\alpha_1,
\ldots,\alpha_n)\in \mathbb{R}^n_{+}.$ Namely, let $\Omega_{+}=\{z=(z_1,\ldots, z_n)\in \Omega\colon\ z_j>0\ (j\in\{1,\ldots , n\})\},$ 
that is a subset  of all points of $\Omega$ with positive real coordinates. We say that an analytic in $\Omega$ function $F$ is a function of bounded index (Krishna-Shah bounded index or $F\in\mathcal{B}(\Omega,\alpha)$)
for $\alpha=(\alpha_1,
\ldots,\alpha_p)\in \Omega_{+}$ in domain $\Omega$ if and only if there exists $N=N(\alpha,F)=(N_1,\ldots , N_n)\in\mathbb{Z}^n_{+}$ such that inequality
$$
\alpha^m T_m(z)\leq\max\{\alpha^pT_p(z)\colon p\leq N\},
$$
is valid for all $z\in\Omega$ and for every $m\in\mathbb{Z}^{n}_{+},$ where
$\alpha^m=\alpha_1^{m_1}\cdots\alpha_n^{m_n},$\ $T_m(z)=|F^{(m)}(z)|/m!$, $F^{(m)}(z)=\frac{\partial^{\|m\|}F}{\partial
z_1^{m_1}\cdots\partial z_n^{m_n}}$ be $\|m\|$-th partial derivative of $F$, $F^{(0,\ldots, 0)}=F,$\ $m!=m_1!\cdots
m_n!,$ $\|m\|=m_1+\ldots m_n,$\ $m=(m_1,\ldots m_n)\in\mathbb{Z}^{n}_{+}.$

For entire functions in two variables, M. Salmassi \cite{salmassiart} generalised bounded index and
proved three criteria of index boundedness. Besides, he researched a system of partial
differential equations and found conditions of bounded index for entire solutions.

To consider the functions of nonexponential type A.~D. Kuzyk and M.~M. Sheremeta \cite{vidlindex}
introduced a bounded $l$-index, replacing $ \frac{ |f^{(p)}(z)|}{p!}$ on $ \frac{ |f^{(p)}(z)|}{p!l^{p}(|z|)}$ in 
\eqref{deflepson}, where $l: \mathbb{R}_{+}\to \mathbb{R}_{+}$  is
a continuous function. Besides, they proved that growth of entire function of L-bounded
index is not higher than a normal type and first order.

Afterwards, S. M. Strochyk and M. M. Sheremeta \cite{strosher} considered bounded l-index
for functions, that are analytic in a disc. Later T. O. Banakh, V. O. Kushnir and
M. M. Sheremeta generalised this term for analytic in arbitrary complex domain $G\subset\mathbb{C}$ functions
(\cite{bankush}, \cite{kusher} -- \cite{kushpro}). Yu. S. Trukhan and M. M. Sheremeta got sufficient conditions
of bounded l-index for infinite products, which are analytic in the unit disc. In particular,
they researched Blaschke product and Naftalevich-Tsuji product (\cite{shertrukh}, \cite{trukhbliashke} -- \cite{Trukh_Sher3}).

M. T. Bordulyak and M. M. Sheremeta (\cite{bagzmin} -- \cite{prostir}) defined a function of bounded
$\mathbf{L}$-index in joint variables, where $\mathbf{L}=\mathbf{L}(z)=\!(l_1(z_1),\ldots,l_n(z_n)\!),$ $l_j(z_j)$ are positive
continuous functions, $j\in\{1,\ldots n\}.$
 If
$\mathbf{L}(z)\equiv\left(\frac{1}{\alpha_1},\ldots,\frac{1}{\alpha_n}\right)$ and $\Omega=C^n$ then a Bordulyak-Sheremeta's definition
 matches with a Krishna - Shah's definition.
If $n=2$ and $\mathbf{L}(z)\equiv (1,1)$ then a
Bordulyak-Sheremeta's definition matches with a Salmassi's definition \cite{salmassiart}.

Methods for investigation of analytic functions in $\mathbb{C}^n$ 
 are divided into several groups.
 One group is based on the study of function $F$ as analytic in each variable separately.
Other methods arise in the study of slice function that is analytic functions of one variable
 $g(\tau)=F(a+b\tau),$
$\tau\in\mathbb{C}.$ This is a restriction of the analytic function $F$ to arbitrary
complex lines $\{z=a+b\tau: \tau\in\mathbb{C} \},$ $a,$ $b \in\mathbb{C}^n.$

Using the first approach, M. T. Bordulyak and M. M. Sheremeta \cite{bagzmin} proved
many properties and criteria of bounded $\mathbf{L}$-index for entire functions in $\mathbb{C}^n.$
They got sufficient conditions of bounded $\mathbf{L}$-index for entire solutions of some
systems of partial differential equations.
However, this approach did not allow to find an equivalent to
a criterion of bounded $\mathbf{L}$-index by the estimate of the logarithmic derivative outside zero
set. In particular, efforts to explore $\mathbf{L}$-index boundedness for some important classes of
entire functions (for example, infinite products with ”plane” zeros) were unsuccessful by
technical difficulties.

% Bordulyak-Sheremeta's definition is well suited to study entire functions of the form
%  $F(z)=f_{1}(z_{1})f_{2}(z_{2})\cdots f_{n}(z_{n}),$
% $F(z)=f(z_{1}+z_{2}+\cdots+z_{n})$  and etc.

For the reasons given above, there was a natural problem to consider and to explore
an entire in $\mathbb{C}^n$ function of bounded $L$-index by a second approach.

Applying this method, we proposed a new approach to introduce an entire in $\mathbb{C}^n$ function of bounded
$L$-index in direction \cite{banduramodified} -- \cite{bandexistsk}.
In contrast to the approach proposed by  M.T. Bordulyak and M. M. Sheremeta, our
definition is based on directional derivative.
It allowed to generalize more results from $\mathbb{C}$ to $\mathbb{C}^n$ and
find new assertions  because a definition contains a directional derivative and it has influence on the
$L$-index.

This success gives possibility of generalisation of bounded $L$-index in direction for
analytic in a ball functions. Besides, analytic in a domain function of bounded index 
by Krishna and Shah is an entire function. It follows from necessary condition
 of $l$-index boundedness for analytic in the unit disc function (\cite{sher},Th.3.3, p.71): $\int_0^rl(t)dt\to\infty$ as $r\to 1.$
 In this paper, we proved criteria of $L$-index boundedness in direction, which
describe a maximum modulus estimate on a larger circle by maximum modulus on a
smaller circle, an analogue of Hayman Theorem, a maximum modulus estimate on circle
by minimum modulus on circle, an estimate of logarithmic directional derivative outside
zero set and an estimate of counting function of zeros. They helped to get conditions
on partial differential equation which provide bounded $L$-index in direction for analytic
solutions. Finally, we describe the growth of analytic in $\mathbb{B}_n$ function of bounded $L$-index in
direction.

\textbf{Remark.} We investigate analytic functions in the unit ball instead the ball of arbitrary radius.

% Below we assume that $R=1.$ Thus, we investigate analytic
% in the unit ball functions of bounded $L$-index in direction.
% It is clear that this case is equivalent to the case of arbitrary ball.

% \begin{center}
\section{Main definition and properties functions of bounded $L$-index in direction} 
% \end{center}  
Let $\mathbf{b}=(b_1,\ldots,b_n)\in\mathbb{C}^n $ be a given direction, $\mathbb{B}_n=\{z\in\mathbb{C}^n: |z|< 1\},$ 
 $\overline{\mathbb{B}}_n =\{z\in\mathbb{C}^n: |z|\leq 1\},$ 
$L: \ \mathbb{B}_n\to \mathbb{R}_{+}$ be a continuous function  that for all
 $z\in\mathbb{B}_n$
 \begin{equation} \label{growthL}
L(z)>\frac{\beta|\mathbf{b}|}{1-|z|}, \ \beta=\mathrm{const}>1, \mathbf{b}\neq 0. 
\end{equation}
For a given $z\in\mathbb{B}_n$ we denote $S_{z}=\{t\in\mathbb{C}: z+t\mathbf{b}\in\mathbb{B}_n\}.$

\begin{remark} \label{kruh}
Notice that if $\eta \in[0,\beta],$\ $z\in\mathbb{B}_n,$\ $z+t_0\mathbf{b}\in \mathbb{B}_n$ and $|t-t_{0}|\leq\frac{\eta}{L(z+t_{0}\mathbf{b})}$ then $z+t\mathbf{b}\in \mathbb{B}_n.$ Indeed,  we have
$|z+t\mathbf{b}|=|z+t_0\mathbf{b}+ (t-t_0)\mathbf{b}|\leq |z+t_0\mathbf{b}|+ |(t-t_0)\mathbf{b}|\leq |z+t_0\mathbf{b}|+ \frac{\eta|\mathbf{b}|}{L(z+t_{0}\mathbf{b})}< |z+t_0\mathbf{b}|+ \frac{\beta|\mathbf{b}|}{\frac{\beta|\mathbf{b}|}{1-|z+t_0\mathbf{b}|}}=1.$
\end{remark}

Analytic in $\mathbb{B}_n$ function $F(z)$ is called a function of {\it bounded $L$-index in a direction $\mathbf{b}\in
\mathbb{C}^{n}$,} if there exists $m_{0}{}\in \mathbb{Z}_{+}$ that for
every $m\in\mathbb{Z}_{+}$ and every $z\in \mathbb{B}_n$
the following inequality is valid
\begin{equation}\label{eq1}
\frac{1}{m!L^{m}(z)}\left|\frac{\partial^{m}F(z)}{\partial
\mathbf{b}^{m}}\right|
\leq\max\left\{\frac{1}{k!L^{k}(z)}\left|\frac{\partial^{k}F(z)}{\partial
\mathbf{b}^{k}}\right|: 0\leq k \leq m_{0} \right\},
\end{equation}
where $\frac{\partial^{0}F(z)}{\partial \mathbf{b}^{0}}=F(z),
 \frac{\partial F(z)}{\partial
\mathbf{b}}=\sum\limits_{j=1}^{n}\frac{\partial F(z)}{\partial
z_{j}}{b_{j}}=\langle \mathbf{grad}\ F,$
$\overline{\mathbf{b}}\rangle, \frac{\partial^{k}F(z)}{\partial
\mathbf{b}^{k}}=\frac{\partial}{\partial
\mathbf{b}}\Big(\frac{\partial^{k-1}F(z)}{\partial
\mathbf{b}^{k-1}}\Big),$ $k \geq 2.$

The least such integer  $m_{0}=m_{0}(\mathbf{b})$
is called the {\it  $L$-index in direction $\mathbf{b}$
of the analytic function $F(z)$ and is denoted by $N_{\mathbf{b}}(F,L)=m_0.$}
If  $n=1,$ $\mathbf{b}=1,$  $L=l,$ $F=f,$ then $N(f,l)\equiv N_{1}(f,l)$ is called the $l$-index of function $f.$

In the case $n=1$ and $\mathbf{b}=1$ \ we have definition of analytic in the unit disc function of bounded $l$-index \cite{strosher}.

% \begin{center}
% \subsection{Elementary properties of $L$-index in direction and a class $Q_{\mathbf{b},\beta}(\mathbb{B}_n)$}
% \end{center}
% Now we state several lemmas that contain the basic properties
% analytic in the unit ball functions of bounded index in direction.
% Denote $\mathbb{B}_n\equiv \mathbb{D}_{1}^n,$ $\mathbb{D}^{z}\equiv \mathbb{D}_1^{z}.$

Now we state several lemmas that contain the basic properties of
analytic in the unit ball functions of bounded $L$-index in direction.
Let $l_{z}(t)=L(z+t\mathbf{b}),$ $g_z(t)=F(z+t\mathbf{b})$ for given $z\in\mathbb{C}^n.$
\begin{lemma}\label{plyash1}
If $F(z)$ is an analytic in $\mathbb{B}_n$ function of bounded $L$-index
$N_{\mathbf{b}}(F,L)$ in direction $\mathbf{b}\in\mathbb{C}^{n}$,
then for every $z^{0}\in\mathbb{B}_n$ the analytic function
$g_{z^{0}}(t),$ $t \in{ S_{z^0}}$,
 is of bounded $l_{z^{0}}$-index and
 $N(g_{z^{0}},l_{z^{0}})\leq N_{\mathbf{b}}(F,L).$
\end{lemma}
\begin{proof}
Let $z^{0}\in \mathbb{B}_n$ be a fixed point and $g(t)\equiv
g_{z^{0}}(t)$, $l(t)\equiv l_{z^{0}}(t)$. Since
%Зауважимо, що
%$g'(t)=\sum_{j=1}^{n}\frac{\partial F}{\partial
%z_{j}}\overline{b_{j}}\Bigr|_{z=z^{0}+t\mathbf{b}}=\frac{\partial
%F(z^{0}+t\mathbf{b})}{\partial \mathbf{b}},$\ а також
for every $p\in\mathbb{N}$
\begin{equation} \label{riv} g^{(p)}(t)=\frac{\partial^{p}
F(z^{0}+t\mathbf{b})}{\partial \mathbf{b}^{p}},
\end{equation}
 then by the definition of bounded $L$-index in direction
$\mathbf{b}$ for all $t\in S_{z^0}$ and for all
$p\in\mathbb{Z}_{+}$\ we obtain
\begin{gather*}
\frac{|g^{(p)}(t)|}{p!l^{p}(t)}=\!\frac{1}{p!L^{p}(z^{0}+t\mathbf{b})}
\Big|\frac{\partial^{p}F(z^{0}+t\mathbf{b})} {\partial
\mathbf{b}^{p}}\Big|\leq\max\Big\{\frac{1}{k!L^{k}(z^{0}+t\mathbf{b})}
\Big|\frac{\partial^{k}F(z^{0}+t\mathbf{b})}{\partial
\mathbf{b}^{k}}\Big|:\ \\  0\leq k\leq
N_{\mathbf{b}}(F,L)\Big\}=\max\Big\{\frac{|g^{(k)}(t)|}{k!l^{k}(t)}:
0\leq k\leq N_{\mathbf{b}}(F,L)\Big\}.\end{gather*} From here,
 $g(t)$ is a function of bounded $l$-index and $N(g,l)\leq
N_{\mathbf{b}}(F,L)$. Lemma \ref{plyash1}  is proved. \end{proof}

An equation (\ref{riv}) implies a following proposition.
\begin{lemma}\label{plyash2}
If $F(z)$ is an analytic in $\mathbb{B}_n$ function of bounded $L$-index in direction
$\mathbf{b}\in\mathbb{C}^{n}$ then $N_{\mathbf{b}}(F,L)=
\max{\left\{N(g_{z^{0}},l_{z^{0}}):z^{0}\in\mathbb{B}_n \right\}}.$
% where $g_{z^{0}}{}(t)=F(z^{0}+t\mathbf{b}),$ $l_{z^{0}}(t)=L(z^{0}+t\mathbf{b}),$  $t\in S_{z^0},$
% $N(g_{z^{0}},l_{z^{0}})$ is the $l_{z^{0}}$-index of function $g_{z^{0}}{}(t)$.
\end{lemma}
However, maximum can be calculated on the subset
 $A$ with points $z^{0}$, which has property $\{z^{0}+t\mathbf{b}: t\in  S_{z^0}, z^{0}\in{A}
\}=\mathbb{B}_n.$ So the following assertion is valid.
\begin{lemma}\label{plyash3}
 If $F(z)$ is an analytic in $\mathbb{B}_n$ function of bounded $L$-index in direction
$\mathbf{b}\in\mathbb{C}^{n}$ and $j_{0}$ is chosen with
$b_{j_{0}}\neq{0}$ then
$N_{\mathbf{b}}(F,L)=\max\{N(g_{z^{0}},l_{z^{0}}):z^{0}\in\mathbb{C}^{n},
z^{0}_{j_{0}}=0\}$ and if $\sum_{j=1}^n b_j\|\neq0$ then
$N_{\mathbf{b}}(F,L)=\max\{N(g_{z^{0}},l_{z^{0}}):z^{0}\in\mathbb{C}^{n},\sum_{j=1}^n z^{0}_j=0\}.$
\end{lemma}
\begin{proof}
We prove that for every $z\in\mathbb{B}_n$ there exist
$z^{0}\in \mathbb{C}^{n}$ and $t\in  S_{z^0}$ with 
 $z=z^{0}+t\mathbf{b}$ and $z^{0}_{j_{0}}=0$. Put $t=z_{j_{0}}/{b}_{j_{0}}$, $z^{0}_{j}=z_{j}-t{\mathbf{b}_{j}}$,
 $j\in\{1, 2, \ldots, n\}.$ Clearly, $z^{0}_{j_{0}}=0$ for this choice.

 However, a point $z^0$ may not be contained in $\mathbb{B}_n.$ But there exists $t\in\mathbb{C}$ that $z^{0}+t\mathbf{b}\in\mathbb{B}_n.$ Let $z^0\notin\mathbb{B}_n$ and $|z| = R_1<1.$ Therefore, $|z^{0}+t\mathbf{b}|= |z- \frac{z_{j_0}}{b_{j_0}}\mathbf{b}+t\mathbf{b}|=
|z+(t-\frac{z_{j_0}}{b_{j_0}})\mathbf{b}|\leq |z|+|t-\frac{z_{j_0}}{b_{j_0}}|\cdot|\mathbf{b}|\leq R_1+|t-\frac{z_{j_0}}{b_{j_0}}|\cdot|\mathbf{b}|< 1.$
Thus, $|t-\frac{z_{j_0}}{b_{j_0}}|<\frac{1-R_1}{|\mathbf{b}|}.$

In second part we prove
for every $z\in\mathbb{B}_n$ there exist $z^{0}\in \mathbb{C}^{n}$
and $t\in  S_{z^0}$ that $z=z^{0}+t\mathbf{b}$ and
$\sum_{j=1}^{n} z^{0}_{j}=0$. Put 
$t=\frac{\sum_{j=1}^{n} z_{j}}{\sum_{j=1}^n b_j}$ and
$z^{0}_{j}=z_{j}-t{\mathbf{b}_{j}}$, $1\leq j\leq n.$
Thus, the following equality is valid
$\sum_{j=1}^n z^{0}_j=\sum_{j=1}^{n}(z_{j}-t{b_{j}})=\sum_{j=1}^{n}z_{j}-\sum_{j=1}^n b_j t
=0.$

% Similarly, as in the proof of the first part of this lemma it can
%be proven that for all $z^0\notin \mathbb{B}_n$ and for all $t$ with a disc
% $|t-\frac{1}{\|\mathbf{b}\|}\sum_{j=1}^{n} z_{j}|<\frac{R-|z|}{|\mathbf{b}|}$  $z^0+t\mathbf{b}\in\mathbb{B}_n.$
Lemma \ref{plyash3} is proved.
\end{proof}

Note that for a given $z\in \mathbb{B}_n$  we can pick uniquely
$z^{0}\in \mathbb{C}^{n}$ and $t\in  S_{z^0}$ such that
$\sum_{j=1}^n z^{0}_j=0$ and $z=z^{0}+t\mathbf{b}.$

\begin{remark}
 If for some $z^0\in\mathbb{C}^n$ \ $\{z^0+t\mathbf{b}: \ t\in\mathbb{C}\}\bigcap \mathbb{B}_n=\emptyset$ then we put $N(g_{z^{0}},l_{z^{0}})=0.$
\end{remark}

\noindent Lemmas \ref{plyash1}--\ref{plyash3} imply the following proposition.

\begin{theorem} \label{te1}
An analytic in $\mathbb{B}_n$ function $F(z)$ is a function of bounded
$L$-index in direction $\mathbf{b}\in\mathbb{C}^{n}$ if and only if
 there exists number $M>0$ such, that for every $z^{0}\in
\mathbb{B}_n$ function $g_{z^{0}}(t)$ is of
bounded $l_{z^{0}}$-index with $N(g_{z^{0}},l_{z^{0}})\leq
M<+\infty,$\ as a function of one variable $t\in S_{z^0},$ and
$N_{\mathbf{b}}(F,L)=\max\{N(g_{z^{0}},l_{z^{0}}):z^{0}\in
\mathbb{B}_n\}.$
\end{theorem}
\begin{proof}
Necessity follows from Lemma \ref{plyash1}.

\noindent {\it We prove sufficiency.}

Since $N(g_{z^{0}},l_{z^{0}})\leq M$ there exists
$\max\{N(g_{z^{0}},l_{z^{0}}): z^{0}\in \mathbb{B}_n\}$.
We denote this maximum by 
$N_{\mathbf{b}}\!(F,L)=\max\{N(g_{z^{0}},l_{z^{0}}): z^{0}\in
\mathbb{B}_n\}<\infty.$ Suppose that $N_{\mathbf{b}}(F)$ is not $L$-index in  direction
$\mathbf{b}$ of function $F(z)$. So there exists
 $n^{*}>N_{\mathbf{b}}(F,L)$ and $z^{*}\in \mathbb{B}_n$
\begin{equation}
\label{eq2}
\frac{1}{n^{*}!L^{n^{*}}(z^{*})}\frac{\big|\partial^{n^{*}}F(z^{*})\big|}{\partial
\mathbf{b}^{n^{*}}}>
\max\bigg\{\frac{1}{k!L^{k}(z^{*})}\frac{\big|\partial^{k}F(z^{*})\big|}{\partial
\mathbf{b}^{k}},
  0\leq k\leq N_{\mathbf{b}}(F,L)\bigg\}.
\end{equation}
But we have $g_{z^{0}}(t)=F(z^{0}+t\mathbf{b}),$
$g_{z^{0}}^{(p)}(t)=\frac{\partial^{p}F(z^{0}+t\mathbf{b})}{\partial
\mathbf{b}^{p}}.$ We can rewrite (\ref{eq2}) as
\begin{equation*}
  \frac{\big|g_{z^{*}}^{(n^{*})}(0)\big|}{n^{*}!l_{z^{*}}^{n^{*}}(0)}>
  \max\bigg\{\frac{|g^{(k)}_{z^{*}}(0)|}{k!l_{z^{*}}^{k}(0)}: 0\leq k\leq
  N_{\mathbf{b}}(F,L)\bigg\}.
\end{equation*}
It contradicts that all $l_{z^{0}}$-indices
$N(g_{z_{0}},l_{z^{0}})$ are bounded by number $N_{\mathbf{b}}(F).$ Thus
$N_{\mathbf{b}}(F)$ is $L$-index in  direction $\mathbf{b}$ of function $F(z)$.
Theorem \ref{te1} is proved.
\end{proof}
From Lemma \ref{plyash3} the following condition is enough in Theorem \ref{te1}:
 {\it there exists $M<+\infty$ that an inequality holds 
$N(g_{z^{0}},l_{z^{0}})\leq M$ for every
 $z^{0}\in \mathbb{C}^{n}$\ with $\sum_{j=1}^{n}{}z^{0}_{j}=0.$}

Since Lemma \ref{plyash3} and \ref{te1} there is a natural {\it question:}\
 what is the least set $A$ that the following equality is valid
$N_{\mathbf{b}}(F,L)=\max\limits_{z^0\in A} N(g_{z^{0}},l_{z^{0}}).$

Below we prove propositions that give a partial answer to this question.
A solution is partial because it is unknown whether our sets are the least which satisfy the mentioned equality.

\begin{theorem}  \label{le1}
Let $\mathbf{b}\in\mathbb{C}^{n}$ be a given direction, $A_{0}$ be an arbitrary set in $\mathbb{C}^{n}$ with 
$\{z+t\mathbf{b}: t\in S_z, \ z\in A_{0}\}=\mathbb{B}_n.$
Analytic in $\mathbb{B}_n$ function $F(z)$ is of bounded
$L$-index in direction $\mathbf{b}\in\mathbb{C}^{n}$ if and only if
 there exists $M>0$ that for all $z^{0}\in A_{0}$ function $g_{z^{0}}(t)$ is of
bounded $l_{z^{0}}$-index $N(g_{z^{0}},l_{z^{0}})\leq
M<+\infty,$\ as a function of variable $t\in S_{z^0}.$ And
$N_{\mathbf{b}}(F,L)=\max\{N(g_{z^{0}},l_{z^{0}}):z^{0}\in
A_{0}\}.$
\end{theorem}

\begin{proof}
By Theorem \ref{te1}, analytic in $\mathbb{B}_n$ function $F(z)$
is of bounded $L$-index in  direction $\mathbf{b}\in\mathbb{C}^{n}$ if and only if  there exists number
$M>0$ such that for every $z^{0}\in\mathbb{B}_n$ function
$g_{z^{0}}(t)$ is of bounded
$l_{z^{0}}$-index $N(g_{z^{0}},l_{z^{0}})\leq M<+\infty,$\ as a function
of variable $t\in S_{z^0}$. But for every $z^{0}+t\mathbf{b}$ by properties
of set $A_{0}$  there exist $\widetilde{z}^{0}\in
A_{0}$ and $\widetilde{t}\in\mathbb{B}_{\widetilde{z}^{0}}$
$$z^{0}+t\mathbf{b}=\widetilde{z}^{0}+\widetilde{t}\mathbf{b}.$$
For all $p\in\mathbb{Z}_{+}$ we have
$$(g_{z_{0}}(t))^{(p)}=(g_{\widetilde{z}_{0}}(\widetilde{t}))^{(p)}.$$ But $\widetilde{t}$ is dependent of $t.$
Therefore, a condition $g_{z^{0}}(t)$ is of bounded
$l_{z^{0}}$-index for all $z^{0}\in\mathbb{B}_n$ is equivalent to a condition
 $g_{\widetilde{z^{0}}}(t)$ is of bounded $l_{\widetilde{z}^{0}}$-index for all $\widetilde{z}^{0}\in A_{0}.$
\end{proof}

\begin{remark} \label{hyperplane} An intersection of arbitrary hyperplane $H=\{z\in\mathbb{C}^n: \langle z, c\rangle=1\}$ and set 
 $\mathbb{B}_n^{\mathbf{b}}=\{z+\frac{1-\langle z,\mathbf{c}\rangle}{\langle \mathbf{b},c\rangle}\mathbf{b}\colon z\in\mathbb{B}_n\},$ where $\langle \mathbf{b}, c\rangle\neq 0,$
satisfies conditions of Theorem \ref{le1}.
\end{remark}

We prove that for every $w\in\mathbb{B}_n$ there exist
$z\in H\bigcap \mathbb{B}_n^{\mathbf{b}}$ and $t\in\mathbb{C}$ such that $w={z}+{t}\mathbf{b}.$

% Obviously, there are
% $\widetilde{z}\in\mathbb{C}^n,$ $\widetilde{t}\in\mathbb{C}$ such that $w=\widetilde{z}+\widetilde{t} \mathbf{b}.$
Choosing $z=w+\frac{1-\langle
w,c\rangle}{\langle\mathbf{b},c\rangle}\mathbf{b} \in  H\bigcap \mathbb{B}_n^{\mathbf{b}},$
$t=\displaystyle\frac{\mathstrut\langle
w,c\rangle-1}{\langle\mathbf{b},c\rangle},$
we obtain
$$z+t\mathbf{b}=w+\frac{1-\langle w,c\rangle}{\langle\mathbf{b},c\rangle}\mathbf{b}+ \frac{\langle w,c\rangle-1}{\langle\mathbf{b},c\rangle}\mathbf{b}=w. $$

\begin{theorem} \label{te2double}
Let $A$ be an everywhere dense set in $\mathbb{B}_n.$
Analytic in $\mathbb{B}_n$ function $F(z)$ is of bounded
$L$-index in  direction $\mathbf{b}\in\mathbb{C}^{n}$ if and
only if there exists number $M>0$ that for every $z^{0}\in A$
function $g_{z^{0}}(t)$ is of bounded $l_{z^{0}}$-index
$N(g_{z^{0}},l_{z^{0}})\leq M<+\infty,$\ as a function of
$t\in S_{z^0},$ and
$N_{\mathbf{b}}(F,L)=\max\{N(g_{z^{0}},l_{z^{0}}):z^{0}\in A\}.$
\end{theorem}

\begin{proof}
The necessity follows from Theorem \ref{te1} (in this theorem
same condition is satisfied for all $z^{0}\in\mathbb{B}_n,$ and we need this condition
 for all $z^{0}\in A,$ that $\overline{A}\cap\mathbb{B}_n=\mathbb{B}_n$).

Now we prove a sufficiency. Since $A$ has been everywhere
dense in $\mathbb{B}_n,$  for every $z^{0}\in\mathbb{B}_n$
there exists a sequence $(z^{m}),$ that $z^{(m)}\to z^{0}$ as 
$m\to+\infty$ and $z^{(m)}\in A$ for all $m\in\mathbb{N}.$ But
$F(z+t\mathbf{b})$ is of bounded $l_{z}$-index for all
$z\in\overline{A}\cap\mathbb{B}_n$ as a function of $t.$
Therefore,  by bounded $l_{z}$-index there exists $M>0$ that
for all $z\in A,$ $t\in\mathbb{C},$ $p\in\mathbb{Z}_{+}$
$$\frac{|g_{z}^{(p)}(t)|}{p!l^{p}(t)}\leq\max\left\{\frac{|g_{z}^{(k)}(t)|}{k!l^{k}_{z}(t)}: 0\leq k\leq M\right\}.$$

After substitution instead of $z$ a sequence $z^{(m)}\in A$ and $z^{(m)}\to z^0,$ for each $m\in\mathbb{N}$ the following inequality holds
$$\frac{|g_{z^{m}}^{(p)}(t)|}{p!l^{p}_{z^{m}}(t)}\leq\max \left\{\frac{|g^{(k)}_{z^{m}}(t)|} {k!l^{k}_{z^{m}}(t)}: 0\leq k\leq M\right\}$$
In other words, we have 
\begin{gather}
\!\frac{1}{p!L^{p}(z^{m}+t\mathbf{b})}\left|\frac{\partial^{p} F(z^{m}+t\mathbf{b})}{\partial\mathbf{b}^{p}}\right|\!\leq\! \max\left\{\frac{1}{k!L^{k}(z^{m}+t\mathbf{b})} \left|\frac{\partial^{k} F(z^{m}+t\mathbf{b})}{\partial\mathbf{b}^{k}}\right|: \right. \nonumber\\ \left. 0\leq k\leq M\right\}. \label{eq6}
\end{gather}

But $F$ is an analytic in $\mathbb{B}_n$ function and $L$ is a
positive continuous. In \eqref{eq6}  we 
calculate a limit $m\to+\infty$ $(z^{m}\to z^0).$ We have
that for all $z^0\in\mathbb{B}_n,$ $t\in S_{z^0},$
$m\in\mathbb{Z}_{+}$
\begin{gather*}
\frac{1}{p!L^p(z^0+t\mathbf{b})}\left|\frac{\partial^p F(z^0+t\mathbf{b})}{\partial\mathbf{b}^p}\right|\leq\max\left\{\frac{1}{k!L^k(z^0+t\mathbf{b})} \left|\frac{\partial^k F(z^0+t\mathbf{b})}{\partial\mathbf{b}^k}\right|:\right.\\ \left. 0\leq k\leq M\right\}.
\end{gather*}
Since this inequality 
 $F(z^0+t\mathbf{b})$ is of bounded $L(z^0+t\mathbf{b})$-index too, as a function of $t,$ for every given $z^0\in\mathbb{B}_n.$
 Applying Theorem \ref{te1} we get a needed conclusion.
Theorem \ref{te2double} is proved.
\end{proof}

Since Remark \ref{hyperplane} and Theorem \ref{te2double} the following corollary is true.
\begin{corollary} \label{corvazhl}
 Let $\mathbf{b}\in\mathbb{C}^{n}$ be a given direction, $A_{0}$ be a set in $\mathbb{C}^{n}$ and its closure is $\overline{A}_{0}=\{z\in\mathbb{C}^n: \langle z, c\rangle=1\}\bigcap \mathbb{B}_n^{\mathbf{b}},$ where $\langle
c,\mathbf{b}\rangle\neq 0,$ $\mathbb{B}_n^{\mathbf{b}}=\{z+\frac{1-\langle z,\mathbf{c}\rangle}{\langle \mathbf{b},c\rangle}\mathbf{b}\colon z\in\mathbb{B}_n\}.$ Analytic in $\mathbb{B}_n$ function $F(z)$ is of bounded $L$-index in  direction
$\mathbf{b}\in\mathbb{C}^{n}$ if and only if there exists number $M>0$ such that for all $z^{0}\in A_{0}$ function
$g_{z^{0}}(t)$ is of bounded $l_{z^{0}}$-index
$N(g_{z^{0}},l_{z^{0}})\leq M<+\infty,$\
as a function of variable $t\in S_{z^0}.$ And
$N_{\mathbf{b}}(F,L)=\max\{N(g_{z^{0}},l_{z^{0}}): z^{0}\in
A_{0}\}.$
\end{corollary}
\begin{proof}
Since Remark \ref{hyperplane} in Theorem \ref{le1} we
can take an arbitrary hyperplane $B_0=\{z\in\mathbb{{B}}^n :
\langle z,c\rangle =1\},$ where $\langle c,\mathbf{b}\rangle\neq
0.$ Let $A_0$ be an everywhere dense set in $B_0,$
 $\overline{A_0}=B_0.$
 Repeating considerations of Theorem \ref{te2double}, we obtain a needed conclusion.
%  (in this theorem
% corresponding condition is satisfied for all $z^{0}\in\mathbb{B}_n,$ and there require this condition
%  for all $z^{0}\in A,$ that $\overline{A}=\mathbb{B}_n$).

 Indeed,  the necessity follows from Theorem \ref{plyash1} (in this theorem
same condition is satisfied for all $z^0\in\mathbb{C}^n,$
and we need this condition for all $z^0\in A_0,$ that
$\overline{A}_0 {\cap\mathbb{B}_n}=\{z\in\mathbb{B}_n:
\langle z,c\rangle =1\}$).

To prove the sufficiency, we use a density of the set $A_0.$ Obviously, for every $z^0\in B_0$ there exists a sequence
$z^{(m)}\to z^0$ and $z^{(m)}\in A_0.$ But $g_z(t)$ is of bounded $l_z$-index for all $z\in A_0$ as a function of $t.$
Since conditions of Corollary \ref{corvazhl}, for some $M>0$ and for all $z\in A_0,$ $t\in\mathbb{C},$ $p\in\mathbb{Z}_{+}$
the following inequality holds 
$$\frac{g_z^{(p)}(t)}{p!l_z^p(t)} \leq\max\left\{\frac{|g_z^{(k)}(t)|}{k!l_z^k(t)}: 0\leq k\leq M \right\}.$$

Substituting an arbitrary sequence $z^{(m)}\in A,$ $z^{(m)}\to z^0$ instead of $z\in A^0,$ we have 
$$\frac{|g_{z^{(m)}}^{(p)}(t)|}{p!l_{z^{(m)}}^{p}(t)}\leq \max\left\{\frac{|g_{z^{(m)}}^{(k)}(t)|}{k!l_{z^{(m)}}^k(t): 0\leq k\leq M} \right\}, $$
i.e.
\begin{gather*}
\!\frac{1}{L^p(z^{(m)}+t\mathbf{b})}\left|\frac{\partial^p
F(z^{(m)}\!+\!t\mathbf{b})}{\partial\mathbf{b}^p}\right| \! \leq 
\!\max\left\{\frac{1}{k!L^k(z^{(m)}\!+\!t\mathbf{b})}\left|\frac{\partial^k
F(z^{(m)}\!+\!t\mathbf{b})}{\partial\mathbf{b}^k}\right|\!: 0\leq k\leq M\right\}.
\end{gather*}
However, $F$ is an analytic in $\mathbb{B}_n$ function, $L$ is a
positive continuous. So we calculate a limit as $m\to+\infty$ $(z^{m}\to z).$
 For all $z^0\in B_0,$ $t\in S_{z^0},$ $m\in\mathbb{Z}_{+}$ we have 
\begin{gather*}
\frac{1}{L^p(z^0+t\mathbf{b})}\left|\frac{\partial^p
F(z^0+t\mathbf{b})}{\partial\mathbf{b}^p}\right| \leq  
\max\left\{\frac{1}{k!L^k(z^0+t\mathbf{b})}\left|\frac{\partial^k
F(z^0+t\mathbf{b})}{\partial\mathbf{b}^k}\right|: \right.\ \left.
0\leq k\leq M\right\}.
\end{gather*}
Therefore, $F(z^0+t\mathbf{b})$ is of bounded $L(z^0+t\mathbf{b})$-index as a function of $t$ at each $z^0\in B^n.$
By Theorem \ref{te2double} and Remark \ref{hyperplane} $F$ is of bounded $L$-index in  direction $\mathbf{b}.$
\end{proof}

\begin{remark}
      Let $H=\{z\in\mathbb{C}^n: \langle z,c\rangle=1\}.$
  The condition $\langle c,\mathbf{b} \rangle\neq 0$ is essential. If $\langle c,\mathbf{b} \rangle= 0$ then for all $z^0\in H$ and
  for all $t\in\mathbb{C}$ the point $z^0+t\mathbf{b}\in H$ because $\langle z^0+t\mathbf{b}, c\rangle=\langle z^0,c\rangle+t\langle \mathbf{b},c\rangle=1.$ Thus, this line $z^0+t\mathbf{b}$ does not describe points outside a hyperplane $H.$
  \end{remark}
 We consider $F(z_1,z_2)=\exp(-z_1^2+z_2^2),$ $\mathbf{b}=(1,1),$ $c=(-1,1).$ On a hyperplane $-z_1+z_2=1$ 
 function $F(z_1,z_2)$ takes a look
 \begin{gather*}
  F(z^0+t\mathbf{b})=F(z^0_1+t,z^0_2+t)=\exp(-(z_1^0+t)^2+(1+z_1^0+t)^2)=\\
  =\exp(1+2z_1^0+2t).
 \end{gather*}
Using definition of $l$-index boundedness and evaluating corresponding derivatives it is easy to prove that  $\exp(1+2z_1^0+2t)$ 
is of bounded index with $l(t)=1$ and $N(g,l)=4.$ 

Thus, $F$ is of unbounded index in  direction $\mathbf{b}.$
On the contrary, we assume $N_{\mathbf{b}}(F)=m$ 
and calculate directional derivatives
%  $$
%  \frac{\partial F}{\partial \mathbf{b}}=2(-z_1+z_2)\exp(-z_1^2+z_2^2),
%  $$
%  $$
%  \frac{\partial^2F}{\partial\mathbf{b}^2}=2^2(-z_1+z_2)^2\exp(-z_1+z_2),
%  $$
 $$
 \frac{\partial^pF}{\partial\mathbf{b}^p}=2^p(-z_1+z_2)^p\exp(-z_1+z_2), \ p\in\mathbb{N}.
 $$
By definition of bounded index,  an inequality holds  $\forall$ $p\in\mathbb{N}$ $\forall z\in\mathbb{C}^n$
 \begin{equation} \label{auxexample}
 2^p|-z_1+z_2|^p|\exp(-z_1+z_2)|\leq \max_{0\leq k\leq m}2^k|-z_1+z_2|^k|\exp(-z_1+z_2)|.
 \end{equation}
  Let $p>m$  and $|-z_1+z_2|=2.$ Dividing equation \eqref{auxexample} by $2^p|\exp(-z_1+z_2)|,$ 
  we get $2^{2p}\leq 2^{2m}.$ It is impossible. Therefore, $F(z)$ is of unbounded index in  direction $\mathbf{b}.$
 
 Using calculated derivatives it can be proven that function $F(z_1,z_2)$ is of bounded $L$-index in  direction $\mathbf{b}$ with $L(z_1,z_2)=2|-z_1+z_2|+1$ and $N_{\mathbf{b}}(F,L)=0.$

  Now we consider another function 
  $$F(z)=(1+\langle z,d \rangle)\prod_{j=1}^{\infty}(1+\langle z,c\rangle \cdot 2^{-j})^j, \ c\neq d.$$ 
The multiplicity of zeros for function $F(z)$ increases to infinity. Below in this paper, we will state Theorem \ref{te6}. By that theorem,  unbounded multiplicity of zeros means that 
   $F(z)$ is of unbounded $L$-index in any direction $\mathbf{b}$ ($\langle \mathbf{b},c\rangle\neq 0$) and for any positive continuous function $L.$ 
   
   We select $\mathbf{b}\in\mathbb{C}^n$ that $\langle \mathbf{b},d\rangle=0.$ Let $H=\{z\in\mathbb{C}^n: \langle z,d\rangle=-1\}.$
  But for $z^0\in H$  we have
   \begin{gather*}
  F(z^0+t\mathbf{b})=(1+\langle z^0,d \rangle+t\langle \mathbf{b},d\rangle)\prod\limits_{j=1}^{\infty}(1+\langle z^0,c\rangle 2^{-j}+t\langle \mathbf{b},c\rangle 2^{-j})^j\equiv 0.
 \end{gather*}
  Thus, $ F(z^0+t\mathbf{b})$ is of bounded index as a function of variable $t.$

  \begin{theorem}
Let $(r_p)$ be a positive sequence such that $r_p\to 1$ as $p\to\infty,$
   $D_p=\{z\in\mathbb{C}^n: \ |z|=r_p\},$ $A_p$ be an everywhere dense set in $D_p$ (i.e. $\overline{A}_p=D_p$) and  $A=\bigcup\limits_{p=1}^{\infty}A_p.$
   Analytic in $\mathbb{B}_n$ function $F(z)$ is of bounded $L$-index in  direction
$\mathbf{b}\in\mathbb{C}^{n}$ if and only if there exists number $M>0$ that for all $z^{0}\in A$ function
$g_{z^{0}}(t)=F(z^{0}+t\mathbf{b})$ is of bounded $l_{z^{0}}$-index
$N(g_{z^{0}},l_{z^{0}})\leq M<+\infty,$\
as a function of variable $t\in S_{z^0},$ where $l_{z^{0}}(t)\equiv
L(z^{0}+t\mathbf{b}).$\ And
$N_{\mathbf{b}}(F,L)=\max\{N(g_{z^{0}},l_{z^{0}}):z^{0}\in
A\}.$
 \end{theorem}
 \begin{proof}
 Theorem \ref{te1} implies the necessity of this theorem.

  \textit{Sufficiency.} It is easy to prove $\{z+t\mathbf{b}: t\in S_z, \ z\in A\}=\mathbb{B}_n.$ Further, we repeat considerations with proof of sufficiency in Theorem \ref{te2double} and obtain a needed conclusion.
 \end{proof}

%  \begin{center}
\section{Auxiliary class $Q^n_{\mathbf{b}}$} 
% \end{center}  
 
%  {\bf{$3^{0}.$ Class $Q^n_{\mathbf{b}}.$}
The positivity and continuity of function $L$ and condition \eqref{growthL} are not enough to explore the behaviour of entire function of bounded  $L$-index in direction. Below we impose the extra condition that function $L$ does not vary as soon.
 
 For $\eta \in[0,\beta],$\ $z\in\mathbb{B}_n,$\
 $t_0\in S_{z}$ such that $z+t_0\mathbf{b}\in \mathbb{B}_n$ \
 we define
\[
\lambda^{\mathbf{b}}_{1}(z,t_{0},\eta,L)=\inf\left\{\frac{L(z+t\mathbf{b})}
{L(z+t_{0}\mathbf{b})}:
|t-t_{0}|\leq\frac{\eta}{L(z+t_{0}\mathbf{b})}\right\},  \]
$\lambda^{\mathbf{b}}_{1}(z,\eta,L)=\inf\{\lambda^{\mathbf{b}}_{1}(z,t_{0},\eta,L):
t_0\in  S_z\},$
$\lambda^{\mathbf{b}}_{1}(\eta,L)=\inf\{\lambda^{\mathbf{b}}_{1}(z,\eta,L):
z\in\mathbb{B}_n\},$ and 
\[
\lambda^{\mathbf{b}}_{2}(z,t_{0},\eta,L)=\sup\left\{\frac{L(z+t\mathbf{b})}
{L(z+t_{0}\mathbf{b})}:
|t-t_{0}|\leq\frac{\eta}{L(z+t_{0}\mathbf{b})}\right\},  \]
$\lambda^{\mathbf{b}}_{2}(z,\eta,L)=\sup\{\lambda^{\mathbf{b}}_{2}(z,t_{0},\eta,L):
t_0\in  S_z\},$
$\lambda^{\mathbf{b}}_{2}(\eta,L)=\sup\{\lambda^{\mathbf{b}}_{2}(z,\eta,L):
z\in\mathbb{B}_n\}.$

If it will not cause misunderstandings, then $\lambda^{\mathbf{b}}_{1}(z,t_{0},\eta)\equiv \lambda^{\mathbf{b}}_{1}(z,t_{0},\eta,L),$
$\lambda^{\mathbf{b}}_{2}(z,$ $t_0,$ $\eta)$ $\equiv$ $\lambda^{\mathbf{b}}_{2}(z,t_0,\eta,L),$ $\lambda^{\mathbf{b}}_{1}(z,\eta)\equiv \lambda^{\mathbf{b}}_{1}(z,\eta,L),$
$\lambda^{\mathbf{b}}_{2}(z,\eta)\equiv\lambda^{\mathbf{b}}_{2}(z,\eta,L),$ $\lambda^{\mathbf{b}}_{1}(\eta)\equiv \lambda^{\mathbf{b}}_{1}(\eta,L),$
$\lambda^{\mathbf{b}}_{2}(z,\eta)\equiv\lambda^{\mathbf{b}}_{2}(\eta,L).$

By ${Q}_{\mathbf{b},\beta}(\mathbb{B}_n)$ we denote the class of all functions $L$ for which the following condition holds for any
$\eta\in[0,\beta]$ \ \ 
$0<\lambda^{\mathbf{b}}_{1}(\eta)\leq\lambda^{\mathbf{b}}_{2}(\eta)<+\infty.$
 Let  $\mathbb{D}\equiv\mathbb{B}^1,$ $Q_{\beta}(\mathbb{D})\equiv{Q}_{1,\beta}(\mathbb{D}).$ 
 
 The following lemma suggests possible approach to compose function with $Q^n_{\mathbf{b}}.$
 
 \begin{lemma}
 Let $L: \overline{\mathbb{B}}_n \to \mathbb{R}_+$ be a continuous function, $m=\min\{L(z): z\in\overline{\mathbb{B}}_n \}.$ 
 Then $\widetilde{L}(z)=\frac{\beta|\mathbf{b}|}{m}\cdot \frac{L(z)}{(1-|z|)^{\alpha}}\in Q^n_{\mathbf{b}}(\mathbb{B}_n)$ for every $\mathbf{b}\in\mathbb{C}^n\setminus\{0\},$ $\alpha\geq 1.$
 \end{lemma}
\begin{proof}
Using definition $Q^n_{\mathbf{b}}$ we have $\forall z\in\mathbb{B}_n$ $\forall t_0\in  S_z$
 \begin{gather*}
 \lambda_{1}^{\mathbf{b}}(z,t_{0},\eta,\widetilde{L})= \\ =
 \inf \left\{\frac{L(z+t\mathbf{b})}{(1-|z+t\mathbf{b}|)^{\alpha}}\cdot \frac{(1-|z+t_0\mathbf{b}|)^{\alpha}}{L(z+t_0\mathbf{b})}:
 |t-t_0|\leq \frac{\eta m(1-|z+t_0\mathbf{b}|)^{\alpha}}{\beta|b|L(z+t_0\mathbf{b})}\right\}\geq \\ 
 \geq \inf\left\{\frac{L(z+t\mathbf{b})}{L(z+t_0\mathbf{b})}: 
 |t-t_0|\leq \frac{\eta m(1-|z+t_0\mathbf{b}|)^{\alpha}}{\beta|b|L(z+t_0\mathbf{b})}
 \right\} \times \\  \inf \left\{ \left(\frac{1-|z+t_0\mathbf{b}|}{1-|z+t\mathbf{b}|}\right)^{\alpha}: 
 |t-t_0|\leq \frac{\eta m(1-|z+t_0\mathbf{b}|)^{\alpha}}{\beta|b|L(z+t_0\mathbf{b})} 
 \right\}
 \end{gather*}
Since Remark \ref{kruh} the first infimum
is not less than some constant $K>0$ which is independent from $z$ and $t_0.$ Besides, we have 
$\forall z\in\mathbb{B}_n$ and $\forall t\in S_z$ $\frac{m}{L(z+t_0\mathbf{b})}\leq 1.$
Thus, for the second infimum the following estimates are valid
\begin{gather*}
 \inf \left\{ \left(\frac{1-|z+t_0\mathbf{b}|}{1-|z+t\mathbf{b}|}\right)^{\alpha}: 
 |t-t_0|\leq \frac{\eta m(1-|z+t_0\mathbf{b}|)^{\alpha}}{\beta|b|L(z+t_0\mathbf{b})} 
 \right\} \geq\\ \geq \inf \left\{ \left(\frac{1-|z+t_0\mathbf{b}|}{1-|z+t\mathbf{b}|}\right)^{\alpha}: 
 |t-t_0|\leq \frac{\eta (1-|z+t_0\mathbf{b}|)^{\alpha}}{\beta|b|} 
 \right\} 
 =  \left(\frac{1-|z+t_0\mathbf{b}|}{1-|z+t^*\mathbf{b}|}\right)^{\alpha}.
\end{gather*}
where $|t^*-t_0|\leq\frac{\eta (1-|z+t_0\mathbf{b}|)}{\beta|b|}.$
Now we find a lower estimate for this fraction
\begin{gather*}
 \frac{1-|z+t_0\mathbf{b}|}{1-|z+t^*\mathbf{b}|}\geq 
 \frac{1-|z+t_0\mathbf{b}|}{1-||z+t_0\mathbf{b}|-|(t^*-t_0)\mathbf{b}||}\geq
 \frac{1-|z+t_0\mathbf{b}|}{1-||z+t_0\mathbf{b}|-\frac{\eta (1-|z+t_0\mathbf{b}|)}{\beta}|}
 \end{gather*}
 Denoting $u=|z+t_0\mathbf{b}|\in[0;1),$ $\gamma=\frac{\eta}{\beta}\in[0,1],$ we consider a function of one real variable
 $s(u)=\frac{1-u}{1-|u-\alpha(1-u)|}=\frac{1-u}{1-|(1+\gamma)u-\gamma|}.$ For $u\in[0,\frac{\gamma}{\gamma+1}]$ the function $s(u)$ 
 strictly decreases and for $t\in[\frac{\gamma}{1+\gamma};1)$ the function $s(u)\equiv \frac{1}{1+\gamma}.$ 
 In fact, we proved that 
  $$\lambda_{1}^{\mathbf{b}}(z,t_{0},\eta,\widetilde{L}) \geq K\cdot \frac{1}{1+\frac{\eta}{\beta}}>0.$$
 Hence, we have $\lambda_1^{\mathbf{b}}(\eta,\widetilde{L})>0.$ By analogy it can be proved that 
%  $\lambda_{1}^{\mathbf{b}}(z,t_{0},\eta,\widetilde{L})\leq K^* \frac{1}{1-\frac{\eta}{\beta}}$ and 
 $\lambda_2^{\mathbf{b}}(\eta,\widetilde{L})<\infty.$
\end{proof}

  We often use the following properties $Q_{\mathbf{b},\beta}(\mathbb{B}_n).$

\begin{lemma}\label{lemma1}\begin{enumerate}
\item If $L\in Q_{\mathbf{b},\beta}(\mathbb{B}_n)$ then
 for every
$\theta\in\mathbb{C}\backslash\{0\}$
$L\in Q_{\theta\mathbf{b},\beta/|\theta|}(\mathbb{B}_n)$ and $|\theta| L\in
Q_{\theta\mathbf{b},\beta}(\mathbb{B}_n)$
\item If  $L\in
Q_{\mathbf{b}_{1},\beta}(\mathbb{B}_n) \bigcap Q_{\mathbf{b}_{2},\beta}(\mathbb{B}_n)$ and for all $z\in\mathbb{B}_n$
$L(z)> \frac{\beta \max\{|\mathbf{b}_1|,|\mathbf{b}_2|,|\mathbf{b}_1+\mathbf{b}_2|\}}{1-|z|}$
then \\ $\min\{\lambda_2^{\mathbf{b}_1}(\beta,L),\lambda_2^{\mathbf{b}_2}(\beta,L)\}L\in
Q_{\mathbf{b}_{1}+\mathbf{b}_{2},\beta}(\mathbb{B}_n).$
\end{enumerate}
\end{lemma}

\begin{proof}

\textbf{1.}
First, we prove that $(\forall
\theta\in\mathbb{C}\backslash\{0\}):$\ $L\in
 Q_{\theta\mathbf{b},\beta}(\mathbb{B}_n).$\ Indeed, we have by definition
\begin{gather*}\lambda_{1}^{\theta\mathbf{b}}(z,t_{0},\eta,L)=\inf
\left\{\frac{L(z+t\theta\mathbf{b})}{L(z+t_{0} \theta\mathbf{b})}:
|t-t_{0}|\leq\frac{\eta}{L(z+t_{0}\theta\mathbf{b})}\right\}=
\\ =\inf \left\{\frac{L(z+(t\theta)\mathbf{b})}{L(z+(t_{0}
\theta)\mathbf{b})}: |\theta t-\theta
t_{0}|\leq\frac{|\theta|\eta}{L(z+(t_{0}\theta)\mathbf{b})}\right\}=
\lambda_{1}^{\mathbf{b}}(z,\theta t_{0},|\theta|\eta,L).
\end{gather*}
Therefore, we get 
\begin{gather*}
\lambda_{1}^{\theta\mathbf{b}}\!(\!\eta,L\!)\!=\!\inf\!\{\!\lambda_{1}^{\theta\mathbf{b}}\!(z,\!\eta,L\!):
z\in\mathbb{B}_n\}=\inf\{\inf\{\lambda_{1}^{\theta\mathbf{b}}(z,t_{0},\eta,L):
 t_0\in S_{z}\}: z\in\mathbb{B}_n\}=\\ \!=\!\inf\{\inf\{
\lambda_{1}^{\mathbf{b}}(z,\theta t_{0},|\theta|\eta,L)\!:
\theta t_0\in S_{z}\}\!: z\in\mathbb{B}_n\}\!=\!
\inf\{\lambda^{\mathbf{b}}_{1}(z,|\theta|\eta,L):
z\in\mathbb{B}_n\}=\\=\lambda_{1}^{\mathbf{b}}(|\theta|\eta,L)>0,
\end{gather*} because $L\in Q_{\mathbf{b},\beta}(\mathbb{B}_n).$ Similarly, we prove that
$\lambda_{2}^{\theta\mathbf{b}}(\eta,L)=\lambda_{2}^{\mathbf{b}}(|\theta|\eta,L)<+\infty.$ But $|\theta|\eta\in[0,\beta].$ So $\eta\in[0,\beta/|\theta|].$
Thus, $L\in  Q_{\theta\mathbf{b},\beta/|\theta|}(\mathbb{B}_n)$.

Let $L^*=|\theta|\cdot L.$
Using definition of $\lambda_1^{\mathbf{b}}(z,t_{0},\eta,L^*)$  we have
\begin{gather*}\lambda_{1}^{\theta\mathbf{b}}(z,t_{0},\eta,L^*)=\inf
\left\{\frac{L^*(z+t\theta\mathbf{b})}{L^*(z+t_{0} \theta\mathbf{b})}:
|t-t_{0}|\leq\frac{\eta}{L^*(z+t_{0}\theta\mathbf{b})}\right\}=\\ =\inf
\left\{\frac{|\theta|L(z+t\theta\mathbf{b})}{|\theta|L(z+t_{0} \theta\mathbf{b})}: 
|t-t_{0}|\leq\frac{\eta}{|\theta|L(z+t_{0}\theta\mathbf{b})}\right\}=
\inf \left\{\frac{L(z+(t\theta)\mathbf{b})}{L(z+(t_{0}
\theta)\mathbf{b})}: 
\right. \\ \left. 
|\theta t-\theta
t_{0}|\leq\frac{\eta}{L(z+(t_{0}\theta)\mathbf{b})}\right\}
=\lambda_{1}^{\mathbf{b}}(z,\theta t_{0},\eta,L).
\end{gather*}

Therefore, we obtain 
\begin{gather*}
\lambda_{1}^{\theta\mathbf{b}}(\eta,L^*)=\inf\{\lambda_{1}^{\theta\mathbf{b}}(z,\eta,L^*):
z\in\mathbb{B}_n\}=\\ =\inf\{\inf\{\lambda_{1}^{\theta\mathbf{b}}(z,t_{0},\eta,L^*):
\theta t_0\in S_{z}\}: z\in\mathbb{B}_n\}=\\ \!=\!\inf\{\inf\{
\lambda_{1}^{\mathbf{b}}(z,\theta t_{0},\eta,L):
\theta t_0\in S_{z}\}: z\in\mathbb{B}_n\}\!=\! \\ = 
\inf\{\lambda^{\mathbf{b}}_{1}(z,\eta,L):
z\in\mathbb{B}_n\}\!=\!\lambda_{1}^{\mathbf{b}}\!(\!\eta,L)\!>\!0,
\end{gather*} because $L\in Q_{\mathbf{b},\beta}(\mathbb{B}_n).$
Similarly, we prove that
$\lambda_{2}^{\theta\mathbf{b}}(\eta,L^*)=\lambda_{2}^{\mathbf{b}}(\eta,L)<+\infty.$
Thus, $L^*=|\theta|\cdot L\in  Q_{\theta\mathbf{b},\beta}(\mathbb{B}_n)$.

\textbf{2.}
It remains to prove a second part.% of Lemma \ref{lemma1}.

If $z^0+t_0(\mathbf{b}_1+\mathbf{b}_2)\in\mathbb{B}_n$ and $|t-t_0|\leq
\frac{\eta}{L(z^0+t_0(\mathbf{b}_1+\mathbf{b}_2))}$ then
$z^0+t\mathbf{b}_1+t_0\mathbf{b}_2\in\mathbb{B}_n$ and $z^0+t_0\mathbf{b}_1+t\mathbf{b}_2\in\mathbb{B}_n.$
Indeed,  we have 
\begin{gather*}
 |z^0+t\mathbf{b}_1+t_0\mathbf{b}_2|\!\leq\! |z^0+t_0\mathbf{b}_1+t_0\mathbf{b}_2|+|t-t_0|\cdot|\mathbf{b}_1|\leq
 |z^0+t_0\mathbf{b}_1+t_0\mathbf{b}_2|+ \\+\frac{\eta|\mathbf{b}_1|}{L(z^0+t_0(\mathbf{b}_1+\mathbf{b}_2))}
  <
|z^0+t_0\mathbf{b}_1+t_0\mathbf{b}_2|+\frac{\beta|\mathbf{b}_1|}{\frac{\beta
 \max\{|\mathbf{b}_1|,|\mathbf{b}_2|,|\mathbf{b}_1+\mathbf{b}_2|\}}{1-|z^0+t_0\mathbf{b}_1+t_0\mathbf{b}_2|}} \leq 1.
\end{gather*}
Thus, $z^0+t\mathbf{b}_1+t_0\mathbf{b}_2\in\mathbb{B}_n.$

Denote $L^*(z)=\min\{\lambda_2^{\mathbf{b}_1}(\beta,L),\lambda_2^{\mathbf{b}_2}(\beta,L)\}\cdot L(z).$
Assume that $$\min\!\{\lambda_2^{\mathbf{b}_1}(\beta,L),\lambda_2^{\mathbf{b}_2}(\beta,L)\}\!=\! \lambda_2^{\mathbf{b}_2}(\beta,L).$$

Using definitions of
$\lambda_1^{\mathbf{b}}(\eta,L),$ $\lambda_2^{\mathbf{b}}(\eta,L)$ and
$Q_{\mathbf{b},\beta}(\mathbb{B}_n)$ we obtain that
\begin{gather}
 \inf\left\{\frac{L^*(z^0+t(\mathbf{b}_1+\mathbf{b}_2))}{L^*(z^0+t_0(\mathbf{b}_1+\mathbf{b}_2))}: \
 |t-t_0|\leq \frac{\eta}{L^*(z^0+t_0(\mathbf{b}_1+\mathbf{b}_2))}\right\}
%  = \nonumber \\
%  =  \inf\left\{ \frac{L(z^0+t_1\mathbf{b}_1+t_2\mathbf{b}_2)}{L(z^0+t_0(\mathbf{b}_1+\mathbf{b}_2))}:
%  |t_2-t_0|\leq \right. \nonumber\\ \leq \left. \frac{\eta}{L(z^0+t_0(\mathbf{b}_1+\mathbf{b}_2))}\right\}= \\ =
%   \inf \left\{
%  \frac{L^*(z^0+t\mathbf{b}_1+t\mathbf{b}_2)}{L^*(z^0+t_0\mathbf{b}_1+t\mathbf{b}_2)} \cdot
% \frac{L^*(z^0+t_0\mathbf{b}_1+t\mathbf{b}_2)}{L^*(z^0+t_0(\mathbf{b}_1+\mathbf{b}_2))}:
%  |t-t_0|\leq
% \frac{\eta}{L^*(z^0+t_0(\mathbf{b}_1+\mathbf{b}_2))}\right\} 
\geq \nonumber\\ \geq
\inf\left\{\frac{L^*(z^0+t\mathbf{b}_1+t\mathbf{b}_2)}{L^*(z^0+t_0\mathbf{b}_1+t\mathbf{b}_2)}: |t-t_0|\leq
\frac{\eta}{L^*(z^0+t_0(\mathbf{b}_1+\mathbf{b}_2))} \right\}\times \nonumber\\ \times
\inf\left\{ \frac{L^*(z^0+t_0\mathbf{b}_1+t\mathbf{b}_2)}{L^*(z^0+t_0(\mathbf{b}_1+\mathbf{b}_2))}:
 |t-t_0|\leq
\frac{\eta}{L^*(z^0+t_0(\mathbf{b}_1+\mathbf{b}_2))}
\right\}= \nonumber \\ \!=\!
\inf\left\{\frac{\lambda_2^{\mathbf{b}_2}(\beta,L)L(z^0\!+\!t\mathbf{b}_1\!+\!t\mathbf{b}_2)}{\lambda_2^{\mathbf{b}_1}
(\beta,L)L(z^0\!+\!t_0\mathbf{b}_1\!+\!t\mathbf{b }_2)}:  |t-t_0|\leq
\frac{\eta}{\lambda_2^{\mathbf{b}_2}(\beta,L)L(z^0\!+\!t_0(\mathbf{b}_1\!+\!\mathbf{b}_2))} \right\}\times \nonumber \\ \times
\inf\left\{
\frac{\lambda_2^{\mathbf{b}_2}(\beta,L)L(z^0\!+\!t_0\mathbf{b}_1\!+\!t\mathbf{b}_2)}{\lambda_2^{\mathbf{b}_1}
(\beta,L)L(z^0\!+\!t_0(\mathbf{b}_1\!+\!\mathbf{b}_2))}:
 |t-t_0|\!\leq\!
\frac{\eta}{\lambda_2^{\mathbf{b}_2}(\beta,L)L(z^0+t_0(\mathbf{b}_1+\mathbf{b}_2))}
\right\}= \nonumber \\ =
\inf\left\{\frac{L(z^0+t\mathbf{b}_1+t\mathbf{b}_2)}{L(z^0+t_0\mathbf{b}_1+t\mathbf{b }_2)}:
|t-t_0|\leq
\frac{\eta}{\lambda_2^{\mathbf{b}_2}(\beta,L)L(z^0+t_0(\mathbf{b}_1+\mathbf{b}_2))} \right\}\times \nonumber \\ \times
\inf\left\{
\frac{L(z^0+t_0\mathbf{b}_1+t\mathbf{b}_2)}{L(z^0+t_0(\mathbf{b}_1+\mathbf{b}_2))}:
 |t-t_0|\leq
\frac{\eta}{\lambda_2^{\mathbf{b}_2}(\beta,L)L(z^0+t_0(\mathbf{b}_1+\mathbf{b}_2))}
\right\}\geq \nonumber \\
% \left|\begin{array}{c}\text{we use}\\ \lambda_2^{\mathbf{b}_2}(\beta,L)\geq 1\end{array}\right|
\geq
\inf\left\{\frac{L(z^0+t\mathbf{b}_1+t\mathbf{b}_2)}{L(z^0+t_0\mathbf{b}_1+t\mathbf{b }_2)}:
|t-t_0|\leq
\frac{\eta}{\lambda_2^{\mathbf{b}_2}(\beta,L)L(z^0+t_0(\mathbf{b}_1+\mathbf{b}_2))} \right\}\times \nonumber \\ \times
\inf\left\{
\frac{L(z^0+t_0\mathbf{b}_1+t\mathbf{b}_2)}{L(z^0+t_0(\mathbf{b}_1+\mathbf{b}_2))}:
 |t-t_0|\leq
\frac{\eta}{L(z^0+t_0(\mathbf{b}_1+\mathbf{b}_2))}
\right\}\geq \nonumber\\
\geq
\inf\left\{\frac{L(z^0+t\mathbf{b}_1+t\mathbf{b}_2)}{L(z^0+t_0\mathbf{b}_1+t\mathbf{b }_2)}:
|t-t_0|\leq
\frac{\eta}{\lambda_2^{\mathbf{b}_2}(\beta,L)L(z^0+t_0(\mathbf{b}_1+\mathbf{b}_2))} \right\} \times \nonumber \\
\times \lambda_1^{\mathbf{b}_2}(z^0+t_0\mathbf{b}_1,t_0,\eta,L)
\geq \lambda_1^{\mathbf{b}_2}(\eta,L) \frac{L(z^0+\hat{t}\mathbf{b}_1+\hat{t}\mathbf{b}_2)}{L(z^0+t_0\mathbf{b}_1+\hat{t}\mathbf{b
}_2)}
\label{kin1}
\end{gather}
where $\hat{t}$ is a point at which infimum is attained 
\begin{gather*}
\frac{L(z^0\!+\!\hat{t}\mathbf{b}_1\!+\!\hat{t}\mathbf{b}_2)}{L(z^0\!+\!t_0\mathbf{b}_1\!+\!\hat{t}\mathbf{b
}_2)} \!=\!\inf\!\left\{\frac{L(z^0\!+\!t\mathbf{b}_1+t\mathbf{b}_2)}{L(z^0\!+\!t_0\mathbf{b}_1+t\mathbf{b }_2)}:
|t-t_0|\!\leq\!
\frac{\eta}{\lambda_2^{\mathbf{b}_2}(\beta,L)L(z^0+t_0(\mathbf{b}_1+\mathbf{b}_2))} \right\}.
\end{gather*}
% because $L(z)$ is a positive continuous function.

But $L\in Q_{\mathbf{b}_{2},\beta}(\mathbb{B}_n),$ then for all $\eta\in[0,\beta]$
$$
\!\sup\!\left\{\frac{L(z^0\!+\!t_0\mathbf{b}_1\!+\!t\mathbf{b}_2)}{L\!(\!z^0\!+\!t_0\mathbf{b}_1\!+\!t_0\mathbf{b}_2\!)}:
|t\!-\!t_0|\!\leq \!
\frac{\eta}{L\!(\!z^0\!+\!t_0(\mathbf{b}_1\!+\!\mathbf{b}_2\!))}\right\}\!%=\!\lambda_2^{\mathbf{b}_2}(z^0+t_0\mathbf{b}_1,t_0,\eta,L)\!
\leq\! \lambda_2^{\mathbf{b}_2}(\eta,L)\!<\!\infty.
$$
Hence, $L(z^0+t_0\mathbf{b}_1+t\mathbf{b}_2)\leq \lambda_2^{\mathbf{b}_2}(\eta,L) \cdot
L(z^0+t_0\mathbf{b}_1+t_0\mathbf{b}_2),$ i.e. for $t=\hat{t}$ we have $L(z^0+t_0\mathbf{b}_1+t_0\mathbf{b}_2)\geq
\frac{L(z^0+t_0\mathbf{b}_1+\hat{t}\mathbf{b}_2)}{\lambda_2^{\mathbf{b}_2}(\eta,L)}.$
Using a proved inequality and \eqref{kin1}, we obtain
\begin{gather}
 \inf\left\{\frac{L^*(z^0+t(\mathbf{b}_1+\mathbf{b}_2))}{L^*(z^0+t_0(\mathbf{b}_1+\mathbf{b}_2))}: \
 |t-t_0|\leq \frac{\eta}{L^*(z^0+t_0(\mathbf{b}_1+\mathbf{b}_2))}\right\} \geq \nonumber \\ \nonumber \!\geq\!
 \lambda_1^{\mathbf{b}_2}(\eta,L)\cdot \inf\left\{
 \frac{L(z^0\!+\!t\mathbf{b}_1\!+\!\hat{t}\mathbf{b}_2)}{L(z^0\!+\!t_0\mathbf{b}_1\!+\!\hat{t}\mathbf{b}_2)}:
 |t\!-\!t_0|\!\leq\! \frac{\eta}{\lambda_2^{\mathbf{b}_2}(\beta,L)L(z^0\!+\!t_0(\mathbf{b}_1\!+\!\mathbf{b}_2))}
 \right\} \geq\\ \! \geq\!
 \lambda_1^{\mathbf{b}_2}(\eta,L)\cdot\inf\left\{
 \frac{L(z^0\!+\!t\mathbf{b}_1\!+\!\hat{t}\mathbf{b}_2)}{L(z^0\!+\!t_0\mathbf{b}_1\!+\!\hat{t}\mathbf{b}_2)}:
 |t\!-\!t_0|\!\leq\! \frac{\eta\lambda_2^{\mathbf{b}_2}(\eta,L)}{\lambda_2^{\mathbf{b}_2}(\beta,L)L(z^0+t_0\mathbf{b}_1+\hat{t}\mathbf{b}_2)}
 \right\}\geq \nonumber
 \\ 
 %\!\left|\begin{array}{c}\text{But } \lambda_2^{\mathbf{b}_2}(\eta,L) \text{ is a}\\  \text{nondecreasing function}\end{array}\right| 
 \!\geq\!
 \lambda_1^{\mathbf{b}_2}(\eta,L)\cdot\inf\left\{
 \frac{L(z^0+t\mathbf{b}_1+\hat{t}\mathbf{b}_2)}{L(z^0+t_0\mathbf{b}_1+\hat{t}\mathbf{b}_2)}:
   |t-t_0|\!\leq\! \frac{\eta}{L(z^0+t_0\mathbf{b}_1+\hat{t}\mathbf{b}_2)}
 \right\} =\nonumber \\ = \lambda_1^{\mathbf{b}_2}(\eta,L) \lambda_1^{\mathbf{b}_1}(z^0+\hat{t}\mathbf{b}_2,t_0,\eta,L)
 \geq \lambda_1^{\mathbf{b}_2}(\eta,L) \lambda_1^{\mathbf{b}_1}(\eta,L).\nonumber
%  \limits_{t_2}\left\{ \inf\limits_{t_1}  \left\{
%  \frac{L(z^0+t_1\mathbf{b}_1+t_2\mathbf{b}_2)}{L(z^0+t_0\mathbf{b}_1+t_2\mathbf{b}_2)}:    |t_1-t_0|\leq
% \frac{\eta\cdot \lambda_2^{\mathbf{b}_2}(\eta)}{L(z^0+t_0\mathbf{b}_1+t_2\mathbf{b}_2)}\right\}: \right. \nonumber \\
% \left. |t_2-t_0|\leq
% \frac{\eta}{L(z^0+t_0(\mathbf{b}_1+\mathbf{b}_2))}\right\}
% \geq \lambda_1^{\mathbf{b}_2}(\eta) \cdot \inf\limits_{t_2}
% \left\{\lambda_1^{\mathbf{b}_1}(z^0+t_2\mathbf{b}_2,t_0,\eta\lambda_2^{\mathbf{b}_2}(\eta)): \right. \nonumber\\ \left.
% |t_2-t_0| \leq \frac{\eta}{L(z^0+t_0(\mathbf{b}_1+\mathbf{b}_2))}\right\} \geq \lambda_1^{\mathbf{b}_2}(\eta) \cdot
% \lambda_1^{\mathbf{b}_1}(\eta\cdot \lambda_2^{\mathbf{b}_2}(\eta)).
 \end{gather}

 Therefore,
$\lambda_{1}^{\mathbf{b}_{1}+\mathbf{b}_{2}}(\eta,L^*)\geq \lambda_1^{\mathbf{b}_2}(\eta,L) \lambda_1^{\mathbf{b}_1}(\eta,L)>0.$
By analogy, we can prove that for
all $\eta\in[0,\beta]$ \ 
$\lambda_{2}^{\mathbf{b}_{1}+\mathbf{b}_{2}}(\eta,L^*)<+\infty.$  Thus, $L^*\in
Q_{\mathbf{b}_{1}+\mathbf{b}_{2},\beta}(\mathbb{B}_n).$

\end{proof}

% \begin{center}
\section{Criteria of $L$-index boundedness in direction,
related to the behaviour of the function $F$.} 
% \end{center}  

% {\bf{$4^{0}.$ }} 

The following theorem is an analogue
of Theorem 2 from \cite{BandSk}.
\begin{theorem}
\label{te2} Let $\beta>1$ and $L\in {Q}_{\mathbf{b},\beta}(\mathbb{B}_n)$. Analytic in $\mathbb{B}_n$ function
$F(z)$ is of bounded $L$-index in  direction $\mathbf{b}\in\mathbb{C}^{n}$ if and only if
 for every $\eta,$ $0<\eta\leq\beta,$ there exist
$n_{0}=n_{0}(\eta)\in \mathbb{Z}_{+}$ and $P_{1}=P_{1}(\eta)\geq 1$
 that for each $z\in
\mathbb{B}_n$ and each $t_{0}\in  S_z$ there exists $k_{0}=k_{0}(t_{0},z)\in \mathbb{Z}_{+},$\ with
$0\leq k_{0}\leq n_{0},$  and the following inequality holds
\begin{equation} \label{poch1}
\max\left\{\left|\frac{\partial^{k_{0}} F(z+t\mathbf{b})}
{\partial\mathbf{b}^{k_{0}}}\right|: |t-t_{0}|\leq\frac{\eta}
{L(z+t_{0}\mathbf{b})}\right\}\leq
P_{1}\left|\frac{\partial^{k_{0}}F(z+t_{0}\mathbf{b})}
{\partial\mathbf{b}^{k_{0}}}\right|.
\end{equation}
\end{theorem}

\begin{proof}

\textbf{Necessity.} Let $F$ be of bounded $L$-index in  direction $\mathbf{b}$ and $N_{\mathbf{b}}(F;L)\equiv N<+\infty$.
We denote $$q(\eta)=[2\eta(N+1)(\lambda_{2}^{\mathbf{b}}(\eta))^{N+1}(\lambda_{1}^{\mathbf{b}}(\eta))^{-N}]+1,$$
where $[a]$ is an entire part of number $a\in\mathbb{R}.$ 
 For $z\in\mathbb{B}_n,$ $t_{0}\in  S_z$ and $p \in\{0,1,\ldots,q(\eta)\}$
we put
\begin{gather*} R^{\mathbf{b}}_{p}(z,t_{0},\eta)\!=\!\max\left\{\frac{1}{k!L^{k}(z\!+\!t\mathbf{b})}\left|\frac{\partial^{k}F(z+t\mathbf{b})}
{\partial\mathbf{b}^{k}}\right|: 
|t-t_{0}|\!\leq\!\frac{p\eta}{q(\eta)L(z+t_{0}\mathbf{b})}, 0\leq k\leq N\right\}.\end{gather*}
and
\begin{gather*} \widetilde{R}^{\mathbf{b}}_{p}(z,t_{0},\eta)\!=\!\max\left\{\frac{1}{k!L^{k}(z+t_{0}\mathbf{b})}\left|\frac{\partial^{k}F(z+t\mathbf{b})}
{\partial\mathbf{b}^{k}}\right|: 
|t\!-\!t_{0}|\!\leq\!\frac{p\eta}{q(\eta)L(z\!+\!t_{0}\mathbf{b})}, 0\leq k\leq N\right\}.\end{gather*}

But $|t-t_{0}|\leq\displaystyle\frac{\mathstrut p\eta}{q(\eta)L(z+t_{0}\mathbf{b})}\leq\displaystyle\frac{\mathstrut \eta}{L(z+t_{0}\mathbf{b})}
\leq\displaystyle\frac{\mathstrut \beta}{L(z+t_{0}\mathbf{b})},$ then  $$\lambda_{1}^{\mathbf{b}}\left(z,t_{0},\displaystyle\frac{\mathstrut p\eta}{q(\eta)}\right)\geq\lambda_{1}^{\mathbf{b}}(z,t_{0},\eta)\geq \lambda_{1}^{\mathbf{b}}(\eta), \ 
 \lambda_{2}^{\mathbf{b}}\left(z,t_{0},\displaystyle\frac{\mathstrut p\eta}{q(\eta)}\right)\leq\lambda_{2}^{\mathbf{b}}(z,t_{0},\eta)\leq \lambda_{2}^{\mathbf{b}}(\eta).$$
Clearly, these quantities $R^{\mathbf{b}}_{p}(z,t_{0},\eta),$ $\widetilde{R}^{\mathbf{b}}_{p}(z,t_{0},\eta)$ are defined. Besides,
\begin{gather}
 R^{\mathbf{b}}_{p}(z,t_{0},\eta)=\max\left\{\frac{1}{k!L^{k}(z+t_{0}\mathbf{b})}\left|\frac{\partial^{k}
 F(z+t\mathbf{b})}{\partial\mathbf{b}^{k}}\right|\left(\frac{L(z+t_{0}\mathbf{b})}{L(z+t\mathbf{b})}\right)^{k}:
\right. \nonumber\\
\left. |t-t_{0}|\leq\frac{p\eta}{q(\eta)L(z+t_{0}\mathbf{b})}, 0\leq k\leq N\right\}\leq\nonumber \\
 \leq\max\left\{\frac{1}{k!L^{k}(z+t_{0}\mathbf{b})}\left|\frac{\partial^{k}
 F(z+t\mathbf{b})}{\partial\mathbf{b}^{k}}\right|\left(\frac{1}{\lambda_{1}^{\mathbf{b}}(z,t_{0},
 \frac{p\eta}{q(\eta)})}\right)^{k}:
\right. \nonumber\\ \left.
|t-t_{0}|\leq\frac{p\eta}{q(\eta)L(z+t_{0}\mathbf{b})}, 0\leq k\leq N\right\}\leq\nonumber\\
 \leq\!\max\left\{\frac{1}{k!L^{k}(z+t_{0}\mathbf{b})}\left|\frac{\partial^{k} F(z+t\mathbf{b})}{\partial\mathbf{b}^{k}}\right|\left(\frac{1}{\lambda_{1}^{\mathbf{b}}(\eta)}\right)^{k}: \nonumber \right.\\ \left.
|t-t_{0}|\leq\frac{p\eta}{q(\eta)L(z+t_{0}\mathbf{b})}, 0\leq k\leq N\right\}\leq\nonumber\\
  \leq\!\left(\frac{1}{\lambda_{1}^{\mathbf{b}}(\eta)}\right)^{N}\max\left\{\frac{1}{k!L^{k}(z+t_{0}\mathbf{b})}\left|
 \frac{\partial^{k} F(z+t\mathbf{b})}{\partial\mathbf{b}^{k}}\right|: \nonumber\right.\\ \left.
 |t-t_{0}|\leq\frac{p\eta}{q(\eta)L(z+t_{0}\mathbf{b})}, 0\leq k\leq N\right\}=
\widetilde{R}^{\mathbf{b}}_{p}(z,t_{0},\eta)(\lambda_{1}^{\mathbf{b}}(\eta))^{-N} \label{eq_th_11_15}
\end{gather}
and
 \begin{gather}
  \widetilde{R}^{\mathbf{b}}_{p}(z,t_{0},\eta)=\max\left\{\frac{1}{k!L^{k}(z+t\mathbf{b})}\left|\frac{\partial^{k}
F(z+t\mathbf{b})}{\partial\mathbf{b}^{k}}\right|\left(\frac{L(z+t\mathbf{b})}{L(z+t_{0}\mathbf{b})}\right)^{k}:
 \nonumber\right.\\ \left. |t-t_{0}|\leq\frac{p\eta}{q(\eta)L(z+t_{0}\mathbf{b})}, 0\leq k\leq N\right\}\leq\nonumber \\
 \leq\max\left\{\frac{1}{k!L^{k}(z+t\mathbf{b})}\left|\frac{\partial^{k}
 F(z+t\mathbf{b})}{\partial\mathbf{b}^{k}}\right|\left(\lambda_{2}^{\mathbf{b}}\left(z,t_{0},
 \frac{p\eta}{q(\eta)}\right)\right)^{k}: \nonumber\right.\\ \left.
|t-t_{0}|\leq\frac{p\eta}{q(\eta)L(z+t_{0}\mathbf{b})}, 0\leq k\leq
N\right\}\leq\nonumber\\
 \leq\max\left\{\frac{\left(\lambda_{2}^{\mathbf{b}}(\eta)\right)^{k}}{k!L^{k}(z+t\mathbf{b})}\left|\frac{\partial^{k}
 F(z+t\mathbf{b})}{\partial\mathbf{b}^{k}}\right|:
|t-t_{0}|\leq\frac{p\eta}{q(\eta)L(z+t_{0}\mathbf{b})}, \right.\nonumber\\ \left. 0\leq k\leq N\right\}
 \leq\left(\lambda_{2}^{\mathbf{b}}(\eta)\right)^{N}\max\left\{\frac{1}{k!L^{k}(z+t\mathbf{b})}\left|\frac{\partial^{k}
F(z+t\mathbf{b})}{\partial\mathbf{b}^{k}}\right|: \right.\nonumber\\ \left. |t-t_{0}|\leq\frac{p\eta}{q(\eta)L(z+t_{0}\mathbf{b})}, 0\leq k\leq N\right\}  =R^{\mathbf{b}}_{p}(z,t_{0},\eta)(\lambda_{2}^{\mathbf{b}}(\eta))^{N}.\label{eq_th_11_16}
 \end{gather}
Let $k^{z}_{p}\in\mathbb{Z},$ $0\leq k_{p}^z\leq N,$ and $t^{z}_{p}\in\mathbb{C},$ $|t^{z}_{p}-t_{0}|\leq\displaystyle\frac{\mathstrut p\eta}{q(\eta)L(z+t_{0}\mathbf{b})},$ be such that
\begin{equation}\label{eq_th_1_1_165}
\frac{1}{k^{z}_{p}!L^{k^{z}_{p}}
(z+t_{0}\mathbf{b})}\left| \frac{\partial^{k_{p}^{z}} F(z+t_{p}^{z}\mathbf{b})}{\partial\mathbf{b}^{k_{p}^{z}}}\right|=\widetilde{R}_{p}^{\mathbf{b}}(z,t_{0},\eta).       \end{equation}
For every given $z\in\mathbb{B}_n$ a function $F(z+t\mathbf{b})$ and its directional derivative are analytic. By the maximum modulus principle an equality (\ref{eq_th_1_1_165}) holds for such $t_{p}^{z},$ that $$|t_{p}^{z}-t_{0}|=\displaystyle\frac{\mathstrut p\eta}{q(\eta)L(z+t_{0}\mathbf{b})}.$$ We put $\widetilde{t_{p}^{z}}=t_{0}+\frac{p-1}{p}(t_{p}^{z}-t_{0}).$
Then
\begin{equation} \label{eq_th_1.1_1.7}
|\widetilde{t_{p}^{z}}-t_{0}|=\frac{(p-1)\eta}{q(\eta)L(z+t_{0}\mathbf{b})}
\end{equation}
and \begin{equation} \label{eq_th_1.1_1.8}
|\widetilde{t_{p}^{z}}-t^z_{p}|=\frac{|t_{p}^{z}-t_{0}|}{p}=\frac{\eta}{q(\eta)L(z+t_0\mathbf{b})}.
  \end{equation}
In view of (\ref{eq_th_1.1_1.7}) and the definition of $\widetilde{R}^{\mathbf{b}}_{p-1}(z,t_{0},\eta),$ we obtain that
$$\widetilde{R}^{\mathbf{b}}_{p-1}(z,t_{0},\eta)\geq\frac{1}{k_{p}^{z}!L^{k_{p}^{z}}(z+t_{0}\mathbf{b})} \left|\frac{\partial^{k_{p}^{z}} F(z+\widetilde{t_{p}^{z}}\mathbf{b})}{\partial\mathbf{b}^{k_{p}^{z}}}\right|.$$
Therefore, this inequality holds 
\begin{gather}
\!0\!\leq\!\widetilde{R}_{p}^{\mathbf{b}}(z,t_{0},\eta)-\widetilde{R}_{p-1}^{\mathbf{b}}(z,t_{0},\eta)\!
\leq\! \frac{\left|\displaystyle\frac{\mathstrut\partial^{k_{p}^{z}} F(z+t_{p}^{z}\mathbf{b})}{\partial\mathbf{b}^{k_{p}^{z}}}\right|-\left|\displaystyle\frac{\mathstrut
\partial^{k_{p}^{z}} F(z+\widetilde{t_{p}^{z}}\mathbf{b})}{\partial\mathbf{b}^{k_{p}^{z}}}\right|}{k_{p}^{z}!L^{k_{p}^{z}}
(z+t_{0}\mathbf{b})}\!=\nonumber
\\  \label{eg_th_11_1_75}
=\frac{1}{k_{p}^{z}!L^{k_{p}^{z}}(z+t_{0}\mathbf{b})}\int_{0}^{1}\frac{d}{ds}\left|\frac{\partial^{k_{p}^{z}} F(z +(\widetilde{t}_{p}^{z}+s(t_{p}^{z}-\widetilde{t}_{p}^{z}))\mathbf{b})}{\partial\mathbf{b}^{k_{p}^{z}}}\right|ds.
\end{gather}
For every analytic complex-valued function of real variable $\varphi(s),$ $s\in\mathbb{R},$ the inequality $\frac{d}{ds}|\varphi(s)|\leq\left|\frac{d}{ds}\varphi(s)\right|$ holds except the points where $\varphi(s)=0.$  Applying this inequality to \eqref{eg_th_11_1_75} and using a mean value theorem, we have 
\begin{gather*}
 \widetilde{R}_{p}^{\mathbf{b}}(z,t_{0},\eta)-\widetilde{R}_{p-1}^{\mathbf{b}}(z,t_{0},\eta)\!\leq 
\frac{|t^{z}_{p}-\widetilde{t}_{p}^{z}|}{k_{p}^{z}!L^{k^{z}_{p}}(z+t_{0}\mathbf{b})}  \int_{0}^{1}
\left|\frac{\partial^{k_{p}^{z}+1} F(z+(\widetilde{t}_{p}^{z}+s(t_{p}^{z}-\widetilde{t}_{p}^{z}))\mathbf{b})}{\partial\mathbf{b}^{k_{p}^{z}+1}}\right|ds=\\
=\frac{|t^{z}_{p}-\widetilde{t}_{p}^{z}|}{k_{n}^{p}!L^{k^{z}_{p}}(z+t_{0}\mathbf{b})}
\left|\frac{\partial^{k_{p}^{z}+1} F(z+(\widetilde{t}_{p}^{z}+s^{*}(t_{p}^{z}-\widetilde{t}_{p}^{z}))\mathbf{b})}{\partial\mathbf{b}^{k_{p}^{z}+1}}\right|=\\=\frac{1}{(k_{p}^{z}+1)!L^{k^{z}_{p}+1}(z+t_{0}\mathbf{b})}
\left|\frac{\partial^{k_{p}^{z}+1} F(z\!+\!(\widetilde{t}_{p}^{z}+s^{*}(t_{p}^{z}\!-\!\widetilde{t}_{p}^{z}))\mathbf{b})}{\partial\mathbf{b}^{k_{p}^{z}+1}}\right| \cdot 
L(z+t_{0}\mathbf{b})(k_{p}^{z}+1)|t^{z}_{p}\!-\!\widetilde{t}_{p}^{z}|,
\end{gather*}
where $s^{*}\in [0,1].$

The point $\widetilde{t}_{p}^{z}+s^{*}(t_{p}^{z}-\widetilde{t}_{p}^{z})$ lies into the set $$\left\{t\in\mathbb{C}: |t-t_{0}|\leq\displaystyle\frac{\mathstrut p\eta}{q(\eta)L(z+t_{0}\mathbf{b})}\leq\displaystyle\frac{\mathstrut\eta}{L(z+t_{0}\mathbf{b})}\right\}.$$

Using $L$-index boundedness in  direction $\mathbf{b}$ of function $F,$ definition $q(\eta),$ inequality (\ref{eq_th_11_15}) and   (\ref{eq_th_1.1_1.8}), for $k_{p}^{z}\leq N$ we have
\begin{gather*}
 \!\widetilde{R}_{p}^{\mathbf{b}}(z,t_{0},\eta)-\widetilde{R}_{p-1}^{\mathbf{b}}(z,t_{0},\eta)\leq
\frac{1}{(k_{n}^{z}+1)!L^{k^{z}_{n}+1}(z+(\widetilde{t}_{p}^{z}+s^{*}(t_{p}^{z}-\widetilde{t}_{p}^{z}))\mathbf{b})}\!
\times\\ \!\times\!
 \left|\frac{\partial^{k_{p}^{z}+1}
F(z+(\widetilde{t}_{p}^{z}+s^{*}(t_{p}^{z}-\widetilde{t}_{p}^{z}))\mathbf{b})}
 {\partial\mathbf{b}^{k_{p}^{z}+1}}\!\right|\!\left(\frac{L(z+(\widetilde{t}_{p}^{z}+s^{*}(t_{p}^{z}-
 \widetilde{t}_{p}^{z}))\mathbf{b})}{L(z+t_{0}\mathbf{b})}\right)^{k_{p}^{z}+1}\!\times\\ \times
  L(z+t_{0}\mathbf{b})(k_{n}^{z}+1)|t^{z}_{p}-\widetilde{t}_{p}^{z}|
  \leq \eta\frac{N+1}{q(\eta)}(\lambda_{2}^{\mathbf{b}}(z,t_{0},\eta))^{N+1}\times\\ \times
 \max\left\{\frac{1}{k!L^{k}(z+(\widetilde{t}_{p}^{z}+s^{*}(t_{p}^{z}-\widetilde{t}_{p}^{z}))\mathbf{b})}
 \left|\frac{\partial^{k} F(z+(\widetilde{t}_{p}^{z}+s^{*}(t_{p}^{z}-\widetilde{t}_{p}^{z}))\mathbf{b})}
 {\partial\mathbf{b}^{k}}\right|: 0\leq k\leq N\right\}\leq \\
 \leq \eta\frac{N+1}{q(\eta)}(\lambda_{2}^{\mathbf{b}}(\eta))^{N+1}R_{p}^{\mathbf{b}}(z,t_{0},\eta)
\leq
 \frac{\eta (N+1)(\lambda_{2}^{\mathbf{b}}(\eta))^{N+1}(\lambda_{1}^{\mathbf{b}}(\eta))^{-N}}{[2\eta(N+1)
 \lambda_{2}^{\mathbf{b}}(\eta)(\lambda_{1}^{\mathbf{b}}(\eta))^{-N}]+1}\widetilde{R}_{p}^{\mathbf{b}}(z,
t_{0},\eta)\leq \\ 
 \leq \frac{1}{2}\widetilde{R}_{p}^{\mathbf{b}}(z,t_{0},\eta)
\end{gather*}
In the last inequality we used that $2a+1\geq [2a+1]=[2a]+1\geq 2a$ for $a\in\mathbb{R}.$ 

It follows that $\widetilde{R}_{p}^{\mathbf{b}}(z,t_{0},\eta)\leq 2\widetilde{R}_{p-1}^{\mathbf{b}}(z,t_{0},\eta).$ Using inequalities (\ref{eq_th_11_15}) and (\ref{eq_th_11_16}),  we deduce for $R_{p}^{\mathbf{b}}(z,t_{0},\eta)$
\begin{gather*} R_{p}^{\mathbf{b}}(z,t_{0},\eta)\leq 2(\lambda_{1}^{\mathbf{b}}(\eta))^{-N}\widetilde{R}^{\mathbf{b}}_{p-1}(z,t_{0},\eta)  \leq 2(\lambda_{2}^{\mathbf{b}}(\eta))^{N}(\lambda_{1}^{\mathbf{b}}(\eta))^{-N}R_{p-1}^{\mathbf{b}}(z,t_{0},\eta).\end{gather*}

Hence, we have 
\begin{gather}
\max\left\{\frac{1}{k!L^{k}(z+t\mathbf{b})}\left|\frac{\partial^{k} F(z+t\mathbf{b})}{\partial\mathbf{b}^{k}}\right|:  |t-t_{0}|\leq\frac{\eta}{L(z+t_{0}\mathbf{b})},
\right. \nonumber\\ \left. 0\leq k\leq N\right\}=
 R_{q(\eta)}^{\mathbf{b}}(z,t_{0},\eta)\leq
2(\lambda_{2}^{\mathbf{b}}(\eta))^{N}(\lambda_{1}^{\mathbf{b}}(\eta))^{-N}R_{q(\eta)-1}^{\mathbf{b}}
(z,t_{0},\eta)\leq
 \nonumber\\ \leq
(2(\lambda_{2}^{\mathbf{b}}(\eta))^{N}(\lambda_{1}^{\mathbf{b}}(\eta))^{-N})^{2}R_{q(\eta)-2}^{\mathbf{b}}(z,t_{0},\eta)\leq\cdots\leq \nonumber\\ \leq
 (2(\lambda_{2}^{\mathbf{b}}(\eta))^{N}(\lambda_{1}^{\mathbf{b}}(\eta))^{-N})^{q(\eta)}R_{0}^{\mathbf{b}}
(z,t_{0},\eta)=
 (2(\lambda_{2}^{\mathbf{b}}(\eta))^{N}(\lambda_{1}^{\mathbf{b}}(\eta))^{-N})^{q(\eta)}
\times \nonumber\\ \times\max\left\{\frac{1}{k!L^{k}(z+t_{0}\mathbf{b})}\left|\frac{\partial^{k}
F(z+t_{0}\mathbf{b})}{\partial\mathbf{b}^{k}}\right|: 0\leq k\leq N\right\}.
\label{eq_th_11_19}
\end{gather}
Let $k^{z}_{0}\in\mathbb{Z},$ $0\leq k^{z}_{0}=k_{0}^{z}(t_{0})\leq N,$ and $\widetilde{t^{z}}\in\mathbb{C},$ $|\widetilde{t^{z}}-t_{0}|=\displaystyle\frac{\mathstrut\eta}{L(z+t_{0}\mathbf{b})},$ be defined  as
\begin{gather*}
\frac{1}{k_{0}^{z}!L^{k_{0}^{z}}(z+t_{0}\mathbf{b})}\left|\frac{\partial^{k^{z}_{0}} F(z+t_{0}\mathbf{b})}{\partial\mathbf{b}^{k_{0}^{z}}}\right|
=\max_{0\leq k\leq N}\left\{\frac{1}{k!L^{k}(z+t_{0}\mathbf{b})}\left|\frac{\partial^{k} F(z+t_{0}\mathbf{b})}{\partial\mathbf{b}^{k}}\right|\right\}
 \end{gather*}
and
 $$\left|\frac{\partial^{k^{z}_{0}}
F(z+\widetilde{t}^{z}\mathbf{b})}{\partial\mathbf{b}^{k_{0}^{z}}}\right|= \max\left\{\left|\frac{\partial^{k^{z}_{0}} F(z+t\mathbf{b})}{\partial\mathbf{b}^{k^{z}_{0}}}\right|:
|t-t_{0}|\leq\frac{\eta}{L(z+t_{0}\mathbf{b})}\right\}. $$
From inequality (\ref{eq_th_11_19}) it follows
\begin{gather*}
 \frac{1}{k_{0}^{z}!L^{k_{0}^{z}}(z+\widetilde{t}^{z}\mathbf{b})}\left|\frac{\partial^{k^{z}_{0}}
F(z+\widetilde{t}^{z}\mathbf{b})}{\partial\mathbf{b}^{k_{0}^{z}}}\right|\leq\\
\leq\max\left\{\frac{1}{k_{0}^{z}!L^{k_{0}^{z}}(z+t\mathbf{b})}\left|\frac{\partial^{k_{0}^{z}} F(z+t\mathbf{b})}{\partial\mathbf{b}^{k_{0}^{z}}}\right|: |t-t_{0}|=\frac{\eta}{L(z+t_{0}\mathbf{b})}\right\}\leq\\
 \leq\max\left\{\frac{1}{k!L^{k}(z+t\mathbf{b})}\left|\frac{\partial^{k}
F(z+t\mathbf{b})}{\partial\mathbf{b}^{k}}\right|: |t-t_{0}|=\frac{\eta}{L(z+t_{0}\mathbf{b})}, \right.\\ \left. 0\leq k\leq N\right\}
 \leq (2(\lambda_{2}^{\mathbf{b}}(\eta))^{N}(\lambda_{1}^{\mathbf{b}}(\eta))^{-N})^{q(\eta)}
\frac{1}{k_{0}^{z}!L^{k_{0}^{z}}(z+t_{0}\mathbf{b})} \cdot \left|\frac{\partial^{k^{z}_{0}} F(z+t_{0}\mathbf{b})}{\partial\mathbf{b}^{k_{0}^{z}}}\right|.
\end{gather*}
Hence, we get
\begin{gather*}
\max\left\{\left|\frac{\partial^{k^{z}_{0}} F(z+t\mathbf{b})}{\partial\mathbf{b}^{k^{z}_{0}}}\right|:
|t-t_{0}|\leq\frac{\eta}{L(z+t_{0}\mathbf{b})}\right\}\leq\\ \leq (2(\lambda_{2}^{\mathbf{b}}(\eta))^{N}(\lambda_{1}^{\mathbf{b}}(\eta))^{-N})^{q(\eta)}
\left(\frac{L(z+\widetilde{t}^{z}\mathbf{b})}{L(z+t_{0}\mathbf{b})}\right)^{k^{z}_{0}}
\left|\frac{\partial^{k^{z}_{0}} F(z+t_{0}\mathbf{b})}{\partial\mathbf{b}^{k_{0}^{z}}}\right|\leq\\
 \leq (2(\lambda_{2}^{\mathbf{b}}(\eta))^{N}(\lambda_{1}^{\mathbf{b}}(\eta))^{-N})^{q(\eta)}
 (\lambda_{2}^{\mathbf{b}}(z,t_{0},\eta))^{N}\left|\frac{\partial^{k^{z}_{0}}
F(z+t_{0}\mathbf{b})}{\partial\mathbf{b}^{k_{0}^{z}}}\right|\leq\\
 \leq (2(\lambda_{2}^{\mathbf{b}}(\eta))^{N}(\lambda_{1}^{\mathbf{b}}(\eta))^{-N})^{q(\eta)}
 (\lambda_{2}^{\mathbf{b}}(\eta))^{N}\left|\frac{\partial^{k^{z}_{0}}
F(z+t_{0}\mathbf{b})}{\partial\mathbf{b}^{k_{0}^{z}}}\right|.
\end{gather*}
We proved (\ref{poch1}) with $n_{0}=N_{\mathbf{b}}(F,L)$ and $$P_{1}(\eta)= (2(\lambda_{2}^{\mathbf{b}}(\eta))^{N}(\lambda_{1}^{\mathbf{b}}(\eta))^{-N})^{q(\eta)}
 (\lambda_{2}^{\mathbf{b}}(\eta))^{N}>1.$$

\noindent\textbf{Sufficiency.} Suppose that for each
$\eta\in(0,\beta]$ there exist $n_{0}=n_{0}(\eta)\in \mathbb{Z}_{+}$ and
$P_{1}=P_{1}(\eta)\geq 1$ that for every $z\in \mathbb{B}_n$ and for every $t_{0}\in
 S_z$ there exists
$k_{0}=k_{0}(t_{0},z)\in \mathbb{Z}_{+},$\ $0\leq k_{0}\leq
n_{0},$ for which inequality (\ref{poch1}) holds. But $\eta$ is arbitrary in $(0,\beta]$ and $\beta>1$ then we can pick $\eta>1.$ 
We select 
$j_{0}\in\mathbb{N}$ satisfying $P_{1}\leq\eta^{j_{0}}.$ For
given $z\in\mathbb{B}_n,$ $t_{0}\in S_z,$ suiting
$k_{0}=k_{0}(t_{0},z)$ and $j\geq j_{0}$ by Cauchy formula
for $F(z+t\mathbf{b})$ as a function of one variable $t$
$$\frac{\partial^{k_{0}+j} F(z+t_{0}\mathbf{b})}{\partial\mathbf{b}^{k_{0}+j}}=\frac{j!}{2\pi i}\int_{|t-t_{0}|=\eta/L(z+t_{0}\mathbf{b})}\frac{1}{(t-t_{0})^{j+1}}\frac{\partial^{k_{0}}F(z+t\mathbf{b})}{\partial\mathbf{b}^{k_{0}}}dt.$$

Since (\ref{poch1}) we have
\begin{gather*}
 \frac{1}{j!}\left|\frac{\partial^{k_{0}+j} F(z+t_{0}\mathbf{b})}{\partial\mathbf{b}^{k_{0}+j}}\right|\leq \frac{L^{j}(z+t_{0}\mathbf{b})}{\eta^{j}}\max\left\{\left|\frac{\partial^{k_{0}}F(z+t\mathbf{b})}{\partial\mathbf{b}^{k_{0}}}
\right|: \right.\\ \left. |t-t_{0}|=\frac{\eta}{L(z+t_{0}\mathbf{b})}\right\}
\leq P_{1}\frac{L^{j}(z+t_{0}\mathbf{b})}{\eta^{j}}\left|\frac{\partial^{k_{0}}F(z+t_0\mathbf{b})}{\partial\mathbf{b}^{k_{0}}}
\right|,
\end{gather*}
that is
\begin{gather*}
 \frac{1!}{(k_{0}+j)!L^{k_{0}+j}(z+t_0\mathbf{b})}\left|\frac{\partial^{k_{0}+j} F(z+t_{0}\mathbf{b})}{\partial\mathbf{b}^{k_{0}+j}}\right|\leq \frac{j!k_{0}!}{(j+k_{0})!}\frac{P_{1}}{\eta^{j}} \times\\ \times
\frac{1}{k_{0}!L^{k_{0}}(z+t_{0}\mathbf{b})}\left|\frac{\partial^{k_{0}} F(z+t_{0}\mathbf{b})}{\partial\mathbf{b}^{k_{0}}}\right|
\leq \eta^{j_{0}-j}\frac{1}{k_{0}!L^{k_{0}}(z+t_{0}\mathbf{b})} \times\\ \times \left|\frac{\partial^{k_{0}} F(z+t_{0}\mathbf{b})}{\partial\mathbf{b}^{k_{0}}}\right|\leq \frac{1}{k_{0}!L^{k_{0}}(z+t_{0}\mathbf{b})}\left|\frac{\partial^{k_{0}} F(z+t_{0}\mathbf{b})}{\partial\mathbf{b}^{k_{0}}}\right|
\end{gather*}
for all $j\geq j_{0}.$

In the above inequality $k_{0}\leq n_{0},$ $n_{0}=n_{0}(\eta)$ and $j_{0}=j_{0}(\eta)$ are independent of $z$ and $t_{0}.$ 
Since  $z\in\mathbb{B}_n$ and $t_{0}\in\mathbb{B}_z$ are arbitrary, this inequality means that  function $F$ is of bounded $L$-index in  direction $\mathbf{b}$ and $N_{\mathbf{b}}(F,L)\leq n_{0}+j_{0}.$
Theorem \ref{te2} is proved.
\end{proof}

\begin{theorem}\label{te20}
Let $\beta>1,$ $L\in{Q}_{\mathbf{b},\beta}(\mathbb{B}_n)$, $\frac{1}{\beta}<\theta_1\leq \theta_2<+\infty,$
$\theta_1L(z)\leq L^*(z)\leq \theta_2 L(z)$. Analytic in $\mathbb{B}_n$ function $F(z),$\ $z\in\mathbb{C}^{n},$ is of
bounded $L^{*}$-index in  direction $\mathbf{b}$ if and only if $F$ is of bounded $L$-index in direction $\mathbf{b}$.
\end{theorem}
\begin{proof}
Obviously, if $L\in{Q}_{\mathbf{b},\beta}(\mathbb{B}_n)$ and $\theta_1L(z)\leq L^*(z)\leq \theta_2 L(z)$, then
$L^{*}\in{Q}_{\mathbf{b},\beta^*}(\mathbb{B}_n),$ $\beta^*\in[\theta_1\beta;\theta_2\beta]$ and $\beta^*>1.$ Let
$N_{\mathbf{b}}(F,L^{*})<+\infty.$ Therefore,  by Theorem \ref{te2} for each
 $\eta^{*},$ $0<\eta^*<\beta\theta_2,$ there exist $n_{0}(\eta^*)\in\mathbb{Z}_{+}$ and
$P_{1}(\eta^*)\geq 1$ that for every $z\in\mathbb{B}_n,$ 
$t_{0}\in S_z$ and some $k_{0},$ $0\leq k_{0}\leq n_{0},$
the inequality (\ref{poch1}) is valid with $L^{*}$ and $\eta^{*}$ instead of
 $L$ and $\eta$. Hence, we put $\eta{*}=\theta_{2}\eta$ and obtain
\begin{multline*}
P_{1}\left|\frac{\partial^{k_{0}}F(z+t_{0}\mathbf{b})}
{\partial\mathbf{b}^{k_{0}}}\right|\geq\max\left\{
\left|\frac{\partial^{k_{0}}
F(z+t\mathbf{b})}{\partial\mathbf{b}^{k_{0}}}\right|:
|t-t_{0}|\leq
\frac{\eta^{*}}{L^{*}(z+t_{0}\mathbf{b})}\right\}\geq\\
%\geq\max\left\{\left|\frac{\partial^{k_{0}}
%F(z+t\mathbf{b})}{\partial\mathbf{b}^{k_{0}}}\right|:
%|t-t_{0}|\leq
%\frac{\eta^{*}}{\theta_{2}L(z+t_{0}\mathbf{b})}\right\}=
\ge\max\left\{\left|\frac{\partial^{k_{0}}
F(z+t\mathbf{b})}{\partial\mathbf{b}^{k_{0}}}\right|:\
|t-t_{0}|\leq\frac{\eta}{L(z+t_{0}\mathbf{b})}\right\}.
\end{multline*}
Therefore,  by Theorem \ref{te2}, %in view of arbitrary $\eta^{*}$
% (and $\eta$ too)
the function $F(z)$ is of bounded $L$-index in
direction $\mathbf{b}$. The converse assertion is obtained by replacing $L$ on
$L^{*}$.
\end{proof}

\begin{theorem} \label{te21}
Let $\beta>1,$ $L\in{Q}_{\mathbf{b},\beta}(\mathbb{B}_n)$, $m\in\mathbb{C}, m\neq0$. Analytic in $\mathbb{B}_n$
function $F(z)$ is of bounded
$L$-index in direction $\mathbf{b}\in\mathbb{C}^{n}$ if and only if
 $F(z)$ is of bounded $L$-index in  direction $m\mathbf{b}$.
\end{theorem}
\begin{proof}
Let $F(z)$ be an analytic in $\mathbb{B}_n$ function of bounded $L$-index in direction
$\mathbf{b}$. By Theorem~\ref{te2} $(\forall \eta>0)$\ $(\exists n_{0}(\eta)\in
 \mathbb{Z}_{+})$\
$(\exists P_{1}(\eta)\geq 1)$ \ $(\forall z\in\mathbb{B}_n)$\  $(\forall t_{0}\in  S_z)$\
$(\exists
k_{0}=k_{0}(t_{0},z)\in\mathbb{Z}_{+},$\ $0\leq k_{0}\leq n_{0}),$
and the following inequality is valid
\begin{equation} \label{poch21}
\max\left\{\left|\frac{\partial^{k_{0}} F(z+t\mathbf{b})}
{\partial \mathbf{b}^{k_{0}}}\right|:
|t-t_{0}|\leq\frac{\eta}{L(z+t_{0} \mathbf{b})}\right\}\leq
P_{1}\left|\frac{\partial^{k_{0}}F(z+t_{0}\mathbf{b})} {\partial
\mathbf{b}^{k_{0}}}\right|.
\end{equation}
Since $\frac{\partial^{k}F}{\partial (m\mathbf{b})^{k}}=
({m})^{k}\frac{\partial^{k}F}{\partial\mathbf{b}^{k}},$\ the inequality (\ref{poch21})
is equivalent to the inequality
\[\max\left\{|{m}|^{k_{0}}\left|\frac{\partial^{k_{0}} F(z+t\mathbf{b})}
{\partial \mathbf{b}^{k_{0}}}\right|:
|t-t_{0}|\leq\frac{\eta}{L(z+t_{0}\mathbf{b})}\right\}\leq
P_{1}|{m}|^{k_{0}}\left|\frac{\partial^{k_{0}}F(z+t_{0}\mathbf{b})}{\partial
 \mathbf{b}^{k_{0}}}\right|\]
or
\[\max\left\{\left|\frac{\partial^{k_{0}}
F(z+\frac{t}{{m}} {m\mathbf{b}})}{\partial
(m\mathbf{b})^{k_{0}}}\right|:
\left|\frac{t-t_{0}}{{m}}\right|\leq\frac{\eta}{|{m}|
L(z+\frac{t_{0}}{{m}}{m\mathbf{b}})}\right\}\!\leq\!
P_{1}\left|\frac{\partial^{k_{0}}F(z+\frac{t_{0}}{{m}}
{m\mathbf{b}})}{\partial (m\mathbf{b})^{k_{0}}}\right|.\]
Denoting $t^{*}=\frac{t}{{m}}, t^{*}_{0}=\frac{t_{0}}{{m}}$,
$\eta^{*}=\frac{\eta}{|{m}|},$ we obtain
$$
\max\left\{\left|\frac{\partial^{k_{0}} F(z+t^{*}{m\mathbf{b}})}
{\partial (m\mathbf{b})^{k_{0}}}\right|: |t^{*}-t^{*}_{0}|
\leq\frac{\eta^{*}}{L(z+t^{*}_{0}{m\mathbf{b}})}\right\}\leq
P_{1}\left| \frac{\partial^{k_{0}}F(z+t_{0}\mathbf{b})} {\partial
\mathbf{b}^{k_{0}}} \right|. $$ By Theorem \ref{te2}  a function
$F(z)$ is of bounded $L$-index in  direction $\mathbf{b}$.
Similarly, the converse assertion can be proved.
\end{proof}

% \begin{center}
\section{Estimate of maximum modulus on a larger circle by maximum modulus
on  a smaller circle and by minimum modulus.} 
% \end{center}  

% {\bf {$5^{0}.$ }}
Now we consider a behaviour of analytic in the unit ball functions
of bounded $L$-index in direction.
Using Theorem \ref{te2}, we prove a criterion of $L$-index  boundedness
in direction.
\begin{theorem} \label{te3}
Let $\beta>1,$ $L\in {Q}_{\mathbf{b},\beta}(\mathbb{B}_n)$. Analytic in $\mathbb{B}_n$ function
$F(z)$ is of bounded $L$-index in  direction
$\mathbf{b}\in\mathbb{C}^{n}$ if and only if for any
 $r_{1}$ and any $r_{2}$ with $0<r_{1}<r_{2}\leq \beta,$
there exists number $P_{1}=P_{1}(r_{1},r_{2})\geq 1$ such that for each
$z^{0}\in \mathbb{B}_n$ and each $t_{0}\in  S_{z^0}$
\begin{gather}  \max\big\{|F(z^{0}+t\mathbf{b})|:
|t-t_{0}|=\frac{r_{2}}{L(z^{0}+t_{0}\mathbf{b})}\big\}\leq \nonumber \\ \label{riv31} \leq P_{1}
\max\big\{|F(z^{0}+t\mathbf{b})|:\
|t-t_{0}|=\frac{r_{1}}{L(z_{0}+t_{0}\mathbf{b})}\big\}.
\end{gather}
\end{theorem}
\begin{proof}

\textbf{Necessity.} Let $N_{\mathbf{b}}(F,L)<+\infty.$
On the contrary, we assume there exists numbers $r_{1}$ and
$r_{2},$ $0<r_{1}<r_{2}\leq\beta,$ that for every $P_{*}\geq 1$
there exist $z^{*}=z^{*}(P_{*})\in \mathbb{B}_n$ and
$t^{*}=t^{*}(P^{*})\in  S_{z^*},$\
the following inequality is valid
\begin{gather*}
 \max\big\{|F(z^{*}+t\mathbf{b})|:
|t-t^{*}|=\frac{r_{2}}{L(z^{*}+t^{*}\mathbf{b})}\big\}>\\ >P_{*}\max
\big\{|F(z^{*}+t\mathbf{b})|:\
|t-t^{*}|=\frac{r_{1}}{L(z^{*}+t^{*}\mathbf{b})}\big\}.  \label{riv32}
 \end{gather*}

By Theorem \ref{te2} there exist $n_{0}=n_{0}(r_{2})\in
\mathbb{Z}_{+}$ and $P_{0}=P_{0}(r_{2})\geq 1$ that for every
$z^{*}\in \mathbb{B}_n,$ $t^{*}\in  S_{z^*}$  and some
$k_{0}=k_{0}(t^{*},z^{*})\in \mathbb{Z}_{+},$  $0\leq k_{0}\leq
n_{0},$ the following inequality holds
\begin{equation} \label{riv33} \max\Big\{\Big|\frac{\partial^{k_{0}}
F(z^{*}+t\mathbf{b})}{\partial\mathbf{b}^{k_{0}}}\Big|:
 |t-t^{*}|=\frac{r_{2}}{L(z^{*}+t^{*}\mathbf{b})}\Big\}\leq P_{0}
\Big|\frac{\partial^{k_{0}}F(z^{*}+t^{*}\mathbf{b})}
{\partial\mathbf{b}^{k_{0}}}\Big|.
\end{equation}
We remark that for  $k_{0}=0$ the proof of necessity is obvious because
(\ref{riv33}) implies
$\max\big\{|F(z^{*}+t\mathbf{b})|:\
|t-t^{*}|={r_{2}}/{L(z^{*}+t^{*}\mathbf{b})}\big\}\leq$\
$P_{0}|F(z^{*}+t^{*}\mathbf{b})|\leq$\
$P_{0}\max\big\{|F(z^{*}+t\mathbf{b})|:\ |t-t^{*}|={r_{1}}/
{L(z^{*}+t^{*}\mathbf{b})}\big\}.$

\noindent We assume that $k_{0}>0,$ and let
\begin{equation}
\label{riv34}
P_{*}=n_{0}!\left(\frac{r_{2}}{r_{1}}\right)^{n_{0}}\left(P_{0}+\frac{r_{1}}{r_{2}-r_{1}}\right)+1.
\end{equation}
Let
$t_{0}\in S_{z^*}$\ be such that
$|t_{0}-t^{*}|={r_{1}}/{L(z^{*}+t^{*}\mathbf{b})}$ and
$$ |F(z^{*}+t_{0}\mathbf{b})|=\max\left\{|F(z^{*}+t\mathbf{b})|:
|t-t^{*}|={r_{1}}/
{L(z^{*}+t^{*}\mathbf{b})}\right\}>0,$$
but $t_{0j}\in S_{z^*},$\
$|t_{0j}-t^{*}|={r_{2}}/{L(z^{*}+t^{*}\mathbf{b})},$\
be such that
$$\left|\frac{\partial^{j}F(z^{*}+t_{0j}\mathbf{b})}{\partial\mathbf{b}^{j}}\right|=\max\left\{\left|\frac{\partial^{j}F(z^{*}+t\mathbf{b})}{\partial\mathbf{b}^{j}}\right|:
|t-t^{*}|={r_{2}}/{L(z^{*}+t^{*}\mathbf{b})}\right\}, \
 j\in\mathbb{Z}_{+}.$$ \\  We remark that in the case
$|F(z^{*}+t_{0}\mathbf{b})|=0$ by the uniqueness theorem for all
$t\in S_{z^*}$ an equality $F(z^{*}+t\mathbf{b})=0$ can be obtained.
However, it contradicts an inequality (\ref{riv32}). By Cauchy inequality we have
\begin{equation} \label{riv35}
\frac{1}{j!}\left|\frac{\partial^{j}F(z^{*}+t^{*}\mathbf{b})}
{\partial\mathbf{b}^{j}}\right|\leq
\left(\frac{L(z^{*}+t^{*}\mathbf{b})}{r_{1}}\right)^{j}|F(z^{*}+t_{0}\mathbf{b})|,
j\in \mathbb{Z}_{+}
\end{equation} and
\begin{eqnarray} \label{riv36}
\left|\frac{\partial^{j}F(z^{*}+t_{0j}\mathbf{b})}{\partial\mathbf{b}^{j}}-
\frac{\partial^{j}F(z^{*}+t^{*}\mathbf{b})}{\partial\mathbf{b}^{j}}\right|=
\left|\int_{t^{*}}^{t_{0j}}\frac{\partial^{j+1}F(z^{*}
+t\mathbf{b})}{\partial\mathbf{b}^{j+1}}\,dt\right|\leq\nonumber\\
\leq\left|\frac{\partial^{j+1}F(z^{*}+t_{0(j+1)}\mathbf{b})}
{\partial\mathbf{b}^{j+1}}\right|\frac{r_{2}}{L(z^{*}+t^{*}\mathbf{b})}.
\end{eqnarray}
The inequalities (\ref{riv35}) and (\ref{riv36}) imply that \begin{gather*}
\left|\frac{\partial^{j\!+\!1}F(z^{*}\!+\!t_{0(j\!+\!1)}\mathbf{b})}
{\partial\mathbf{b}^{j\!+\!1}}\right|\geq\frac{L(z^{*}\!+\!t^{*}\mathbf{b})}{r_{2}}
\left\{\left|\frac{\partial^{j}F(z^{*}\!+\!t_{0j}\mathbf{b})}{\partial\mathbf{b}^{j}}\right|
\!-\!\left|\frac{\partial^{j}F(z^{*}\!+\!t^{*}\mathbf{b})}{\partial\mathbf{b}^{j}}\right|\right\}\geq
% =\\
% =\frac{L(z^{*}+t^{*}\mathbf{b})}{r_{2}}\left|\frac{\partial^{j}F(z^{*}+t_{0j}\mathbf{b})}
% {\partial\mathbf{b}^{j}}\right|-\frac{L(z^{*}+t^{*}\mathbf{b})}
% {r_{2}}\left|\frac{\partial^{j}F(z^{*}+t^{*}\mathbf{b})}
% {\partial\mathbf{b}^{j}}\right|\geq
\\ \geq
\frac{L(z^{*}\!+\!t^{*}\mathbf{b})}{r_{2}}
\left|\frac{\partial^{j}F(z^{*}+t_{0j}\mathbf{b})}{\partial\mathbf{b}^{j}}\right|
-\frac{j!L^{j+1}(z^{*}+t^{*}\mathbf{b})}{r_{2}(r_{1})^{j}}|F(z^{*}+t_{0}\mathbf{b})|,
j\in \mathbb{Z}_{+}.
\end{gather*}
Hence, for $k_{0}\geq 1$ we get
\begin{gather}
\nonumber\left|\frac{\partial^{k_{0}}F(z^{*}+t_{0k_{0}}\mathbf{b})}
{\partial\mathbf{b}^{k_{0}}}\right|\geq\frac{L(z^{*}+t^{*}\mathbf{b})}{r_{2}}
\left|\frac{\partial^{k_{0}-1}F(z^{*}+t_{0(k_{0}-1)}\mathbf{b})}
{\partial\mathbf{b}^{k_{0}-1}}\right|-\!\\ \nonumber
\!-\frac{(k_{0}-1)!L^{k_{0}}(z^{*}+t^{*}\mathbf{b})}{r_{2}
(r_{1})^{k_{0}-1}}|F(z^{*}+t_{0}\mathbf{b})|\geq
\ldots\geq\frac{L^{k_{0}}(z^{*}+t^{*}\mathbf{b})}
{(r_{2})^{k_{0}}}|F(z^{*}+t_{00}\mathbf{b})|-
\\ \nonumber -\left(\frac{0!}{(r_{2})^{k_{0}}}+\frac{1!}{(r_{2})^{k_{0}-1}r_{1}}+
\ldots+\frac{(k_{0}-1)!}{r_{2}(r_{1})^{k_{0}-1}}\right)
L^{k_{0}}(z^{*}+t^{*}\mathbf{b})\times\!\\
\times|F(z^{*}+t_{0}\mathbf{b})|
\!=\frac{L^{k_{0}}(z^{*}+t^{*}\mathbf{b})}{(r_{2})^{k_{0}}}|F(z^{*}+t_{0}\mathbf{b})|
\left(\frac{|F(z^{*}+t_{00}\mathbf{b})|}{|F(z^{*}+t_{0}\mathbf{b})|}
-\!\sum_{j=0}^{k_{0}-1}j!\left(\frac{r_{2}}{r_{1}}\right)^{j}\right).\label{riv37}
\end{gather}
Since (\ref{riv32}) we have 
${|F(z^{*}+t_{00}\mathbf{b})|}/{|F(z^{*}+t_{0}\mathbf{b})|}>P_{*}.$ 
Besides, this inequality holds 
\begin{gather*}
\sum_{j=0}^{k_{0}-1}j!\left(\frac{r_{2}}{r_{1}}\right)^{j}\leq
k_{0}!\left(\frac{\left(r_{2}/r_{1}\right)^{k_{0}}-1}
{r_{2}/r_{1}-1}\right)\leq
%k_{0}!\left(\frac{r_{2}}{r_{1}}\right)^{k_{0}}\frac{1}{r_{2}/r_{1}-1}\leq\\
%\leq
n_{0}!\frac{r_{1}}{r_{2}-r_{1}}\left(\frac{r_{2}}{r_{1}}\right)^{n_{0}}.
\end{gather*}
Applying (\ref{riv34}), we obtain
\begin{eqnarray*} \frac{|F(z^{*}\!+\!t_{00}\mathbf{b})|}
{|F(z^{*}\!+\!t_{0}\mathbf{b})|}\!-\!\sum_{j=0}^{k_{0}\!-\!1}j!
\left(\frac{r_{2}}{r_{1}}\right)^{j}\!>\!P_{*}\!-\!n_{0}!\frac{r_{1}}{r_{2}-r_{1}}
\left(\frac{r_{2}}{r_{1}}\right)^{n_{0}}
=n_{0}!\left(\frac{r_{2}}{r_{1}}\right)^{n_{0}}P_{0}+1.
\end{eqnarray*}
From (\ref{riv37}), in view of (\ref{riv33}) and (\ref{riv35}),
it follows that
\begin{gather*}
\left|\frac{\partial^{k_{0}}F(z^{*}\!+\!t_{0k_{0}}\mathbf{b})}
{\partial\mathbf{b}^{k_{0}}}\right|\!>\!\frac{L^{k_{0}}(z^{*}\!+\!t^{*}\mathbf{b})}{(r_{2})^{k_{0}}}\left(P_{*}-n_{0}!
\frac{r_{1}}{r_{2}-r_{1}}\left(\frac{r_{2}}{r_{1}}\right)^{n_{0}}\right)
\left(\frac{r_{1}}{L(z^{*}\!+\!t^{*}\mathbf{b})}\right)^{k_{0}}\!\times\! \\ \times \frac{1}{k_{0}!}\!
\left|\frac{\partial^{k_{0}}F(z^{*}\!+\!t^{*}\mathbf{b})}
{\partial\mathbf{b}^{k_{0}}}\right|\geq\left(\frac{r_{1}}{r_{2
}}\right)^{n_{0}}\frac{1}{n_{0}!P_{0}}\left(P_{*}\!-\!n_{0}!\frac{r_{1}}{r_{2}\!-\!r_{1}}
\left(\frac{r_{2}}{r_{1}}\right)^{n_{0}}\right) 
\left|\frac{\partial^{k_{0}}F(z^{*}\!+\!t_{0k_{0}}\mathbf{b})}
{\partial\mathbf{b}^{k_{0}}}\right|.
\end{gather*}
Hence, $P_{*}<n_{0}!\left(\frac{r_{2}}{r_{1}}\right)^{n_{0}}
\left(P_{0}+\frac{r_{1}}{r_{2}-r_{1}}\right)$ and it contradicts
(\ref{riv34}).

\textbf{Sufficiency.} We choose any two numbers
$r_{1}\in(0,1)$ and $r_{2}\in(1,\beta)$. For given
$z^{0}\in\mathbb{B}_n,$ $t_{0}\in S_{z^0}$ we expand a function
$F(z^{0}+t\mathbf{b})$ in the power series by powers $t-t_{0}$
\begin{equation*}
 F(z^{0}+t\mathbf{b})=\sum_{m=0}^{\infty}b_{m}(z^{0}+t_0\mathbf{b})(t-t_{0})^{m}, \
b_{m}(z^{0}+t_0\mathbf{b})=\frac{1}{m!}\frac{\partial^{m} F(z^{0}+t_{0}\mathbf{b})}{\partial\mathbf{b}^{m}}
\end{equation*}
in a disc $\left\{t: |t-t_{0}|\leq\displaystyle\frac{\mathstrut \beta}{L(z^{0}+t_{0}\mathbf{b}}\right\}\subset  S_{z^0}.$ For $r\leq\displaystyle\frac{\mathstrut \beta}{L(z^{0}+t_{0}\mathbf{b})}$ we denote
$$M_{\mathbf{b}}(r,z^{0},t_{0},F)=\max\{|F(z^{0}+t\mathbf{b})|: |t-t_{0}|=r\},$$
 $$\mu_{\mathbf{b}}(r,z^{0},t_{0},F)=\max\{|b_{m}(z^{0}+t_0\mathbf{b})|r^{m}: m\geq 0\},$$
$$\nu_{\mathbf{b}}(r,z^{0},t_{0},F)=\max\{|b_{m}(z^{0})|r^{m}: |b_{m}(z^{0}+t_0\mathbf{b})|r^{m}=\mu_{\mathbf{b}}(r,z^{0},t_{0},F)\}.$$

By Cauchy inequality $\mu_{\mathbf{b}}\!(r,z^{0},t_{0},F)\!\leq\! M_{\mathbf{b}}(r,z^{0},t_{0},F).$
 But for $\displaystyle r=\frac{1}{L(z^{0}+t_{0}\mathbf{b})}$ we have
\begin{gather*}
M_{\mathbf{b}}(r_{1}r,z^{0},t_{0},F)\leq \sum_{m=0}^{\infty}|b_{m}(z^{0}+t_0\mathbf{b})|r^{m}r_{1}^{m}\leq \mu_{\mathbf{b}}(r,z^{0},t_{0},F)\sum_{m=0}^{\infty}r_{1}^{m}=\\ =\frac{1}{1-r_{1}}
\mu_{\mathbf{b}}(r,z^{0},t_{0},F)\end{gather*}
and, applying a monotone of $\nu_{\mathbf{b}}(r,z^{0},t_{0},F)$ by $r,$ we get  
\begin{gather*}
 \ln\mu_{\mathbf{b}}(r_{2}r,z^{0},t_{0},F)-\ln\mu_{\mathbf{b}}(r,z^{0},t_{0},F)=
\int_{r}^{r_{2}r}\frac{\nu_{\mathbf{b}}(t,z^{0},t_{0},F)}{t}dt 
\geq \nu_{\mathbf{b}}(r,z^{0},t_{0},F)\ln r_{2}.
\end{gather*}

Hence, we get 
\begin{gather} \nonumber
\nu_{\mathbf{b}}(r,z^{0},t_{0},F)\leq\frac{1}{\ln r_{2}}( \ln\mu_{\mathbf{b}}(r_{2}r,z^{0},t_{0},F)-\ln\mu_{\mathbf{b}}(r,z^{0},t_{0},F))\leq\\ \nonumber\leq\frac{1}{\ln r_{2}}\{\ln M_{\mathbf{b}}(r_{2}r,z^{0},t_{0},F)-\ln((1-r_{1})M_{\mathbf{b}}(r_{1}r,z^{0},t_{0},F))\}=\\ =-\frac{\ln(1-r_{1})}{\ln r_{2}}+\frac{1}{\ln r_2}\{\ln M_{\mathbf{b}}(r_{2}r,z^{0},t_{0},F)-\ln M_{\mathbf{b}}(r_{1}r,z^{0},t_{0},F))\}\label{eq_th_12_19}
\end{gather}

Let $N_{\mathbf{b}}(z^{0}+t_{0}\mathbf{b},L,F)$ be $L$-index in direction of function $F$ at a point $z^{0}+t_{0}\mathbf{b},$ i.e.  $N_{\mathbf{b}}(z^{0}+t_{0}\mathbf{b},L,F)$ is the smallest number $m_{0}$ for which  an inequality (\ref{eq1}) holds with $z=z^{0}+t_{0}\mathbf{b}.$
It is obvious that $$N_{\mathbf{b}}(z^{0}+t_{0}\mathbf{b},L,F)\leq \nu_{\mathbf{b}}(1/L(z^{0}+t^{0}\mathbf{b}),z^{0},t_{0},F)=\nu_{\mathbf{b}}(r,z^{0},t_{0},F).$$
However, an inequality (\ref{riv31}) can be written in the following form
$$M_{\mathbf{b}}\left(\frac{r_{2}}{L(z^{0}+t_{0}\mathbf{b})},z^{0},t_{0},F\right)\leq P_{1}(r_{1},r_{2}) M_{\mathbf{b}}\left(\frac{r_{1}}{L(z^{0}+t_{0}\mathbf{b})},z^{0},t_{0},F\right).$$
Thus, from (\ref{eq_th_12_19}) we have 
$N_{\mathbf{b}}(z^{0}+t_{0}\mathbf{b},L,F)\leq -\frac{\ln(1-r_{1})}{\ln r_{2}}+\frac{\ln P_{1}(r_{1},r_{2})}{\ln r_{2}} $
for every $z^{0}\in\mathbb{C}^{n},$ $t_{0}\in\mathbb{C},$ i.e.
$$N_{\mathbf{b}}(F,L)\leq -\frac{\ln(1-r_{1})}{\ln r_{2}}+\frac{\ln P_{1}(r_{1},r_{2})}{\ln r_{2}}.$$
Theorem \ref{te3} is proved.
\end{proof}

In view of proof of Theorem \ref{te3} the following theorem is true.
\begin{theorem} \label{te4}
Let $\beta>1$ and $L\in {Q}_{\mathbf{b},\beta}(\mathbb{B}_n).$ Analytic in $\mathbb{B}_n$ function
$F(z)$ is of bounded $L$-index in direction
$\mathbf{b}\in\mathbb{C}^{n}$ if and only if there exist numbers
 $r_{1}$ and $r_{2}$, $0<r_{1}<1<r_{2}\leq \beta,$ and $P_{1}\geq
1$ that for every $z^{0}\in \mathbb{B}_n$ and $t_{0}\in
 S_{z^0}$ inequality (\ref{riv31}) holds.
\end{theorem}

Here is another criterion that is an analogue of Hayman Theorem \cite{Hayman}.
\begin{theorem}\label{te11}
Let $\beta>1$ and $L\in{Q}_{\mathbf{b},\beta}(\mathbb{B}_n)$. An analytic in $\mathbb{B}_n$ function $F(z)$ is
of bounded $L$-index in
direction $\mathbf{b}\in\mathbb{C}^{n}$ if and only if there
exist $p\in\mathbb{Z}_{+}$ and $C>0$ such that for every
$z\in\mathbb{B}_n$ the following inequality holds
\begin{equation} \label{riv101}
\left|\frac{1}{L^{p+1}(z)}\frac{\partial^{p+1}
F(z)}{\partial\mathbf{b}^{p+1}}\right| \leq
C\max\left\{\left|\frac{1}{L^{k}(z)}\frac{\partial^{k}F(z)}
{\partial\mathbf{b}^{k}}\right|:\  0\leq k\leq p\right\}.
\end{equation}
\end{theorem}
\begin{proof}
\textbf{Necessity.} If $N_{\mathbf{b}}(F,L)<+\infty$ then
 by definition of $L$-index boundedness in the direction 
we obtain an inequality (\ref{riv101}) with $p=N_{\mathbf{b}}(F,L)$ and
$C=(N_{\mathbf{b}}(F,L)+1)!$

\textbf{Sufficiency.} Let an inequality (\ref{riv101}) holds,
$z^{0}\in\mathbb{B}_n,$ $t_{0}\in S_{z^0}$ and
$$K=\left\{t\in\mathbb{C}: |t-t_{0}|\leq
\frac{1}{L(z^{0}+t_{0}\mathbf{b})}\right\}.$$ Thus, using $L\in Q_{\mathbf{b},\beta}(\mathbb{B}_n),$ for every
 $t\in K$ with (\ref{riv101}) we have 
\begin{gather}
 \!\frac{1}{L^{p+1}(z^{0}+t_{0}\mathbf{b})}\left|\frac{\partial^{p+1}
F(z^{0}+t\mathbf{b})}{\partial\mathbf{b}^{p+1}}\right| \!\leq\!\left(\frac{L(z^{0}+t\mathbf{b})}{L(z^{0}+t_{0}\mathbf{b})}\right)^{p+1} \frac{1}{L^{p+1}(z^{0}+t\mathbf{b})}\!\times\!\nonumber\\ \times
 \left|\frac{\partial^{p+1}
 F(z^{0}+t\mathbf{b})}{\partial\mathbf{b}^{p+1}}\right|
 \leq (\lambda_{2}^{\mathbf{b}}(1))^{p+1}\frac{1}{L^{p+1}(z^{0}+t\mathbf{b})}\left|\frac{\partial^{p+1}
 F(z^{0}+t\mathbf{b})}{\partial\mathbf{b}^{p+1}}\right|\leq\nonumber\\ \leq C(\lambda_{2}^{\mathbf{b}}(1))^{p+1}
 \max\left\{\left|\frac{1}{L^{k}(z^{0}+t\mathbf{b})}\frac{\partial^{k}F(z^{0}+t\mathbf{b})}
 {\partial\mathbf{b}^{k}}\right|:\  0\leq k\leq p\right\}\leq\nonumber\\ \leq C(\lambda_{2}^{\mathbf{b}}(1))^{p+1}
 \max\left\{\left(\frac{L(z^{0}+t_{0}\mathbf{b})}{L(z^{0}+t\mathbf{b})}\right)^{k}\left|\frac{1}{L^{k}(z^{0}+t_{0}\mathbf{b})}\frac{\partial^{k}F(z^{0}+t\mathbf{b})}
 {\partial\mathbf{b}^{k}}\right|: \right.\nonumber\\ \left. 0\leq k\leq p\right\}
 \leq C(\lambda_{2}^{\mathbf{b}}(1))^{p+1}
 \max\left\{ \left|\frac{1}{L^{k}(z^{0}+t_{0}\mathbf{b})}\frac{\partial^{k}F(z^{0}+t\mathbf{b})}
 {\partial\mathbf{b}^{k}}\right|\times\right.\nonumber \\ \left.\times (\lambda_{1}^{\mathbf{b}}(1))^{-k} :  0\leq k\leq p\right\}\leq Bg_{z^{0}}(t_{0},t),
\label{eq_th_1_24}
\end{gather}
where $B=C(\lambda_{2}^{\mathbf{b}}(1))^{p+1}(\lambda_{1}^{\mathbf{b}}(1))^{-p}$ and $$g_{z^{0}}(t_{0},t)\!=\!\max\!\left\{\left|\displaystyle\frac{\mathstrut 1}{L^{k}(z^{0}+t_{0}\mathbf{b})}\frac{\partial^{k}F(z^{0}+t\mathbf{b})}
{\partial\mathbf{b}^{k}}\right|\!: 0\leq k\leq p\right\}.$$

We introduce denotations
$$\gamma_{1}=\left\{t\in\mathbb{C}: |t-t_{0}|=\frac{1}{2\beta L(z^{0}\!+\!t_{0}\mathbf{b})}\right\}, \
\gamma_{2}=\left\{t\in\mathbb{C}: |t-t_{0}|=\frac{\beta}{L(z^{0}\!+\!t_{0}\mathbf{b})}\right\}. $$
We choose arbitrary points $t_{1}\in\gamma_{1},$ $t_{2}\in\gamma_{2}$ and join them by a piecewise-analytic curve $\gamma = (t=t(s), 0\leq s\leq T),$ that $g_{z^{0}}(t_{0},t)\neq 0$  with $t\in\gamma.$
We choose a curve $\gamma$ that its length $|\gamma|$ does not exceed $\displaystyle\frac{2\beta^2+1}{\beta L(z^{0}+t_{0}\mathbf{b})}.$

The function $g_{z^{0}}(t_{0},t(s))$ is continuous on $[0,T].$
Without loss of generality, we consider that function
 $t=t(s)$ is analytic on $[0,T].$
 Otherwise, we can consider separately the intervals of analyticity for this function and repeat similar arguments which below we present for $ [0, T].$
 First, we prove that the function $g_{z^{0}}(t_{0},t(s))$ is continuously differentiable on $ [0, T] $ except, perhaps, a finite set of points.
For arbitrary $k_{1}, k_{2}, 0\leq k_{1}\leq k_{2}\leq p,$ either 
$$\frac{1}{L^{k_{1}}(z^{0}+t_{0}\mathbf{b})}\left|\frac{\partial^{k_{1}} F(z^{0}+t(s)\mathbf{b})}{\partial\mathbf{b}^{k_{1}}}\right|\equiv
\frac{1}{L^{k_{2}}(z^{0}+t_{0}\mathbf{b})}\left|\frac{\partial^{k_{2}} F(z^{0}+t(s)\mathbf{b})}{\partial\mathbf{b}^{k_{2}}}\right| $$
or the equality
$$\frac{1}{L^{k_{1}}(z^{0}+t_{0}\mathbf{b})}\left|\frac{\partial^{k_{1}} F(z^{0}+t(s)\mathbf{b})}{\partial\mathbf{b}^{k_{1}}}\right|=
\frac{1}{L^{k_{2}}(z^{0}+t_{0}\mathbf{b})}\left|\frac{\partial^{k_{2}} F(z^{0}+t(s)\mathbf{b})}{\partial\mathbf{b}^{k_{2}}}\right| $$
holds only for a finite set of points $s_{k}\in[0,T].$ Thus, we can split the segment $[0,T]$ on a finite number of segments that on each segment $$g_{z^{0}}(t_{0},t(s))\equiv\frac{1}{L^{k}(z^{0}+t_{0}\mathbf{b})} \left|\frac{\partial^{k}F(z^{0}+t(z)\mathbf{b})}{\partial\mathbf{b}^{k}}\right|$$
for some $k,$ $0\leq k\leq p.$ This means that a function $g_{z^{0}}(t_{0},t(s))$ is continuously differentiable except, perhaps, a finite set of points. Since (\ref{eq_th_1_24}) we obtain
\begin{gather*}
 \frac{dg_{z^{0}}(t_{0},t(s))}{ds} \leq\max\left\{\frac{d}{ds}\left(\frac{1}{L^{k}(z^{0}+t_{0}\mathbf{b})}
\left|\frac{\partial^{k}F(z^{0}+t(s)\mathbf{b})}{\partial\mathbf{b}^{k}}\right|\right):  0\leq k\leq p\right\}\leq\\
\leq \max\left\{\frac{1}{L^{k}(z^{0}+t_{0}\mathbf{b})}\left|\frac{\partial^{k+1}
F(z^{0}+t(s)\mathbf{b})}{\partial\mathbf{b}^{k+1}}\right| |t'(s)|: 0\leq k\leq p \right\}=\\
=L(z^{0}+t_{0}\mathbf{b})|t'(s)|\max\left\{\frac{1}{L^{k+1}(z^{0}+t_{0}\mathbf{b})}\left|\frac{\partial^{k+1}
F(z^{0}+t(s)\mathbf{b})}{\partial\mathbf{b}^{k+1}}\right|:  0\leq k\leq p \right\} \leq\\ \leq
 Bg_{z^{0}}(t_{0},t(s))|t'(s)|L(z^{0}+t_{0}\mathbf{b}).
\end{gather*}

Hence, we have 
\begin{gather*} \left|\ln\frac{g_{z^{0}}(t_{0},t_{2})}{g_{z^{0}}(t_{0},t_{1})}\right|=
\left|\int_{0}^{T}\frac{dg_{z^{0}}(t_{0},t(s))}{g_{z^{0}}(t_{0},t(s))}\right|\leq B L(z^{0}+t_{0}\mathbf{b})
\int_{0}^{T}|t'(s)|ds=\\ =B L(z^{0}+t_{0}\mathbf{b})|\gamma|\leq 2B \frac{\beta^2+1}{\beta}.\end{gather*}

If we pick a point $t_{2}\in\gamma_{2},$ for which
$$|F(z^{0}+t_{2}\mathbf{b})|=\max\left\{|F(z^{0}+t\mathbf{b})|:
 |t-t_{0}|=\frac{\beta}{L(z^{0}+t_{0}\mathbf{b})}\right\},$$
then we have
\begin{gather}
 \max\left\{|F(z^{0}+t\mathbf{b})|:
 |t-t_{0}|=\frac{2}{L(z^{0}+t_{0}\mathbf{b})}\right\}\leq g_{z^{0}}(t_{0},t_{2}) %\nonumber\leq\\
 \leq g_{z^{0}}(t_{0},t_{1})\exp\{2B \frac{\beta^2+1}{\beta}\}. \label{eq_th_125}\end{gather}

Applying Cauchy inequality and using $t_{1}\in\gamma_{1},$ for all $j=1,\ldots, p$ we have
\begin{gather*}
 \left|\frac{\partial^{j} F(z^{0}\!+\!t_{1}\mathbf{b})}{\partial\mathbf{b}^{j}}\right|\leq%\\ \leq 
 j!(2\beta L(z^{0}\!+\!t_{0}\mathbf{b}))^{j}\max\left\{|F(z^{0}\!+\!t\mathbf{b}):
 |t\!-\!t_{1}|=\frac{1}{2\beta L(z^{0}\!+\!t_{0}\mathbf{b})}\right\}\leq \\ \leq
j!(2\beta L(z^{0}+t_{0}\mathbf{b}))^{j} \max\left\{|F(z^{0}+t\mathbf{b}):
 |t-t_{0}|=\frac{1}{\beta L(z^{0}+t_{0}\mathbf{b})}\right\},
\end{gather*}
i.e.
$$g_{z^{0}}(t_{0},t_{1})\leq p!(2\beta)^{p}\max\left\{|F(z^{0}+t\mathbf{b}):
 |t-t_{0}|=\frac{1}{\beta L(z^{0}+t_{0}\mathbf{b})}\right\} .$$
Thus, (\ref{eq_th_125}) implies
\begin{gather*}
 |F(z^{0}+t_{2}\mathbf{b})|=\max\left\{|F(z^{0}+t\mathbf{b})|:
|t-t_{0}|=\frac{\beta}{L(z^{0}+t_{0}\mathbf{b})}\right\}\leq
g_{z^{0}}(t_{0},t_{2}) \leq \\
\leq g_{z^{0}}(t_{0},t_{1})\exp\bigg\{2B\frac{\beta^2+1}{\beta}\bigg\} \leq p!(2\beta)^{p}\exp\bigg\{2B\frac{\beta^2+1}{\beta}\bigg\} \times \\ \times
\max\left\{|F(z^{0}+t\mathbf{b})|: |t-t_{0}|=\frac{1}{\beta L(z^{0}+t_{0}\mathbf{b})}\right\}.
\end{gather*}
By Theorem  \ref{te4} this inequality  implies that a function $F$ is of bounded $L$-index in the direction $\mathbf{b}\in\mathbb{C}^{n}.$ Theorem \ref{te11} is proved.
\end{proof}

The following theorem gives an estimate of maximum modulus by minimum modulus.
\begin{theorem} \label{te5}
Let  $\beta>1$ and $L\in{Q}_{\mathbf{b},\beta}(\mathbb{B}_n)$. Analytic in $\mathbb{B}_n$ function  $F(z)$ is of bounded
 $L$-index in  direction $\mathbf{b}$
if and only if for every $R,$ $0<R\leq \beta,$ there exist numbers $P_{2}(R)\geq
1$ and $\eta(R)\in(0,R)$ that for each $z^{0}\in
\mathbb{B}_n,$ $t_{0}\in  S_{z^0}$ and some
$r=r(z^{0},t_{0})\in[\eta(R),R]$ the following inequality is valid
 \begin{gather} 
\max\Big\{|F(z^{0}+t\mathbf{b})|: |t-t_{0}|=\frac{r}
{L(z^{0}+t_{0}\mathbf{b})}\Big\}\leq \nonumber \\ \leq 
P_{2}\min\Big\{|F(z^{0}+t\mathbf{b})|:\
|t-t_{0}|=\frac{r}{L(z^{0}+t_{0}\mathbf{b})}\Big\}.\label{riv51}
\end{gather}
\end{theorem}

\begin{proof}

\textbf{Necessity.}
Let $N_{\mathbf{b}}(F,L)=N<+\infty$ and $R\geq 0.$ We put
$$R_{0}=1, r_{0}=\frac{R}{8(R+1)},\
R_{j}=\frac{R_{j-1}}{4N}r_{j-1}^{N},\ r_{j}=\frac{1}{8}R_{j}
(j=1,2,\ldots,N).$$
Let $z^{0}\in\mathbb{B}_n,$ $t_{0}\in S_{z^0}$ and $N_{0}=N_{\mathbf{b}}(z^{0}+t_{0}\mathbf{b},L,F)$ be $L$-index in  direction $\mathbf{b}$ of function $F$ at point $z^{0}+t_{0}\mathbf{b},$ i.e. $N_{\mathbf{b}}(z^{0}+t_{0}\mathbf{b},L,F)$ is the smallest number $m_{0},$ for which   inequality (\ref{eq1}) holds with $z=z^{0}+t_{0}\mathbf{b}.$ The maximum in the right part of \eqref{eq1} is attained at $m_0.$
But $0\leq N_{0}\leq N.$ For given $z^{0}\in\mathbb{B}_n,$ $t_{0}\in S_{z^0}$ a function $F(z^{0}+t\mathbf{b})$
expands in power series by powers $t-t_{0}$
\begin{equation*}
 F(z^{0}+t\mathbf{b})=\sum_{m=0}^{\infty}b_{m}(z^{0}+t_0\mathbf{b})(t-t_{0})^{m}, \ 
%  \end{equation*}
% \begin{equation*}
b_{m}(z^{0}+t_0\mathbf{b})=\frac{1}{m!}\frac{\partial^{m} F(z^{0}+t_{0}\mathbf{b})}{\partial\mathbf{b}^{m}}.
\end{equation*}
We put $$a_{m}(z^{0})=\frac{|b_{m}(z^{0}+t_0\mathbf{b})|}{L^{m}(z^{0})}=\frac{1}{m!L^{m}(z^{0})}
\left|\frac{\partial^{m}F(z^{0}+t_{0}\mathbf{b})}{\partial\mathbf{b}^{m}}\right|.$$
For any $m\in
\mathbb{Z}_{+}$ inequality holds $$a_{N_{0}}(z^{0})\geq
a_{m}(z^{0})=R_{0}a_{m}(z^{0}).$$ There exists the smallest number
$n_{0}\in\{0,1,\ldots,N_{0}\}$ that for all $m\in \mathbb{Z}_{+}$ $a_{n_{0}}(z^{0})\geq
a_{m}(z^{0})R_{N_{0}-n_{0}}.$  Thus, $a_{n_{0}}(z^{0})\geq a_{N_{0}}(z^{0})R_{N_{0}-n_{0}}$ and 
$a_{j}(z^{0})<a_{N_{0}}(z^{0})R_{N_{0}-j}$
for $j<n_{0}$, because if $a_{j_{0}}(z^{0})\geq a_{N_{0}}(z^{0})R_{N_{0}-j_{0}}$
for some $j_{0}<n_{0}$, then $a_{j_{0}}(z^0)\geq a_{m}(z^0)R_{N_{0}-j_{0}}$
for all $m\in \mathbb{Z}_{+}$ and it contradicts the choice of $n_{0}.$
Since inequalities
$a_{j}(z^0)<a_{N_{0}}(z^0)R_{N_{0}-j}$ $(j<n_{0})$ and $a_{m}(z^0)\leq a_{N_{0}}(z^0)$
 $(m>n_{0})$  for $t\in S_{z^0}$ and $|t-t_{0}|=\frac{1}{L(z^{0}+t_{0}\mathbf{b})}r_{N_{0}-n_{0}}$ 
 we have 
\begin{gather}
|F(z^{0}+t\mathbf{b})|=|b_{n_{0}}(z^{0}+t_0\mathbf{b})(t-t_{0})^{n_{0}}+\sum_{m\neq n_{0}}b_{m}(z^{0}+t_0\mathbf{b})(t-t_{0})^{m}|\geq\nonumber\\
\geq|b_{n_{0}}(z^{0})||t-t_{0}|^{n_{0}}
-\sum_{m\neq n_{0}}|b_{m}(z^{0})||t-t_{0}|^{m}=a_{n_{0}}(z^{0})r_{N_{0}-n_{0}}^{n_{0}}-
\nonumber\\
-\sum_{m\neq 0}a_{m}(z^{0})r^{m}_{N_{0}-n_{0}} =a_{n_{0}}(z^{0})r_{N_{0}-n_{0}}^{n_{0}}-\sum_{j<n_{0}}a_{j}(z^{0})r^{j}_{N_{0}-n_{0}}-\nonumber\\ \!-\!\sum_{m>n_{0}}a_{m}(z^{0})r^{m}_{N_{0}-n_{0}}\! \geq\! a_{N_{0}}(z^{0})R_{N_{0}-n_{0}}r_{N_{0}-n_{0}}^{n_{0}}-\sum_{j<n_{0}}a_{N_{0}}(z^{0})R_{N_{0}-j}r^{j}_{N_{0}-n_0}\!-\!\nonumber\\ -
 \sum_{m>n_{0}}a_{N_{0}}(z^{0})r^{m}_{N_{0}-n_{0}}
 \geq
a_{N_{0}}(z^{0})R_{N_{0}-n_{0}}r_{N_{0}-n_{0}}^{n_{0}}-n_{0}a_{N_{0}}(z^{0})R_{N_{0}-n_{0}+1}-
\nonumber\\ -a_{N_{0}}(z^0)r_{N_{0}-n_{0}}^{n_{0}+1}\frac{1}{1-r_{N_{0}-n_{0}}} = a_{N_{0}}(z^{0})
  \left(R_{N_{0}-n_{0}}r_{N_{0}-n_{0}}^{n_{0}}-\frac{n_{0}}{4N}R_{N_{0}-n_{0}}
 r^{N}_{N_{0}-n_{0}}-\nonumber\right.\\ \left. -r_{N_{0}-n_{0}}^{n_{0}}\frac{r_{N_{0}-n_{0}}}{1-r_{N_{0}-n_{0}}}\right)
   \geq a_{N_{0}}(z^{0})
  \left(R_{N_{0}-n_{0}}r_{N_{0}-n_{0}}^{n_{0}}-\frac{1}{4}R_{N_{0}-n_{0}}r_{N_{0}-n_{0}}^{n_{0}}-
 \nonumber\right.\\ \left. -\frac{1}{4}R_{N_{0}-n_{0}}r_{N_{0}-n_{0}}^{n_{0}}\right) =
  \frac{1}{2}a_{N_{0}}(z^{0})R_{N_{0}-n_{0}}r_{N_{0}-n_{0}}^{n_{0}}.
\label{eq_th_1321}
\end{gather}

Besides, for $t\in S_{z^0}$ the following inequality holds 
\begin{gather}
 \nonumber
|F(z^{0}+t\mathbf{b})|\leq \sum_{m=0}^{+\infty}|b_{m}(z^{0}+t_0\mathbf{b})||t-t_{0}|^{m}=\sum_{m=0}^{\infty}a_{m}(z^{0})r_{N_{0}-n_{0}}^{m}\leq
\\ \leq a_{N_{0}}(z^{0})\sum_{m=0}^{+\infty}r_{N_{0}-n_{0}}^{m} =\frac{a_{N_{0}}(z^{0})}{1-r_{N_{0}-n_{0}}}\leq
\frac{a_{N_{0}}(z^{0})}{1-1/8}=\frac{8}{7}a_{N_{0}}(z^{0}).\label{eq_th_1322}
\end{gather}

From (\ref{eq_th_1321}) and (\ref{eq_th_1322}) we have
\begin{gather*}
\max\left\{|F(z^{0}+t\mathbf{b})|: |t-t_{0}|=\frac{r_{N_{0}-n_{0}}}{L(z^{0}+t_{0}\mathbf{b})}\right \}\leq\frac{8}{7}a_{N_{0}}(z^{0})\leq\\ \leq \frac{16}{7}\frac{1}{R_{N_{0}-n_{0}}}r^{-n_{0}}_{N_{0}-n_{0}}\min\left\{|F(z^{0}+t\mathbf{b})|: |t-t_{0}|=\frac{r_{N_{0}-n_{0}}}{L(z^{0}+t_{0}\mathbf{b})}\right\}\leq\\
\leq \frac{16}{7}\frac{1}{R_{N}}r^{-N}_{N}\min\left\{|F(z^{0}+t\mathbf{b})|: |t-t_{0}|=\frac{r_{N_{0}-n_{0}}}{L(z^{0}+t_{0}\mathbf{b})}\right\},
\end{gather*}
i.e.  (\ref{riv51}) holds with $\displaystyle P_{2}(R)=\frac{16}{7R_{N}r_{N}^{N}},$ $\eta(R)=r_{N}=\displaystyle\frac{1}{8R_{N}}$ and  $r=r_{N_{0}-n_{0}}.$

\textbf{Sufficiency.} In view of Theorem \ref{te4} it is
enough to prove there exists number $P_{1}$ that for every $z^{0}\in
\mathbb{B}_n,$ $t_{0}\in  S_{z^0}$
\begin{multline}\label{riv54}
\max\left\{|F(z^{0}+t\mathbf{b})|: |t-t_{0}|=\frac{\beta+1}
{2L(z^{0}+t_{0} \mathbf{b})}\right\}\leq\\
%\end{equation*}
%\begin{equation} \label{riv54}
 \leq
P_{1}\max\left\{|F(z^{0}+t\mathbf{b})|:
|t-t_{0}|=\frac{\beta-1}{4\beta L(z^{0}+t_{0}\mathbf{b})}\right\}.
\end{multline}
Let $\widetilde{R}=\frac{\beta-1}{4\beta}$. Then there exist
$P^{*}_{2}=P_{2}\big(\widetilde{R}\big)$ and
$\eta=\eta\big(\widetilde{R}\big)\in\big(0,\widetilde{R}\big)$
 that for every $z^{*}\in \mathbb{B}_n,$  $t^{*}\in
 S_{z^*}$ and some $r\in\big[\eta,\widetilde{R}\big]$
the following inequality is valid
\begin{gather*}
\max\big\{|F(z^{*}+t\mathbf{b})|:
|t-t^{*}|=\frac{r}{L(z^{0}+t^{*}\mathbf{b})}\big\} \leq \\ \leq 
P^{*}_{2}\min\big\{|F(z^{*}+t\mathbf{b})|:
|t-t^{*}|\!=\!\frac{r}{L(z^{0}+t^{*}\mathbf{b})}\big\}.
\end{gather*}
 Let $L^{*}\!=\!\max\!\{L(z^{0}+t\mathbf{b}):
|t-t_{0}|\leq{\beta}/{L(z^{0}\!+\!t_{0} \mathbf{b})}\!\},$
$\rho_{0}\!=\!{(\beta\!-\!1)}/{(4\beta L\!(z^{0}\!+\!t_{0}\mathbf{b}))},$\
$\rho_{k}=\rho_{0}+{k\eta}/{L^{*}},$\ $ k\in \mathbb{Z}_{+}$.
Hence, $\frac{\eta}{L^{*}}<\frac{\beta-1}
{4\beta L(z^{0}+t_{0}\mathbf{b})}<\frac{\beta}{L(z^{0}+t_{0}\mathbf{b})}
-\frac{\beta+1}{2L(z^{0}+t_{0}\mathbf{b})}.$\ Therefore,  there exists
$n^{*}\in \mathbb{N}$, which does not depend on $z^{0}$ and $t_{0}$  that
$\rho_{p-1}<\frac{\beta+1}{2L(z^{0}+t_{0}
\mathbf{b})}\leq\rho_{p}\leq\frac{\beta}{L(z^{0}+t_{0}\mathbf{b})}$
for some $p=p(z^{0},t_{0})\leq n^{*}.$

Let $c_{k}=\{t\in\mathbb{C}: |t-t_{0}|=\rho_{k}\},$\
$|F(z^{0}+t_{k}^{**}\mathbf{b})|=\max\{|F(z^{0}+t\mathbf{b})|:
t\in c_{k}\}$ and $t_{k}^{*}$ be a point of intersection of the segment
$[t_{0},t^{**}_{k}]$ with the circle $c_{k-1}$. Then for every $r>\eta$
 the following inequality holds $|t^{**}_{k}-t^{*}_{k}|={\eta}/{L^{*}}
\leq{r}/{L(z^{0}+t^{*}_{k}\mathbf{b})}$. Hence, for some
$r\in[\eta,\widetilde{R}]$ the following inequality is valid
\begin{gather*}
|F(z^{0}+t_{k}^{**}\mathbf{b})|\leq\max\left\{|F(z^{0}+t\mathbf{b})|:
|t-t_{k}^{*}|=\frac{r}{L(z^{0}+t^{*}_{k}\mathbf{b})}\right\}\leq
\\ \leq P_{2}^{*}\min\left\{|F(z^{0}+t\mathbf{b})|:
|t-t_{k}^{*}|=\frac{r}{L(z^{0}+t^{*}_{k}\mathbf{b})}\right\}\leq 
P_{2}^{*}\max\{|F(z^{0}+t\mathbf{b})|: t\in c_{k-1}\}.
\end{gather*}
Therefore, we get inequality (\ref{riv54}) with
 $P_{1}^{*}=(P_{2}^{*})^{n^{*}}$
\begin{multline*} \max\left\{|F(z^{0}+t\mathbf{b})|:
|t-t_{0}|=\frac{\beta+1}{2L(z^{0}+t_{0}\mathbf{b})}\right\}\leq
\max\{|F(z^{0}+t\mathbf{b})|: t\in c_{p}\}\leq \\\leq
P_{2}^{*}\max\{|F(z^{0}+t\mathbf{b})|: t\in
c_{p-1}\}\leq\ldots\leq(P_{2}^{*})^{p}\max\{|F(z^{0}+t\mathbf{b})|:
t\in c_{0}\}\leq\\
\leq(P_{2}^{*})^{n^{*}}\max\left\{|F(z^{0}+t\mathbf{b})|:
|t-t_{0}|=\frac{\beta-1}{4\beta L(z^{0}+t_{0}\mathbf{b})}\right\}.\end{multline*}
Theorem \ref{te5} is proved.
\end{proof}

% \begin{center}
\section{Logarithmic derivative and zeros.} 
% \end{center}  

% {\bf {$6^{0}.$ }}

Below we prove another criterion of $L$-index boundedness in a direction that describes behaviour of the directional
logarithmic derivative and distribution of zeros.

We need some additional denotations.

Denote $g_{z^0}(t):=F(z^{0}+t\mathbf{b}).$
If for a given $z^0\in\mathbb{B}_n$ \ \
$g_{z^0}(t)\not= 0$\  for all $t\in S_{z^0}$,
 then  $G^{\mathbf{b}}_{r}(F,z^{0}):=\emptyset;$
if for a given $z^0\in\mathbb{B}_n$ \ \ $g_{z^0}(t)\equiv 0,$ then
$G^{\mathbf{b}}_r(F,z^0):=\{z^0+t\mathbf{b}\colon
t\in S_{z^0}\}.$
 And if for a given $z^{0}\in
\mathbb{B}_n$\ \ $g_{z^0}(t)\not\equiv 0$ and $a_{k}^{0}$
are  zeros of $g_{z^0}(t)$,
% i.e. $F(z^{0}+a_{k}^{0}\mathbf{b})=0,$
then %we denote
\[G^{\mathbf{b}}_{r}(F,z^{0}):=\bigcup_{k}\left\{z^{0}+t\mathbf{b}\colon
|t-a_{k}^{0}|\leq\frac{r}{L(z^{0}+a^{0}_{k}\mathbf{b})} \right\},\
\ r>0.\]
% }

Let
\begin{equation} \label{riv60}
G^{\mathbf{b}}_{r}(F)=\bigcup_{z^{0}\in\mathbb{B}_n}
{G^{\mathbf{b}}_{r}(F,z^{0})}.
\end{equation}

We remark that if $L(z)\equiv 1$, then
$G^{\mathbf{b}}_{r}(F)\subset\left\{z\in\mathbb{B}_n:
dist(z,\mathbb{Z}_{F})<r|\mathbf{b}|\right\},$ where $\mathbb{Z}_{F}$ is a zero set of
function $F$. By
$n\big(r,z^{0},t_{0},{1}/{F}\big)=\sum_{|a_{k}^{0}-t_{0}|\leq
r}\,1$ we denote a counting function of zeros $a_{k}^{0}$.
\begin{theorem} \label{te6}
Let $F(z)$ be an analytic in $\mathbb{B}_n$ function,
$L\in{Q}_{\mathbf{b},\beta}(\mathbb{B}_n)$ and $\mathbb{B}_n\setminus
G^{\mathbf{b}}_{\beta}(F)\neq \emptyset$. $F(z)$ is of bounded
 $L$-index in  direction $\mathbf{b}\in\mathbb{C}^{n}$
if and only if
\begin{enumerate} \renewcommand{\labelenumi}{\arabic{enumi})}
\item \label{teum61} for every $r\in(0,\beta]$ there exists  $P=P(r)>0$  that
for each $z\in \mathbb{B}_n\backslash G^{\mathbf{b}}_{r}(F)$
\begin{equation} \label{riv61}
\left|\frac{1}{F(z)}\frac{\partial
F(z)}{\partial\mathbf{b}}\right|\leq PL(z);
\end{equation}
\item \label{teum62} for every $r\in(0,\beta]$ there exists $\widetilde{n}(r)\in
\mathbb{Z}_{+}$ that for each $z^{0}\in \mathbb{B}_n$ with
 $F(z^{0}+t\mathbf{b})\not \equiv 0,$ and for each
 $t_{0}\in  S_{z^0}$
\begin{equation} \label{riv62}
n\left(\frac{r}{L(z^{0}+t^{0}\mathbf{b})},z^{0},t_{0},\frac{1}{F}\right)
\leq\widetilde{n}(r).
\end{equation}
\end{enumerate}
\end{theorem}
\begin{proof}

\textbf{Necessity.} First, we prove that if  function $F(z)$ is of
bounded $L$-index in a direction, then for every
$\widetilde{z}^{0}=z^{0}+t_{0}\mathbf{b}\in
\mathbb{B}_n\backslash G^{\mathbf{b}}_{r}(F)$ $(r\in(0,\beta])$ and for every
$\widetilde{a}^{k}=z^{0}+a^{0}_{k}\mathbf{b}$ the following inequality holds
\begin{equation} \label{riv63} |\widetilde{z}^{0}-\widetilde{a}_{k}|>
\frac{r|\mathbf{b}|}{2L(\widetilde{z}^{0})\lambda^{\mathbf{b}}_{2}\left(z^{0},
r\right)}.\end{equation}
On the contrary,  we assume that there exists
 $\widetilde{z}^{0}=z^{0}+t_{0}\mathbf{b}\in
\mathbb{B}_n\backslash G^{\mathbf{b}}_{r}(F)$ and
$\widetilde{a}^{k}=z^{0}+a^{0}_{k}\mathbf{b}$ that
$$|\widetilde{z}^{0}-\widetilde{a}_{k}|\leq\displaystyle\frac{\mathstrut r |\mathbf{b}|}{2L(\widetilde{z}^{0})
\lambda^{\mathbf{b}}_{2}\left(z^{0},r\right)} \leq
\frac{r|\mathbf{b}|}{2L(\widetilde{z}^{0})}<\frac{r|\mathbf{b}|}{L(\widetilde{z}^{0})}.$$
Hence, $|t_{0}-a^{0}_{k}|<\frac{r}{L(\widetilde{z}^{0})}.$
But for $\lambda^{\mathbf{b}}_{2}$ the following estimate holds 
$$L(\widetilde{a}^{k})\leq\lambda^{\mathbf{b}}_{2}\left(z^{0},
r\right)L(\widetilde{z}^{0}),$$ and therefore
$$|\widetilde{z}^{0}-\widetilde{a}^{k}|=
|\mathbf{b}|\cdot |t_{0}-a^{0}_{k}|\leq\frac{r|\mathbf{b}|}{2L(\widetilde{a}^{k})},$$ i.e.
$|t_{0}-a^{0}_k|
\leq\displaystyle\frac{\mathstrut r}{2L(\widetilde{a}^{k})}.$
 We obtained a contradiction with $\widetilde{z}^{0}\in \mathbb{C}^{n}\backslash
G^{\mathbf{b}}_{r}(F)$. In fact, in (\ref{riv63}) instead of
$\lambda^{\mathbf{b}}_{2}\left(z^{0},r\right)$
we can take $\lambda^{\mathbf{b}}_{2}\left(r\right)$.

We choose in Theorem \ref{te5} $R=\displaystyle\frac{\mathstrut
r}{2 \lambda_{2}^{\mathbf{b}}\left(r\right)}.$
Then there exists $P_{2}\geq 1$ and $\eta\in(0,R)$ that for every
$\widetilde{z^{0}}=z^{0}+ t_{0}\mathbf{b}\in\mathbb{B}_n$  and
some $r^{*}\in[\eta,R]$ inequality (\ref{riv51}) holds with $r^{*}$
instead of $r$. Therefore,  by Cauchy inequality
\begin{gather}
\left|\frac{\partial
F(z^{0}+t_{0}\mathbf{b})}{\partial\mathbf{b}}\right|
 \leq\frac{L(z^{0}+t_{0}\mathbf{b})}{r^{*}}\max\Big\{|F(z^{0}+t\mathbf{b}):
 |t-t_{0}|=\frac{r^{*}}{L(z^{0}+t_{0}\mathbf{b})}\Big\}
\leq \nonumber\\ \leq
P_{2}\frac{L(z^{0}+t_{0}\mathbf{b})}{\eta}\,\min\{|F(z^{0}+t\mathbf{b})|:
|t-t_{0}|=\frac{r^{*}}{L(z^{0}+t_{0}\mathbf{b})}\} \label{riv64}
\end{gather}

Since (\ref{riv63}) for every $z^{0}+t_{0}\mathbf{b}\in \mathbb{B}_n\backslash
G^{\mathbf{b}}_{r}(F),$\ a set
$$\left\{z^{0}+t\mathbf{b}: |t-t_{0}|\leq\frac{r}
{2 \lambda^{\mathbf{b}}_{2}\left(r\right)
L(z^{0}+t_{0}\mathbf{b})}\right\}$$ does not contain zeros of function $F(z^{0}+t\mathbf{b}).$ Therefore, 
applying to $1/F$, as a function of variable $t$, a maximum principle, we have
\begin{equation} \label{riv65} |F(z^{0}+t_{0}\mathbf{b})|
\geq\min\left\{|F(z^{0}+t\mathbf{b})|: |t-t_{0}|=\frac{r^{*}}
{L(z^{0}+t_{0}\mathbf{b})}\right\}
\end{equation}
 The inequalities (\ref{riv64}) and (\ref{riv65}) imply (\ref{riv61}) with
$P=\displaystyle\frac{P_{2}}{\eta}.$

Now we prove that if $F$ is of bounded $L$-index in direction
$\mathbf{b}$ then there exists $P_{3}>0$ that for every  $z^{0}\in
\mathbb{B}_n,$ $t_{0}\in  S_{z^0},$ $r\in(0,1]$
\begin{gather}
n\Big(\frac{r}{L(z^{0}+t_{0}\mathbf{b})},z^{0},t_{0},1/F\Big)
\min\Big\{|F(z^{0}+t\mathbf{b})|:
|t-t_{0}|=\frac{r}{L(z^{0}+t_{0}\mathbf{b})}\Big\} \leq\nonumber\\
\leq P_{3}\max\Big\{|F(z^{0}+t\mathbf{b})|: |t-t_{0}|=\frac{1}
{L(z^{0}+t_{0}\mathbf{b})}\Big\}. \label{riv66}\end{gather}
By Cauchy inequality and Theorem \ref{te3} for all $t\in S_{z^0}$ with circle 
 $|t-t_{0}|=\displaystyle\frac{1}{L(z^{0}+t_{0}\mathbf{b})}$ we have
\begin{gather}
\Big|\frac{\partial F(z^{0}+t\mathbf{b})}{\partial\mathbf{b}}\Big|
\leq\frac{L(z^{0}+t_{0}\mathbf{b})}{\beta-1}\max\Big\{|F(z^{0}+\theta\mathbf{b})|:
|\theta-t|=\frac{\beta-1}{L(z^{0}+t_{0}\mathbf{b})}\Big\}\leq\nonumber\\
\leq\frac{L(z^{0}+t_{0}\mathbf{b})}{\beta-1}\max\Big\{|F(z^{0}+t\mathbf{b})|:
|t-t_{0}|=\frac{\beta}{L(z^{0}+t_{0}\mathbf{b})}\Big\}\leq\nonumber\\
\!\leq\!\frac{P_{1}(1,\beta)}{\beta-1}L(z^{0}+t_{0}\mathbf{b})\max\left\{|F(z^{0}\!+\!t\mathbf{b})|:
|t-t_{0}|=\frac{1}{L(z^{0}+t_{0}\mathbf{b})}\!\right\}. \label{dopeql1}
\end{gather}
If $F(z^{0}+t\mathbf{b})\neq 0$ on a circle $\left\{t\in S_{z^0}:
|t-t_{0}|=\frac{r}{L(z^{0}+t_{0}\mathbf{b})}\right\},$ then
\begin{gather}
n\left(\frac{r}{L(z^{0}+t_{0}\mathbf{b})},z^{0},t_{0},\frac{1}{F}\right) =
\left|\frac{1}{2\pi\,i}\int\limits_{|t-t_{0}|=\frac{r}{L(z^{0}+t_{0}\mathbf{b})}}
\frac{\partial F(z^{0}+t\mathbf{b})}{\partial\mathbf{b}}\frac{1}
{F(z^{0}+t\mathbf{b})}dt\right|\leq\nonumber\\
\leq\frac{\max\Big\{\Big|\frac{\partial F(z^{0}+t\mathbf{b})}
{\partial\mathbf{b}}\Big|:|t-t_{0}|=\frac{r}{L(z^{0}+t_{0}\mathbf{b})}\Big\}}
{\min\left\{|F(z^{0}+t\mathbf{b})|:
|t-t_{0}|=\frac{r}{L(z^{0}+t_{0}\mathbf{b})}\right\}}\,
\frac{r}{L(z^{0}+t_{0}\mathbf{b})}. \label{dopeql2}
\end{gather}

From \eqref{dopeql1} and \eqref{dopeql2} we have  \begin{gather*}
n\left(\frac{r}{L(z^{0}+t_{0}\mathbf{b})},z^{0},t_{0},1/F\right) %\times\\ \times
\min\left\{|F(z^{0}+t\mathbf{b})|:|t-t_{0}|=\frac{r}{L(z^{0}+t_{0}\mathbf{b})}\right\}
\leq\\\leq\frac{r}{L(z^{0}+t_{0}\mathbf{b})}\max\left\{\left|\frac{\partial
F(z^{0}+t\mathbf{b})} {\partial\mathbf{b}}\right|:
|t-t_{0}|=\frac{r}{L(z^{0}+t_{0}\mathbf{b})}\right\}\leq\\
 \leq\frac{1}{L(z^{0}+t_{0}\mathbf{b})}
\max\left\{\left|\frac{\partial F(z^{0}+t\mathbf{b})}
{\partial\mathbf{b}}\right|: |t-t_{0}|=\frac{1}{L(z^{0}+t_{0}\mathbf{b})}\right\}\leq\\
\leq\frac{P_{1}(1,\beta)}{\beta-1}\max\left\{|F(z^{0}+t\mathbf{b})|:
|t-t_{0}| =\frac{1}{L(z^{0}+t_{0}\mathbf{b})}\right\}.
\end{gather*}
Thus, we obtain \eqref{riv66} with $\displaystyle P_{3}=\frac{P_{1}(1,\beta)}{\beta-1}$. If function
$F(z^{0}+t\mathbf{b})$ has zeros on the circle
 $\left\{t\in D_R^{z^0}:\right.$ $\left.|t-t_{0}|=\frac{r}{L(z^{0}+t_{0}\mathbf{b})}\right\}$  then an inequality \eqref{riv66} is obvious.

Now we put $R=1$ in Theorem \ref{te5}. Then there exists
$P_{2}=P_{2}(1)\geq 1$ and $\eta\in (0,1)$ that for each
$z^{0}\in \mathbb{B}_n,$ $t_{0}\in  S_{z^0}$ and some
$r^{*}=r^{*}(z^{0},t_{0})\in[\eta,1]$
\begin{gather*} \max\left\{|F(z^{0}+t\mathbf{b})|:
|t-t_{0}|=\frac{r^{*}}{L(z^{0}+t_{0}\mathbf{b})}\right\}\leq\\ \leq
P_{2}\min\left\{|F(z ^{0}+t\mathbf{b})|:
|t-t_{0}|=\frac{r^{*}}{L(z^{0}+t_{0}\mathbf{b})}\right\}.
\end{gather*}
Besides, by
Theorem \ref{te3} there exists $P_{1}\geq 1$ such that for all $z^{0}\in
\mathbb{B}_n,$ $t_{0}\in  S_{z^0}$
\begin{gather*} \max\left\{|F(z^{0}+t\mathbf{b})|: |t-t_{0}|=\frac{1}
{L(z^{0}+t_{0}\mathbf{b})}\right\}\leq\\ \leq
P_{1}(1,\eta)\max\left\{|F(z^{0}+t\mathbf{b})|:\
 |t-t_{0}|=\frac{\eta}{L(z^{0}+t_{0}\mathbf{b})}\right\}
\leq\\ \leq P_{1}(1,\eta)\max\left\{|F(z^{0}+t\mathbf{b})|: |t-t_{0}|=
\frac{r^{*}}{L(z^{0}+t_{0}\mathbf{b})}\right\} \leq \\
\leq P_{1}(1,\eta) P_{2}\min\left\{|F(z^{0}+t\mathbf{b})|:
|t-t_{0}|=\frac{r^{*}}{L(z^{0}+t_{0}\mathbf{b})}\right\}.
\end{gather*}

Since \eqref{riv66}, we have
\begin{gather*}
n\left(\frac{r^{*}}{L(z^{0}+t_{0}\mathbf{b})},z^{0},t_{0}, \frac{1}{F}\right) %\times\\ \times
\min\left\{|F(z^{0}+t\mathbf{b})|:|t-t_{0}|=\frac{r^{*}}
{L(z^{0}+t_{0}\mathbf{b})}\right\}\leq\\ \leq P_{3}P_{1}(1,\eta)
P_{2}\min\left\{|F(z^{0}+t\mathbf{b})|: |t-t_{0}|=\frac{r^{*}}
{L(z^{0}+t_{0}\mathbf{b})}\right\},\end{gather*}
i.e.
$n\left(\frac{r^{*}}{L(z^{0}+t_{0}\mathbf{b})},z^{0},t_{0},\frac{1}{F}\right)\leq
P_{1}(1,\eta)P_{2}P_{3}.$ Hence,
\begin{gather*}
n\left(\frac{r^{*}}{L(z^{0}+t_{0}\mathbf{b})},z^{0},t_{0},\frac{1}{F}\right)\leq
P_{4}=
P_{1}(1,\eta)P_{2}P_{3}%=\\
=\frac{P_{1}(1,\eta)P_{2}(1)P_{1}(1,r+1)}{r}.\end{gather*}
If $r\in(0,\eta]$ then property (\ref{riv62}) is proved.

Let  $r\in(\eta,\beta]$ and $L^{*}=\max\left\{L(z^{0}+t\mathbf{b}):
|t-t_{0}|=\frac{r}{L(z^{0}+t_{0}\mathbf{b})}\right\}.$
Using properties of $Q^n_{\mathbf{b}},$ we have $L^{*}\leq \lambda_{2}^{\mathbf{b}}(r)L(z^{0}+t_{0}\mathbf{b}).$
Put
$\rho=\frac{\eta}{L(z^{0}+t_{0}\mathbf{b})\lambda^{\mathbf{b}}_{2}(r)},$
 $R=\frac{r}{L(z^{0}+t_{0}\mathbf{b})}.$ We can cover every set
 $\overline{K}=\{z^{0}+t\mathbf{b}: |t-t_{0}|\leq R\}$
 by a finite number $m=m(r)$ of closed sets
$\overline{K}_{j}=\{z^{0}+t\mathbf{b}: |t-t_{j}|\leq \rho\},$ where
$t_{j}\in\overline{K}$. Since
$$\frac{\eta}{\lambda^{\mathbf{b}}_{2}(r)L(z^{0}+t_{0}\mathbf{b})}
\leq\frac{\eta}{L^{*}}\leq\frac{\eta}{L(z^{0}+t_{j}\mathbf{b})}$$
  in each $\overline{K}_{j}$ there are at most $[P_{4}]$ zeros of
function $F(z^{0}+t\mathbf{b})$. Thus,
$$n\left(\frac{r}{L(z^{0}+t_{0}\mathbf{b})},z^{0},t_{0},
{1}/{F}\right)\leq\widetilde{n}(r)=[P_{4}]\,m(r)$$ and property
(\ref{riv62}) is proved.

\noindent\textbf{Sufficiency.} On the contrary, suppose that conditions
(\ref{riv61}) and (\ref{riv62}) hold. By condition (\ref{riv62})
for every $R\in(0,\beta]$  there exists $\widetilde{n}(R)\in \mathbb{Z}_{+}$
 that in each set
$$\overline{K}=\left\{z^{0}+t\mathbf{b}: |t-t_{0}|\leq
\frac{R}{L(z^{0}+t_{0}\mathbf{b})}\right\}$$ the number of zeros of
$F(z^{0}+t\mathbf{b})$ does not exceed $\widetilde{n}(r)$.

We put $a=a(R)=\frac{R\lambda^{\mathbf{b}}_{1}(R)}
{2(\widetilde{n}(R)+1)}.$ By condition (\ref{riv61}) there exists
$P=P(a)=\widetilde{P}(R)\geq 1$ that $\left|\displaystyle\frac{\mathstrut\partial
F(z)}{\partial\mathbf{b}}\frac{1}{F(z)}\right|\leq PL(z)$ for all
$z\in \mathbb{B}_n\backslash G^{\mathbf{b}}_{a}$, that is for all $z\in
\overline{K}$ lying  outside the sets
$$b^{0}_{k}=\left\{z^{0}+t\mathbf{b}:\
|t-a^{0}_{k}|<\displaystyle\frac{\mathstrut a(R)}{ L(z^{0}+a_{k}^{0}\mathbf{b})}\right\},$$
where $a^{0}_{k}\in\overline{K}$ are zeros of function $F(z^{0}+t\mathbf{b})\not \equiv 0.$  
Since properties $\lambda^{\mathbf{b}}_{1}$ we have 
$$\lambda^{\mathbf{b}}_{1}(R)L(z^{0}+t_{0}\mathbf{b})
\leq\lambda^{\mathbf{b}}_{1}(R,z^{0})L(z^{0}+t_{0}\mathbf{b})\leq
L(z^{0}+ a^{0}_{k}\mathbf{b}).$$ Therefore, 
$\left|\displaystyle\frac{\mathstrut 1}{F(z)}\frac{\partial F(z)}
{\partial\mathbf{b}}\right|\leq PL(z)$ for all $z\in\mathbb{B}_n$, lying
outside the sets
\begin{gather*}
c^{0}_{k}=\left\{z^{0}+t\mathbf{b}: |t-a^{0}_{k}|\leq\frac{a(R)}
{ \lambda^{\mathbf{b}}_{1}(R)L(z^{0}+t_{0}\mathbf{b})}= %\right.\\ \left. =
\frac{R}{2 (\widetilde{n}(R)+1)L(z^{0}+t_{0}\mathbf{b})}\right\}.
\end{gather*}
Obviously, the sum of diameters of sets $c^{0}_{k}$ does not
exceed
$$\frac{R\widetilde{n}(R)}{
(\widetilde{n}(R)+1)L(z^{0}+t_{0}\mathbf{b})}<\frac{R}
{ L(z^{0}+t_{0}\mathbf{b})}.$$ Therefore,  there exist a set
 $\widetilde{c}^{0}=\left\{z^{0}+t\mathbf{b}: |t-t_{0}|=\frac{r}
{L(z^{0}+t_{0}\mathbf{b})}\right\},$ where
$$\displaystyle\frac{\mathstrut R}{2 (\widetilde{n}(R)+1)}=\eta(R)<r<
R,$$ such that for all $z\in
\widetilde{c}^{0}$ the following inequality is valid
\begin{gather*} \left|\frac{1}{F(z)}\frac{\partial F(z)}{\partial\mathbf{b}}\right|\leq PL(z)
\leq P\lambda^{\mathbf{b}}_{2}(r)L(z^{0}+t_{0}\mathbf{b})\leq %\\ \leq
P\lambda^{\mathbf{b}}_{2}\left(R\right)
L(z^{0}+t_{0}\mathbf{b}).\end{gather*}
 For any points
$z_{1}=z^{0}+t_{1}\mathbf{b}$ and $z_{2}=z^{0}+t_{2}\mathbf{b}$ with
$\widetilde{c}^{0}$ we have
\begin{gather*}
\ln\left|\frac{F(z^{0}+t_{1}\mathbf{b})}{F(z^{0}+t_{2}\mathbf{b})}
\right|\leq\int_{t_{1}}^{t_{2}}\Big|\frac{1}{F(z^{0}+t\mathbf{b})}\frac{\partial
F(z^{0}+t\mathbf{b})} {\partial\mathbf{b}}\Big||dt| \leq\\ \leq
P\lambda^{\mathbf{b}}_{2}\left(R\right)L(z^{0}+t_{0}\mathbf{b}) \frac{2r}
{L(z^{0}+t_{0}\mathbf{b})}\leq
2R\,P(R)
\lambda^{\mathbf{b}}_{2}\left(R\right).\end{gather*}
Hence, we get  
\begin{gather*}\max\left\{|F(z^{0}+t\mathbf{b})|: |t-t_{0}|=\frac{r}
{L(z^{0}+t_{0}\mathbf{b})}\right\}\leq \\ \leq
P_{2}\!\min\!\left\{|F(z^{0}+
t\mathbf{b})|:\ |t-t_{0}|=\frac{r}{L(z^{0}+t_{0}\mathbf{b})}\right\}, \end{gather*}
where $P_{2}=\exp\left\{\mathstrut 2R\,P(R)
\lambda^{\mathbf{b}}_{2}\left(R\right)\right\}.$
Thus, by Theorem \ref{te5} the function
$F(z)$ is of bounded $L$-index in  direction
$\mathbf{b}$. Theorem \ref{te6} is proved.
\end{proof}

% \begin{center}
\section{Boundedness $L$-index in the direction of analytical solutions of some partial
differential equations.} 
% \end{center}  

% {\bf{$7^{0}.$ }}

We consider a partial differential equation
\begin{equation} \label{riv140}
g_{0}(z)\frac{\partial^{p}w}{\partial\mathbf{b}^{p}}+g_{1}(z)\frac{\partial^
{p-1}w}{\partial\mathbf{b}^{p-1}}+\ldots+g_{p}(z)w=h(z).
\end{equation}
First, we prove an auxiliary assertion.
\begin{lemma} \label{te15}
Let $\beta>1,$ $L\in Q_{\mathbf{b},\beta}(\mathbb{B}_n),$ $F(z)$ be an analytic in $\mathbb{B}_n$ function of bounded $L$-index in direction
 $\mathbf{b}\in \mathbb{C}^{n},$  $\mathbb{B}_n\backslash
G^{\mathbf{b}}_{\beta}(F)\neq \emptyset$. Then for every $r\in(0,\beta]$ and for every
$m\in\mathbb{N}$ there exists $P=P(r,m)>0$
 such that for all $z\in\mathbb{B}_n\backslash G^{\mathbf{b}}_{r}(F)$ inequality holds
$$\left|\frac{\partial^{m}F(z)}{\partial\mathbf{b}^{m}}\right|\leq PL^{m}(z)|F(z)|.$$
\end{lemma}
\begin{proof}
In Theorem~\ref{te6} we proved that if an entire function $F(z)$ is of bounded $L$-index in direction  $\mathbf{b},$ then (\ref{riv63}) holds, i.e. for each
$\widetilde{z}^{0}=z^{0}+t_{0}\mathbf{b}\in
\mathbb{B}_n\backslash G^{\mathbf{b}}_{r}(F)$ $(r\in(0,\beta])$ and 
$\widetilde{a}^{k}=z^{0}+a^{0}_{k}\mathbf{b}$ an inequality
 holds
\begin{equation} \label{riv635} |\widetilde{z}^{0}-\widetilde{a}_{k}|>
\frac{r|\mathbf{b}|}{2L(\widetilde{z}^{0})\lambda^{\mathbf{b}}_{2}(z^{0},
r/( ))}.
\end{equation}
We put in Theorem \ref{te5} $R=\displaystyle\frac{
r}{2 \lambda_{2}^{\mathbf{b}}(r)}.$ Then there exist
$P_{2}=P_{2}\left(\displaystyle\frac{
r}{2 \lambda_{2}^{\mathbf{b}}(r)}\right)\geq 1$
and $\eta\left(\displaystyle\frac{\mathstrut
r}{2 \lambda_{2}^{\mathbf{b}}(r)}\right)\in\left(0,\displaystyle\frac{\mathstrut
r}{2 \lambda_{2}^{\mathbf{b}}(r)}\right)$
 that for all $z^{0}\in \mathbb{B}_n,$
 $t_{0}\in  S_{z^0}$ and some
$r^{*}=r^{*}(z^{0},t_{0})\in\left[\eta\left(\displaystyle\frac{\mathstrut
r}{2 \lambda_{2}^{\mathbf{b}}(r)}\right),\displaystyle\frac{\mathstrut
r}{2 \lambda_{2}^{\mathbf{b}}(r)}\right]$
an inequality (\ref{riv51}) holds with $r^{*}$ instead of $r.$
Using Cauchy inequality, we get 
\begin{gather*}
 \frac{1}{m!}\left|\frac{\partial^{m} F(z^{0}+t_{0}\mathbf{b})}{\partial \mathbf{b}^{m}}\right|\leq%\\ \leq 
 \left(\frac{L(z^{0}+t_{0}\mathbf{b})}{r^{*}}\right)^{m}\max\left\{|F(z^{0}+t\mathbf{b})|:
|t-t_{0}|=\frac{r^{*}}{L(z^{0}+t_{0}\mathbf{b})}\right\}\leq\\
\leq P_{2}\left(\frac{L(z^{0}+t_{0}\mathbf{b})}{\eta}\right)^{m}\min\left\{|F(z^{0}+t\mathbf{b})|: |t-t_{0}|=\frac{r^{*}}{L(z^{0}+t_{0}\mathbf{b})}\right\}.
\end{gather*}
From (\ref{riv635}) for every $z^{0}\in \mathbb{B}_n\backslash
G^{\mathbf{b}}_{r}(F)$ the set
$$\left\{z^{0}+t\mathbf{b}: |t-t_{0}|\leq\displaystyle\frac{\mathstrut r}{2 \lambda_{2}^{\mathbf{b}}(r)L(z^{0}+t_{0}\mathbf{b})}\right\}$$  does not contain zeros of function $F(z^{0}+t\mathbf{b}).$ Therefore, applying to $\frac{1}{F(z^0+t\mathbf{b})}$ a maximum modulus principle in variable $t\in S_{z^0},$ we have
$$|F(z^{0}+t_{0}\mathbf{b})|\geq\min\left\{|F(z^{0}+t\mathbf{b})|:
 |t-t_{0}|=\frac{r^{*}}{L(z^{0}+t_{0}\mathbf{b})}\right\}.$$

Thus,
$$\left|\frac{\partial^{m} F(z^{0}+t_{0}\mathbf{b})}{\partial\mathbf{b}^{m}}\right|\leq
m!\frac{P_{2}}{\eta^{m}}L^{m}(z^{0}+t_{0}\mathbf{b})|F(z^{0}+t_{0}\mathbf{b})|.$$
Hence, we proved a needed inequality with $P=P_{2}m!\eta^{-m}.$
\end{proof}
Using Lemma \ref{te15}, we deduce a following theorem.
\begin{theorem} \label{te14}
Let $\beta>1,$ $L\in{Q}_{\mathbf{b},\beta}(\mathbb{B}_n),$ $g_{0}(z),\ldots,g_{p}(z),h(z)$
be analytic in $\mathbb{B}_n$ functions of bounded $L$-index in  direction $\mathbf{b},$ $\mathbb{B}_n\backslash
G^{\mathbf{b}}_{\beta}(g_{0})\neq \emptyset$ and
for every $r\in(0;\beta]$ there exists $T=T(r)>0$ that for each
$z\in\mathbb{B}_n \backslash G_{r}^{\mathbf{b}}(g_{0})$ and
$j=1,\ldots,p$ inequality holds
\begin{equation} \label{riv1410} |g_{j}(z)|\leq TL^{j}(z)|g_{0}(z)| .\end{equation}
Then an analytic function $F(z),$ $z\in\mathbb{B}_n,$ which satisfies
an equation (\ref{riv140}), is of bounded $L$-index in  direction $\mathbf{b}$.
\end{theorem}

\begin{proof}
For every given $z^{0}\in\mathbb{B}_n$ let $b_{k}^{0}$ be zeros of function $g_{0}(z^{0}+t\mathbf{b})$ and  $\{c^{0}_{k}\}$ be a set of zeros of all functions $g_{0}(z^{0}+t\mathbf{b}),$ $g_{1}(z^{0}+t\mathbf{b}),$ $\ldots, g_{p}(z^{0}+t\mathbf{b})$ and  $h(z^{0}+t\mathbf{b}),$ as functions of one variable $t\in S_{z^0}.$ Obviously, this inclusion is valid $\{b_{k}^{0}\}\subset\{c_{k}^{0}\}.$ We put
$$G_{r}^{\mathbf{b}}(z^{0})=\bigcup_{k}\left\{z^{0}+t\mathbf{b}:
 |t-c^{0}_{k}|\leq\frac{r}{L(z^{0}+c^{0}_{k}\mathbf{b})}\right\}, \ G_{r}^{\mathbf{b}}=\bigcup_{z^{0}}G_{r}^{\mathbf{b}}(z^{0}).$$

It is easy  to see that $G^{\mathbf{b}}_{r}=G^{\mathbf{b}}_{r}(h)\cup\bigcup_{j=1}^{p}G^{\mathbf{b}}_{r}(g_{j}).$
Suppose that $\mathbb{B}_n\setminus G^{\mathbf{b}}_{r}(g_{0})\neq \emptyset.$ Lemma \ref{te15} and equation (\ref{riv1410}) implies that for every $r\in(0,\beta]$ there exists $T^{*}=T^{*}(r)>0$ such that for all $z\in\mathbb{B}_n\setminus G_{r}^{\mathbf{b}}$ the following inequalities hold
$$\left|\frac{\partial h(z)}{\partial\mathbf{b}}\right|\leq T^{*}|h(z)|L(z),
 |g_{j}(z)|\leq T^{*}|g_{0}(z)|L^{j}(z),\  j\in\{1,2,\ldots,p,\}$$
$$\left|\frac{\partial g_{j}(z)}{\partial\mathbf{b}}\right|\leq P(r) L(z)|g_{j}(z)|\leq T^{*}(r)|g_{0}(z)|L^{j+1}(z), \  j\in\{0,1,2,\ldots,p,\}.$$
 In equation (\ref{riv140}) we evaluate a derivative in  direction $\mathbf{b}:$
$$g_{0}(z)\frac{\partial^{p+1} F(z)}{\partial\mathbf{b}^{p+1}}+\sum_{j=1}^{p}g_{j}(z)\frac{\partial^{p+1-j} F(z)}{\partial\mathbf{b}^{p+1-j}}+\sum_{j=0}^{n}\frac{\partial g_{j}(z)}{\partial\mathbf{b}}\frac{\partial^{p-j}F(z)}{\partial\mathbf{b}^{p-j}}=\frac{\partial h(z)}{\partial\mathbf{b}}.$$
This obtained equality implies that for all $z\in\mathbb{B}_n\setminus G_{r}^{\mathbf{b}}:$
\begin{gather*}
 |g_{0}(z)|\left|\frac{\partial^{p+1} F(z)}{\partial\mathbf{b}^{p+1}}\right|\leq
 \left|\frac{\partial h(z)}{\partial\mathbf{b}}\right|+\sum_{j=1}^{p}|g_{j}(z)|\left|\frac{\partial^{p+1-j}
 F(z)}{\partial\mathbf{b}^{p+1-j}}\right|+\\ +\sum_{j=0}^{p}\left|\frac{\partial
 g_{j}(z)}{\partial\mathbf{b}}\right|\left|\frac{\partial^{p-j} F(z)}{\partial\mathbf{b}^{p-j}}\right|
 \!\leq\! T^{*}|h(z)|L(z)+\sum_{j=1}^{p}|g_{j}(z)|\left|\frac{\partial^{p+1-j}
 F(z)}{\partial\mathbf{b}^{p+1-j}}\right|\!+\!\\ +\sum_{j=0}^{p}\left|\frac{\partial
 g_{j}(z)}{\partial\mathbf{b}}\right|\left|\frac{\partial^{p-j} F(z)}{\partial\mathbf{b}^{p-j}}\right|
\leq T^{*}L(z)\sum_{j=0}^{p}|g_{j}(z)|\left|\frac{\partial^{p-j} F(z)}{\partial\mathbf{b}^{p-j}}\right|+
 \\ +\sum_{j=1}^{p}|g_{j}(z)|\left|\frac{\partial^{p+1-j}
  F(z)}{\partial\mathbf{b}^{p+1-j}}\right|+\sum_{j=0}^{p}\left|\frac{\partial
 g_{j}(z)}{\partial\mathbf{b}}\right|\left|\frac{\partial^{p-j} F(z)}{\partial\mathbf{b}^{p-j}}\right|\leq\\
 \!\leq T^{*}|g_{0}(z)|\left(T^{*}L(z)\sum_{j=0}^{p}L^{j}(z)\left|\frac{\partial^{p-j} F(z)}{\partial\mathbf{b}^{p-j}}\right|+\sum_{j=1}^{p}L^{j}(z)\left|\frac{\partial^{p+1-j} F(z)}{\partial\mathbf{b}^{p+1-j}}\right|\!+\!\right.\\
 \left. +\sum_{j=0}^{p}L^{j+1}(z)\left|\frac{\partial^{p-j} F(z)}{\partial\mathbf{b}^{p-j}}\right|\right)=
 T^{*}|g_{0}(z)|L^{p+1}(z)|\left((T^{*}+1)\times\right.\\ \left. \times\sum_{j=0}^{p}\frac{1}{L^{p-j}(z)}\left|\frac{\partial^{p-j} F(z)}{\partial\mathbf{b}^{p-j}}\right|+\sum_{j=1}^{p}\frac{1}{L^{p+1-j}(z)}\left|\frac{\partial^{p+1-j} F(z)}{\partial\mathbf{b}^{p+1-j}}\right|\right)\leq\\ \leq\!
  T^{*}((T^{*}+1)(p+1)+p)|g_{0}(z)|L^{p+1}(z) %\times\\ \times
  \max\left\{\frac{1}{L^{j}(z)}
 \left|\frac{\partial^{j}F(z)}{\partial\mathbf{b}^{j}}\right|: 0\leq j\leq p\right\}.
\end{gather*}

Thus, for every $r>0$ there exists $P_{3}=P_{3}(r)>0$ that for all $z\in\mathbb{B}_n\setminus
G_{r}^{\mathbf{b}}$ inequality holds
\begin{equation}
\frac{1}{L^{p+1}(z)}\left|\frac{\partial^{p+1}F(z)}{\partial\mathbf{b}^{p+1}}\right|\leq
P_{3}\max\left\{\frac{1}{L^{j}(z)}\left|\frac{\partial^{j}F(z)}{\partial\mathbf{b}^{j}}\right|: 0\leq j\leq p\right\}. \label{eq54}
\end{equation}

Let $z^{0}+t_{0}\mathbf{b}$ be an arbitrary point with $\mathbb{B}_n$  and $$K^{0}\!=\left\{z^{0}+t_{0}\mathbf{b}:
 |t-t_{0}|\leq\frac{\beta}{L(z^{0}+t_{0}\mathbf{b})}\right\}.$$ But $g_{0}, g_{1},\ldots, g_{p}, h$ are analytic in $\mathbb{B}_n$ functions  of bounded $L$-index in  direction $\mathbf{b}.$ Hence, by Theorem \ref{te6} the set $K^{0}$ contains at most $N<+\infty$ elements of the set  $\{c^{0}_{k}\}$ and $N$ is independent of $z^{0}$ and $t_{0}.$

Let $\widetilde{K}^{0}_{k}=\left\{z^{0}+t\mathbf{b}: |t-c^0_{k}|\leq \displaystyle\frac{\mathstrut
\lambda_{1}^{\mathbf{b}}(\beta)(\beta-1)}{8(N+1)L(z^{0}+c^{0}_{k}\mathbf{b})}\right\}.$ 
From condition $L\!\in Q_{\mathbf{b},\!\beta}(\!\mathbb{B}_n)$ it follows $L(z^{0}+c^{0}_{k}\mathbf{b})\geq \lambda_{1}^{\mathbf{b}}(1)L(z^{0}+t_{0}\mathbf{b}).$  If $c^{0}_{k}\in K^{0}$ then $\widetilde{K}_{k}^{0}$ is a subset  $K^{0}_{k}$
$$\widetilde{K}_{k}^{0}\subset K^{0}_{k}=\left\{z^{0}+t\mathbf{b}: |t-c_{k}^{0}|\leq\frac{\beta-1}{8(N+1)L(z^{0}+t_{0}\mathbf{b})}\right\}. $$

From the presented considerations, we deduce that for $z^{0}+t\mathbf{b}\in K^{0}\setminus \bigcup_{c_{k}^0\in K^{0}}K_{k}^{0}$ the inequality (\ref{eq54}) holds with $P_{3}\!=P_{3}\left(\displaystyle\frac{
\lambda_{1}^{\mathbf{b}}(\beta)(\beta-1)}{8(N+1)}\right).$

Again for these $z^{0}+t\mathbf{b}\in K^{0}\setminus \bigcup_{c_{k}^0\in K^{0}}K_{k}^{0}$ inequality holds  $L(z^{0}+t_{0}\mathbf{b})\geq\displaystyle\frac{\mathstrut L(z^{0}+t\mathbf{b})}{\lambda_{2}^{\mathbf{b}}(\beta)}.$ Using (\ref{eq54}),
we have 
\begin{gather} \nonumber
 \frac{1}{L^{p+1}(z^{0}+t_{0}\mathbf{b})}
\left|\frac{\partial^{p+1}F(z^{0}+t\mathbf{b})}{\partial\mathbf{b}^{p+1}}\right|\leq
\lambda_{2}^{\mathbf{b}}(\beta)\frac{1}{L^{p+1}(z^{0}+t\mathbf{b})}\times\nonumber\\ \!\times\!
\left|\frac{\partial^{p+1}F(z^{0}+t\mathbf{b})}{\partial\mathbf{b}^{p+1}}\right|
\leq\! P_{3}(\lambda_{2}^{\mathbf{b}}(\beta))^{p+1}\max\left\{\frac{1}{L^{j}(z^{0}+t\mathbf{b})}
\left|\frac{\partial^{j}F(z^{0}+t\mathbf{b})}{\partial\mathbf{b}^{j}}\right|: \nonumber\right.\\ \left. 0\leq j\leq p\right\}
\leq P_{3}(\lambda_{2}^{\mathbf{b}}(\beta))^{p+1}\max\left\{\frac{1}{L^{j}(z^{0}+t_{0}\mathbf{b})}
\left|\frac{\partial^{j}F(z^{0}+t\mathbf{b})}{\partial\mathbf{b}^{j}}\right| \times\nonumber\right.\\ \times \left.
\left(\frac{1}{\lambda^{\mathbf{b}}_{1}(\beta)}\right)^{j}: 0\leq j\leq p\right\} \leq
P_{3}\left(\frac{\lambda_{2}^{\mathbf{b}}(\beta)}{\lambda_{1}^{\mathbf{b}}(\beta)}\right)^{p}
\lambda_{1}^{\mathbf{b}}(\beta)\max\left\{\frac{1}{L^{j}(z^{0}+t_{0}\mathbf{b})} \times\nonumber\right.\\ \times\left.
\left|\frac{\partial^{j}F(z^{0}+t\mathbf{b})}{\partial\mathbf{b}^{j}}\right|\!: 0\leq\! j\leq\! p\right\}=
P_{4}g_{z^{0}}(t_{0},t),\label{riv_eq_51}
\end{gather}
where $P_{4}=P_{3}\lambda_{2}^{\mathbf{b}}(\beta)\left(\displaystyle\frac{\mathstrut\lambda_{2}^{\mathbf{b}}(\beta)}
{\lambda_{1}^{\mathbf{b}}(\beta)}\right)^{p}$ and
$$g_{z^{0}}(t_{0},t)=\max\left\{\frac{1}{L^{j}(z^{0}+t_{0}\mathbf{b})}\left|\frac{\partial^{j}F(z^0+t\mathbf{b})}{\partial\mathbf{b}^{j
}}\right|: 0\leq j\leq p\right\}. $$ Let $D$ be a sum of diameters of sets $K_{k}^{0}.$ Then
$$D\leq\displaystyle\frac{\mathstrut
2|\mathbf{b}|(\beta-1)N}{8(N+1)L(z^{0}+t_{0}\mathbf{b})}\leq\displaystyle\frac{\mathstrut
|\mathbf{b}|(\beta-1)}{4L(z^{0}+t_{0}\mathbf{b})}.$$ Therefore, 
there exist radii $r_{1}\in\left[\frac{\mathstrut
\beta}{4},\frac{\beta}{2}\right]$ and
$r_{2}\in\left[\displaystyle\frac{\mathstrut
\beta+1}{2},\beta\right]$ with property: if either 
$$z^{0}+t\mathbf{b}\in C_{1}=\left\{z^{0}+t\mathbf{b}:
|t-t_{0}|=\displaystyle\frac{\mathstrut
r_{1}}{L(z^{0}+t_{0}\mathbf{b})}\right\}$$
 $$\text{ or }
z^{0}+t\mathbf{b}\in C_{2}=\left\{z^{0}+t\mathbf{b}:
|t-t_{0}|=\displaystyle\frac{\mathstrut
r_{2}}{L(z^{0}+t_{0}\mathbf{b})}\right\}$$
  then  $z^{0}+t\mathbf{b}\in 
K^{0}\setminus\bigcup_{c_k^0\in K^{0}}K_{m}^{0}.$ We take 
two any points $z^{0}+t_{1}\mathbf{b}\in C_{1}$ and
$z^{0}+t_{2}\mathbf{b}\in C_{2}$ and connect them by a smooth curve
$\gamma=\{z^{0}+t\mathbf{b}: t=t(s), 0\leq s\leq T\}$  that
$F(z^{0}+t(s)\mathbf{b})\neq 0$ and $\gamma\subset
K^{0}\setminus\bigcup_{c_k^0\in K^{0}}K_{m}^{0}.$ This curve
can be selected such that for its length the following estimate holds
\begin{gather} \nonumber
\!|\gamma|\!\leq\! |\mathbf{b}|\left(\frac{\pi r_{1}}
{L(z^{0}+t_{0}\mathbf{b})}+\frac{r_{2}-r_{1}}{L(z^{0}+t_{0}\mathbf{b})}+
\frac{\pi N(\beta-1)}{8(N+1)L(z^{0}+t_{0}\mathbf{b})}\right)\!\leq\\
 \leq |\mathbf{b}|\left(\frac{r_{2}+(\pi-1)r_{1}}{L(z^{0}+t_{0}\mathbf{b})}+
\frac{\pi(\beta-1)}{8L(z^{0}+t_{0}\mathbf{b})}\right)\leq\nonumber\\ \leq |\mathbf{b}|
\frac{1}{L(z^{0}+t_{0}\mathbf{b})}\left(\frac{(\pi-1)\beta}{2}+\beta+\frac{\pi(\beta-1)}{8}\right) < \frac{3\pi\beta|\mathbf{b}|}{L(z^{0}+t_{0}\mathbf{b})}.\label{eq_riv39}
\end{gather}
Then an inequality (\ref{riv_eq_51}) holds on $\gamma$ that is
$$\frac{1}{L^{p+1}(z^{0}+t_{0}\mathbf{b})}\left|\frac{\partial^{p+1} F(z^{0}+t(s)\mathbf{b})}{\partial\mathbf{b}^{p+1}}\right|\leq P_{4} g_{z^0}(t_0,t(s)), \ 0\leq s\leq T.$$

In the proof of Theorem \ref{te11} we obtained that the function $g_{z^{0}}(t_{0},t(s))$ is continuous on $[0,T]$ and continuously differentiable except, perhaps, finite number of points.
 Besides, for complex-valued function of real variable inequality holds
  $\displaystyle\frac{\mathstrut d}{ds}|\varphi(s)|\leq
\left|\displaystyle\frac{d}{ds}\varphi(s)\right|$  except the points, where $\varphi(s)=0.$

Then, in view of \eqref{riv_eq_51}, we have
\begin{gather*}
 \frac{d}{ds}g_{z^{0}}(t_{0},t(s))\leq\max\left\{\!\frac{d}{ds}\frac{1}{L^{j}(z^{0}+t_{0}\mathbf{b}))}
\left|\frac{\partial^{j}F(z^{0}+t(s)\mathbf{b})}{\partial\mathbf{b}^{j}}\right|:  0\leq j\leq p\right\}\leq \\
\leq\max\left\{\frac{1}{L^{j+1}(z^{0}+t_{0}\mathbf{b})}\left|\frac{\partial^{j+1}F(z^{0}+t(s)\mathbf{b})}
{\partial\mathbf{b}^{j+1}}\right| |t'(s)| L(z^{0}+t_{0}\mathbf{b}): 0\leq j\leq p\right\}\leq \\
\leq\max\left\{\frac{1}{L^{j+1}(z^{0}+t_{0}\mathbf{b})}\left|\frac{\partial^{j+1}F(z^{0}+t(s)\mathbf{b})}
{\partial\mathbf{b}^{j+1}}\right|: 0\!\leq\! j\!\leq\! p; \!\left|\!\frac{\partial^{p+1}F(z^{0}+t(s)\mathbf{b})}
{\partial\mathbf{b}^{p+1}}\right|\right\}\!\times\!\\ \times|t'(s)|L(z^{0}+t_{0}\mathbf{b})\leq
P_{5}g_{z^{0}}(t_{0},t(s))|t'(s)|L(z^{0}+t_{0}\mathbf{b}).
\end{gather*}
where $P_5=\max\{1,P_4\}.$
But (\ref{eq_riv39}) is true, then
\begin{gather*}
 \left|\ln\frac{g_{z^{0}}(t_{0},t_{2})}{g_{z^{0}}(t_{0},t_{1})}\right|=
\left|\int_{0}^{T}\frac{1}{g_{z^{0}}(t_{0},t(s))}\frac{d}{ds}g_{z^{0}}(t_{0},t(s))ds\right|\leq\\
\leq P_{5}L(z^{0}+t_{0}\mathbf{b})\int_{0}^{T}|t'(s)|ds\leq P_{5}L(z^{0}+t_{0}\mathbf{b})|\gamma|\leq
3\pi\beta|\mathbf{b}|P_{5},
\end{gather*}
i.e. $$g_{z^{0}}(t_{0},t_{2})\leq g_{z^{0}}(t_{0},t_{1})\exp\{3\pi\beta|\mathbf{b}|P_{5}\}. $$
We can choose $t_{2}$  that
$|F(z^{0}+t_{2}\mathbf{b})|=\max\{|F(z^{0}+t\mathbf{b})|: z^{0}+t\mathbf{b}\in C_{2}\}.$ Hence,
\begin{gather} \nonumber
 \max\left\{|F(z^{0}+t\mathbf{b})|: |t-t_{0}|=\frac{\beta+1}{2L(z^{0}+t_{0}\mathbf{b})}\right\}\leq
|F(z^{0}+t_{2}\mathbf{b})|\leq \\ \leq g_{z^{0}}(t_{0},t_{2}) \leq g_{z^{0}}(t_{0},t_{1}) \exp\{3\pi\beta|\mathbf{b}|P_{5}\}.\label{eq_riv_55}
\end{gather}

Since $z^{0}+t_{1}\mathbf{b}\in C_{1}$ we apply Cauchy inequality in variable $t$ for all $j=1,2,\ldots,p,$ and obtain
 \begin{gather*}
  \left|\frac{\partial^{j} F(z^{0}+t_{1}\mathbf{b})}{\partial\mathbf{b}^{j}}\right|
%   =  \frac{j!}{2\pi}
%  \left|\int\limits_{\ |t-t_{1}|=1/(10L(z^{0}+t_{0}\mathbf{b}))}\frac{F(z^{0}+t\mathbf{b})}
%  {(t-t_{1})^{j+1}}dt\right|
%  \leq\\
%  \leq \\
 \leq j!\left(10L(z^{0}+t_{0}\mathbf{b})\right)^{j}\max\left\{|F(z^{0}+t\mathbf{b})|:
|t-t_{1}|=\frac{1}{2\beta L(z^{0}+t_{0}\mathbf{b})}\right\}\!\leq\!\\
 \leq p!\left(10L(z^{0}+t_{0}\mathbf{b})\right)^{j}
  \max\left\{|F(z^{0}+t\mathbf{b})|:
|t-t_{0}|=\frac{1}{\beta L(z^{0}+t_{0}\mathbf{b})}\right\}
 \end{gather*}
It follows that 
\begin{equation}
 g_{z^{0}}(t_{0},t_{1})\leq p!10^{p}
\max\left\{|F(z^{0}+t\mathbf{b})|:
|t-t_{0}|=\frac{1}{\beta L(z^{0}+t_{0}\mathbf{b})}\right\} \label{eq_riv_56}
\end{equation}
 From inequalities (\ref{eq_riv_55}) and (\ref{eq_riv_56}) we have 
\begin{gather*}
 \!\max\left\{|F(z^{0}+t\mathbf{b})|: |t-t_{0}|=\frac{\beta+1}{2L(z^{0}+t_{0}\mathbf{b})}\right\}\leq
p!10^{p}\exp\{|\mathbf{b}|P_{5}\}\!\times\!\\
\times\max\left\{|F(z^{0}+t\mathbf{b})|:
|t-t_{0}|=\frac{1}{\beta L(z^{0}+t_{0}\mathbf{b})}\right\}.
\end{gather*}
 Therefore, by Theorem \ref{te4} an analytic function $F(z)$ is of bounded $L$-index in  direction $\mathbf{b}.$
\end{proof}

% \begin{center}
\section{Growth of analytic in $\mathbb{B}_n$ functions of bounded $L$-index in the direction.} 
% \end{center}  

% \textbf{$8^0.$ Growth of analytic in $\mathbb{B}_n$ functions of bounded $L$-index in the direction.}
We denote $a^+=\max\{a,0\}.$
\begin{theorem}
 Let $L: \mathbb{B}_n\to\mathbb{R}_{+},$ for every $z^0\in\mathbb{B}_n,$ $\theta\in[0,2\pi]$ a function $L(z^0+re^{i\theta}\mathbf{b})$ be a continuously differentiable function of real variable
 $r\in[0,R),$ where $R=\min\{t\in\mathbb{R}_+: |z^0+te^{i\theta}\mathbf{b}|=1\}.$
   If an analytic in $\mathbb{B}_n$ function $F$ is of bounded $L$-index in  direction $\mathbf{b}$ then
for every $z^0\in\mathbb{B}_n,$ $\theta\in[0,2\pi],$ $r\in[0,R)$ and every integer $p\geq 0$
\begin{gather}
\ln\left(\frac{1}{p!L^p(z^0+re^{i\theta}\mathbf{b})}\left|\frac{\partial^p F(z^0+re^{i\theta}\mathbf{b})}{\partial\mathbf{b}^p}\right|\right)\!\leq\! %\nonumber \\ \leq
\ln \max\left\{\frac{1}{k!L^k(z^0)}\left|\frac{\partial^k F(z^0)}{\partial\mathbf{b}^k}\right|: 0\leq k\leq N\right\}\!+\!\nonumber \\+
\int_0^r \left\{(N+1)L(z^0+te^{i\theta}\mathbf{b})+N\frac{(-L'_t(z^0+te^{i\theta}\mathbf{b}))^+}{L(z^0+te^{i\theta}\mathbf{b})} \right\} dt \label{eqim1}
\end{gather}
But if, in addition, for every $z^0\in\mathbb{B}_n$ and $\theta\in[0,2\pi]$
$\left(-\frac{\partial L(z^0+re^{i\theta}\mathbf{b})}{\partial\mathbf{b}}\right)^+/(L^2(z^0+re^{i\theta}\mathbf{b}))\rightrightarrows 0$  as  $|z^0+re^{i\theta}\mathbf{b}|\to 1$  then
for every $z^0\in\mathbb{B}_n$ and $\theta\in[0,2\pi]$
\begin{equation}
\label{eq344}
\varlimsup_{|z^0+re^{i\theta}\mathbf{b}|\to 1} \frac{\ln |F(z^0+re^{i\theta}\mathbf{b})|}{\int_0^r L(z^0+te^{i\theta}\mathbf{b})dt}\leq N_{\mathbf{b}}(F,L)+1,
\end{equation}
holds.
\end{theorem}

\begin{proof}

 We remark that $R\geq \frac{1-|z^0|}{|\mathbf{b}|},$ because
 $|z^0+te^{i\theta}\mathbf{b}|\leq |z^0|+|t|\cdot|\mathbf{b}|\leq |z^0|+\frac{1-|z^0|}{|\mathbf{b}|}\cdot |\mathbf{b}|\leq 1.$ The condition
 $r\in [0,R)$ provides $z^0+re^{i\theta}\mathbf{b}\in\mathbb{B}_n.$

 Denote $N=N_{\mathbf{b}}(F,L).$ For fixed $z^0\in\mathbb{B}_n$ and $\theta\in[0,2\pi]$ we consider the function
 \begin{equation}
 g(r)=\max\left\{\frac{1}{k!L^{k}(z^0+re^{i\theta}\mathbf{b})}\left|\frac{\partial^k F(z^0+re^{i\theta}\mathbf{b})}{\partial\mathbf{b}^k}\right|: 0\leq k\leq N\right\}. \label{defgr}
\end{equation}

 Since the function $\frac{1}{k!L^{k}(z^0+re^{i\theta}\mathbf{b})}\left|\frac{\partial^k F(z^0+re^{i\theta}\mathbf{b})}{\partial\mathbf{b}^k}\right|$ is a continuously differentiable of real $r\in[0,R),$ the function $g$ is continuously differentiable on
 $[0,R),$ exception, perhaps, a finite set of points, and
 \begin{gather*}
  g'(r)\leq \max\left\{\frac{d}{dt}\left(\frac{1}{k!L^{k}(z^0+re^{i\theta}\mathbf{b})}\left|\frac{\partial^k F(z^0+re^{i\theta}\mathbf{b})}{\partial\mathbf{b}^k}\right|\right): 0\leq k\leq N\right\}\leq \\
  \leq \max\left\{\frac{1}{k!L^k(z^0+re^{i\theta}\mathbf{b})}\left|\frac{\partial^{k+1} F(z^0+re^{i\theta}\mathbf{b})}{\partial\mathbf{b}^{k+1}}\right|-\frac{1}{k!L^k(z^0+re^{i\theta}\mathbf{b})}\times\right.\\ \left.\times\left|\frac{\partial^k F(z^0+re^{i\theta}\mathbf{b})}{\partial\mathbf{b}^k}\right|k 
  \frac{L'_r(z^0+re^{i\theta}\mathbf{b})}{L(z^0+re^{i\theta}\mathbf{b})}: \ 0\leq k\leq N\right\} \leq \\ 
  \leq   \max\left\{\frac{1}{(k+1)!L^{k+1}(z^0+re^{i\theta}\mathbf{b})}\left|\frac{\partial^{k+1} F(z^0+re^{i\theta}\mathbf{b})}{\partial\mathbf{b}^{k+1}}\right| (k+1)L(z^0+re^{i\theta}\mathbf{b})+\right.
  \\ \left. +\frac{1}{k!L^k(z^0+re^{i\theta}\mathbf{b})}\left|\frac{\partial^k F(z^0+re^{i\theta}\mathbf{b})}{\partial\mathbf{b}^k}\right|k\frac{(-L'_r(z^0+re^{i\theta}\mathbf{b}))^+}{L(z^0+re^{i\theta}\mathbf{b})}:
  \right. \\ \left. 0\leq k\leq N\right\}\leq g(r)\left((N+1)L(z^0+re^{i\theta}\mathbf{b})+N\frac{(-L'_r(z^0+re^{i\theta}\mathbf{b}))^+}{L(z^0+re^{i\theta}\mathbf{b})}\right).
 \end{gather*}
 Thus, we have
 $$
 \frac{d}{dr} \ln g(r)\leq (N+1)L(z^0+re^{i\theta}\mathbf{b})+N\frac{(-L'_r(z^0+re^{i\theta}\mathbf{b}))^+}{L(z^0+re^{i\theta}\mathbf{b})}.
 $$
 Since $F$ is a function of bounded $L$-index in the direction then $g(0)\neq 0$ and
  \begin{equation*}
  g(r)\!\leq\! g(0)\exp\left\{
\int_0^r \left(  (N+1)L(z^0+te^{i\theta}\mathbf{b})\!+\!N\frac{(-L'_t(z^0\!+\!te^{i\theta}\mathbf{b}))^+}{L(z^0\!+\!te^{i\theta}\mathbf{b})}\right)dt
  \right\}, \ r\to\! R,
 \end{equation*}
 so that
 $$
 \ln g(r)\leq \ln g(0) \!+\! \int_0^r \left((N\!+\!1)L(z^0+te^{i\theta}\mathbf{b})+N\frac{(-L'_t(z^0+te^{i\theta}\mathbf{b}))^+}{L(z^0+te^{i\theta}\mathbf{b})}\right)dt, \ r\to R.
$$
 Using \eqref{defgr}, we obtain \eqref{eqim1}.
If, in addition,  for every $z^0\!\in\!\mathbb{B}_n$ and $\theta\!\in\![0,2\pi\!]$ 
$\!\left(\!-\frac{\partial L(z^0+re^{i\theta}\mathbf{b})}{\partial\mathbf{b}}\!\right)^+/(L^2(z^0+re^{i\theta}\mathbf{b}))\rightrightarrows 0$  when  $|z^0+re^{i\theta}\mathbf{b}|\to 1$  then
 \begin{gather*}
  g(r)\leq g(0)\exp\left\{
  (N+1)\int_0^r \left(L(z^0+te^{i\theta}\mathbf{b})+\frac{(-L'_r(z^0+te^{i\theta}\mathbf{b}))^+}{L(z^0+te^{i\theta}\mathbf{b})}\right)dt
  \right\}=\\ =g(0)\exp\left\{(N+1)(1+o(1))\int_0^rL(z^0+te^{i\theta}\mathbf{b})dt\right\}, \ r\to R,
 \end{gather*}
so that
$$
|F(z^0\!+\!re^{i\theta}\mathbf{b})|\!\leq\! g(r)\leq g(0)\exp\left\{(N+1)(1+o(1))\int_0^rL(z^0+te^{i\theta}\mathbf{b})dt\right\}, \ r\to R,
$$
for $\theta\in[0,2\pi],$ $z^0\in\mathbb{B}_n,$ whence
\begin{equation} \label{enheq}
\ln |F(z^0+re^{i\theta}\mathbf{b})|\leq g(0)+(N+1)(1+o(1))\int_0^rL(z^0+te^{i\theta}\mathbf{b})dt, \ r\to R.
\end{equation}
Moreover, for every $z^0\in\mathbb{B}_n$ and $\theta\in[0,2\pi]$ we have 
\begin{equation*}
\varlimsup_{|z^0+re^{i\theta}\mathbf{b}|\to 1} \frac{\ln |F(z^0+re^{i\theta}\mathbf{b})|}{\int_0^r L(z^0+te^{i\theta}\mathbf{b})dt}\leq N_{\mathbf{b}}(F,L)+1.
\end{equation*}
\end{proof}

\begin{remark}
The equations \eqref{eqim1} and \eqref{eq344} can be written in more convenient forms:
\begin{gather}
\ln\max\limits_{|t|=r}\left(\frac{1}{p!L^p(z^0+t\mathbf{b})}\left|\frac{\partial^p F(z^0+t\mathbf{b})}{\partial\mathbf{b}^p}\right|\right)\leq%\nonumber\\ \leq 
\ln\max\left\{\frac{1}{k!L^k(z^0)}\left|\frac{\partial^k F(z^0)}{\partial\mathbf{b}^k}\right|: 0\leq k\leq N\right\}+\nonumber \\+
\max\limits_{\theta\in[0,2\pi]} \int_0^r \left\{(N+1)L(z^0+te^{i\theta}\mathbf{b})+N\frac{(-L'_t(z^0+te^{i\theta}\mathbf{b}))^+}{L(z^0+te^{i\theta}\mathbf{b})} \right\} dt \label{eqim1c}
\end{gather}
and
\begin{equation}
\label{eq344c}
\varlimsup_{|z^0+re^{i\theta}\mathbf{b}|\to 1} \max\limits_{\theta\in[0,2\pi]} \frac{\ln |F(z^0+re^{i\theta}\mathbf{b})|}{\int_0^r L(z^0+te^{i\theta}\mathbf{b})dt}\leq N_{\mathbf{b}}(F,L)+1.
\end{equation}
Besides, if we put $z^0=0$ then the estimate \eqref{enheq} implies a following inequality
\begin{equation}
 \label{enhlim}
 \varlimsup_{R\to1/|\mathbf{b}|}  \frac{\ln \max\{|F(t\mathbf{b})|: \ |t|=R\}}{\max\limits_{\theta\in[0,2\pi]}\int_0^R L(re^{i\theta}\mathbf{b})dr}\leq N_{\mathbf{b}}(F,L)+1.
\end{equation}
\end{remark}

For $n=1$ we deduce corollaries.
\begin{corollary}\label{improvcor}
Let $l: \mathbb{D}\to \mathbb{R}_+,$ $\mathbb{D}=\{z\in\mathbb{C}: |z|<1\}$ and for $\theta\in[0,2\pi]$ a function $l(re^{i\theta})$ be a  continuously differentiable function of real variable $t\in[0,1).$ If $f(z)$ is an analytic function of bounded $l$-index then for every integer $p\geq 0$
\begin{gather}
\ln\frac{|f^{(p)}(re^{i\theta})|}{p!l^p(re^{i\theta})}\leq\nonumber \\ \leq 
\ln\max\left\{\frac{|f^{(k)}(0)|}{k!l^k(0)}: 0\leq k\leq N\right\}+
\int_0^r  \left\{(N+1)l(te^{i\theta})+N\frac{(-L'_t(te^{i\theta}))^+}{L(te^{i\theta})} \right\} dt \label{croeq2}
\end{gather}
If, in addition,
$(-l'(re^{i\theta}))^+/l^2(re^{i\theta})\rightrightarrows 0$ as $r\to 1$ then
\begin{equation} \label{croeq1}
\varlimsup_{r\to 1} \frac{\ln |f(re^{i\theta})|}{\int_0^r l(te^{i\theta})dt}\leq N(f,l)+1, \ \theta\in[0,2\pi]
\end{equation}
holds, where $N(f,l)$ is $l$-index of function $f.$
\end{corollary}

\begin{remark}
The equations \eqref{croeq2} and \eqref{croeq1} can be written in more convenient forms
\begin{gather}
\!\ln\max_{|t|=r}\!\frac{|f^{(p)}(t)|}{p!l^p(t)}\!\leq\! \nonumber\\ \leq 
\ln\max\left\{\frac{|f^{(p)}(0)|}{p!l^p(0)}: 0\!\leq\! k\!\leq\! N\!\right\}\!+\!\max\limits_{\theta\in[0,2\pi]}
\!\int_0^r\!  \left\{\!(\!N\!+\!1\!)l(te^{i\theta})\!+\!N\frac{(\!-L'_t(te^{i\theta}))^+}{L(te^{i\theta})} \right\} dt
\end{gather}
and
\begin{equation}
\varlimsup_{r\to 1} \max_{\theta\in[0,2\pi]} \frac{\ln |f(re^{i\theta})|}{\int_0^r l(te^{i\theta})dt}\leq N(f,l)+1, \
\end{equation}
\end{remark}

The Corollary \ref{improvcor} is an improvement of similar result of Sheremeta and Strochyk \cite{strosher} because
we do not assume that $l=l(|z|).$

\begin{corollary}
 Let $F: \mathbb{B}_n\to \mathbb{C}$ be an analytic function of bounded $L$-index in direction $\mathbf{b},$ $N=N_{\mathbf{b}}(F,L),$
  $z^0$ be a fixed point in $\mathbb{B}_n$ that $F(z^0)=1.$ Then for every $r\in[0,R),$ where $R=\min\{t\in\mathbb{R}_+: |z^0+te^{i\theta}\mathbf{b}|=1\},$ the next inequality
  \begin{gather*}
  \int_0^r \frac{n(t,z^0,0,1/F)}{t}dt\leq \ln\max\{|F(z^0+t\mathbf{b})|: \ |t|=r\}\leq\\ \leq 
  \ln\max\left\{\frac{1}{p!L^p(z^0)}\left|\frac{\partial^p F(z^0)}{\partial\mathbf{b}^p}\right|:  0\leq k\leq N\right\}+\\ +
\max\limits_{\theta\in[0,2\pi]} \int_0^r \left\{(N+1)L(z^0+te^{i\theta}\mathbf{b})+N\frac{(-L'_t(z^0+te^{i\theta}\mathbf{b}))^+}{L(z^0+te^{i\theta}\mathbf{b})} \right\} dt
  \end{gather*}
  holds.
\end{corollary}
\begin{proof}
 We consider a function $F(z^0+t\mathbf{b})$ as a function of one variable $t.$ Thus, the first inequality follows from the classical Jensen Theorem. In addition, the second inequality follows from \eqref{eqim1c} for $p=0.$
\end{proof}

\end{document}